\documentclass[12pt]{article}

\usepackage{amsfonts,amssymb,amstext}
\usepackage{subcaption}
\usepackage{dsfont}
\usepackage{epsfig}
\usepackage{graphicx}
\usepackage{enumerate}
\usepackage{empheq}
\usepackage{comment}
\usepackage{color}
\usepackage{hyperref}
\hypersetup{colorlinks=true,citecolor=blue,urlcolor=blue,linkcolor=blue}
\usepackage{cite} 
\usepackage{amsmath}
\allowdisplaybreaks

\newcommand{\fr}{Fr}

\usepackage[titletoc, title]{appendix}

\newcommand{\beginsupplement}{%
        \setcounter{table}{0}
        \renewcommand{\thetable}{S\arabic{table}}%
        \setcounter{figure}{0}
        \renewcommand{\thefigure}{S\arabic{figure}}%
        \setcounter{section}{0}
        \renewcommand{\thesection}{S.\arabic{section}}%
        \setcounter{equation}{0}
        \renewcommand{\theequation}{S\arabic{equation}}%
     }

\makeatletter
  \newcounter{savesection}
\newcounter{apdxsection}
\renewcommand\appendix{\par
  \setcounter{savesection}{\value{section}}%
  \setcounter{section}{\value{apdxsection}}%
  \renewcommand\theHsection{appendix.\thesection}
  \setcounter{subsection}{0}%
  \gdef\thesection{\@Alph\c@section}}
  
\newcommand\unappendix{\par
  \setcounter{apdxsection}{\value{section}}%
  \setcounter{section}{\value{savesection}}%
  \renewcommand\theHsection{\thesection}
  \setcounter{subsection}{0}%
  \gdef\thesection{\thechapter.\arabic{section}}}
\makeatother

\title{A transient depth-averaged lava flow model with a Herschel–Bulkley rheology accounting for three phases}

\author{J. Binard$^1$, A. Burgisser$^2$, E.D. Fernández-Nieto$^1$, G. Narbona-Reina$^1$}

\begin{document} 
\date{}
\maketitle

{\scriptsize $^1$ Dpto. Matemática Aplicada I E.T.S Arquitectura, Universidad de Sevilla, Sevilla, Spain}

{\scriptsize $^2$ Univ. Grenoble Alpes, Univ. Savoie Mont Blanc, CNRS, IRD, Univ. Gustave Eiffel, ISTerre, Grenoble, France}

\bigskip

\begin{abstract}
    This study presents a three-phase suspension lava flow model with a Herschel–Bulkley rheology. The suspension contains crystals and gas bubbles, and two closures for the evolution of the crystal volume fraction are considered. The first closure minimizes the complexity of the system by treating crystal fraction as a transported quantity subject to relaxation towards an equilibrium state. This closure avoids a parametrization of the many heat transfer mechanisms a lava flow is subjected to. The other closure is on the lava temperature considering four heat transfer mechanisms and a prescribed temperature--crystallinity relationship. We deduce from this system and solve numerically a one-dimensional depth-averaged model. A comparison with a pre-existing model based on real lava flow data suggests that the prediction of flow parameters done with the multi-parametric evolution of temperature yields more accurate results than those obtained with the other closure. The transient nature of our model correctly predicts that confined lava traveling down an irregular steep slope yields a series of cascading fill-then-breakout lumps that causes the overall flow to be pulsatory. These pulses dominate the local dynamics and preclude a strict steady state to be reached. Theoretically, taking gas bubbles into account is best done with a general rheological relationship valid at any capillary number. In the conditions explored herein, bubbles modulate viscosity within a factor 2 with a shear thinning behavior, decelerating slow flows and accelerating fast flows. When a simplified rheology treating bubbles as hard spheres was used, the only dynamic parameter affected was bulk viscosity.\\

\end{abstract}

\paragraph{Keywords}
lava flow ; 
three-phase ;
bubble suspension ;
rheology ;
crystallization ;
Mauna Ulu ;

\newpage

\tableofcontents

\section{Introduction}

Lava flows have surprising rich dynamics because of large variations of rheology; the magma is emitted at high temperature at the vent system, and the cooling during flow causes some of it to become solid. Such crossing of rheological boundaries causes some of the lava to solidify as small blocks that accumulate on the flow sides and front, creating levees that channel flows. Sometimes, theses levees modify flow path by fostering braiding. Such creation of diverging flow branches increases significantly the area that will be covered by the lava. Sometimes levees rise high enough to cause a full encasing of the flow within its channel, forming a lava tube. Such channeling considerably lengthen the reach of the downstream flow because the channel limits lava cooling \cite{harris2022anatomy}.\\

Prediction of the emplacement of lava flows is required for hazard and risk assessment, and for the planning of risk-mitigation measures. The main physical parameters controlling lava emplacement are topography (or slope in 1D), eruptive input conditions (volume effusion rate, source geometry (vent or fissure) and effusive temperature), thermal boundary conditions atop and at the base of the flow (insulation, convection, radiation and conduction), and the lava physical properties (density, thermal conductivity, rheology). Numerical models have various capabilities of approximating the motion of a given volume of lava from its source vent to the deposition area. Such models have undergone major developments over the past years  \cite{scott_2026}, and they can be separated into the two broad categories of probabilistic and deterministic approaches \cite{cordonnier:2016}. Because of their fast computational time compared to more complex deterministic models, numerical simulations of lava flows based on the probability of lava presence at each discretized location played a role in the monitoring and mitigation of three recent eruptions  \cite{cappello_2016, carracedo_2022, pedersen_2023}.\\

Deterministic models, on the other hand, are intended to mimic the behavior of the natural systems by modeling physical processes based on a set of conservation equations and including parameters such as lava discharge rate, lava temperature and heat transfer. Models including more physical processes will take longer to produce a solution. At the low end-member of complexity, one finds 1D deterministic models, which are restricted to the simulation of channeled flows, such as FLOWGO and PyFLOWGO \cite{chevrel:2018, harris2001flowgo}, which only consider temperature evolution. A notable example is the model of \cite{crisp1994influence}, which describes the influence of crystallization on the evolution of lava flow temperature. Both the temperature and the solid concentration are assumed to depend only on time and their joint evolution is combined into a single differential equation. One difficulty when modeling temperature evolution is the many mechanisms of heat transfer that this single temperature equation has to take into account. A remarkably exhaustive list of these mechanisms has been proposed by \cite{patrick_2004}. Heat transfer at the lava surface is affected by radiation, natural and forced convection, whereas within the flow conduction, convection, and latent heat release operate. Most lava flow models differ by their selection of these mechanisms. For instance, \cite{crisp1994influence} included four mechanisms and \cite{harris2001flowgo} considered seven mechanisms.\\

A level of complexity above, shallow-water models offer an interesting trade-off between calculation speed and physical completeness and nicely bridge the gap between probabilistic models and a complete 3D resolution of a deterministic model. First introduced by De Saint-Venant in 1864, these depth-averaged  models are used in many applications involving geophysical flows under the assumption of shallow flow. Their main limitation is the loss of vertical description in the involved variables. In the case of lava flows, which have a strong feedback between temperature, viscosity, and flow dynamics \cite{zeinalova:2025}, this restriction limits the influence of rheology on flow behaviour.\\

The most recent depth-averaged models of lava flows can be found in \cite{costa_macedonio_2005, kelfoun:2015, bernabeu:2013, biagioli:2023}. These studies propose two-dimensional models for non-Newtonian fluids, treating lava as a single phase. Temperature evolution is considered by three models \cite{costa_macedonio_2005, bernabeu:2013, biagioli:2023}, but crystallization and lava composition are not explicitly accounted for; they are implicitly considered with temperature-dependent viscosity relationships. In \cite{biagioli:2023}, the lack of vertical viscosity variation is addressed by considering a vertical temperature profile that results in a vertical viscosity profile. Although extending such models to explicitly treat solids may increase computational cost, it would allow one to explore processes such as disequilibrium crystallization. Another limitation of the single-phase approach is that the presence of gas is neglected. Lava flows are often emitted at the vent with a bubble volume fraction of $>70$ vol\% that then decreases drastically downflow, significantly increasing bulk flow viscosity \cite{soldati:2020, bokharaeian:2023}.\\

Here we present a multiphase model for lava flows that explicitly accounts for melt, crystals, and gas bubbles. The lava flow is simulated with a depth-averaged one-dimensional model with Herschel-Bulkley rheology. Constitutive laws describing the rheology of two- and three-phase suspensions are notoriously difficult to establish and measure in the laboratory. Our rheology is built on the constitutive relationships proposed by \cite{truby:2015} for three-phase flows and by \cite{mueller:2010} for a particle-fluid mixture. We propose a new closure for the flow index $n$ that match previous experimental results. We also consider two possible closures for the evolution of the crystal volume fraction. The first is an attempt to minimize the complexity of the system by treating crystal fraction as a transported quantity subject to a simple relaxation towards an equilibrium state. This closure is based on the characteristic time of crystallization and thus is not linked to the temperature evolution, thereby avoiding a parametrization of the many heat sources and sinks a lava flow is subjected to. The other closure is on the multi-parametric evolution of the lava temperature considering four heat transfer mechanisms with a prescribed temperature--crystallinity relationship and thus only simulates equilibrium crystallization. We present several numerical verification results and an inter-model comparison based on the case of a channeled flow emitted during the Mauna Ulu eruption.\\

\begin{table}
    \centering
    \begin{tabular}{c|c}
    $\eta$     &  Apparent viscosity\\
    $\eta_0$     & Fluid viscosity\\
    $n$ & Flow index\\
    $\tau_c$ & Critical shear stress\\
    $\tau_c^{\max}$ & Maximum shear stress\\
    $\Phi_b, \Phi_b^*$ & Bubbles volume fraction (crystal-free, bulk)\\
    $\Phi_p, \Phi_p^*$ & Crystal volume fraction (bulk, bubble-free) \\
    $\Phi_m^*$ & Maximum packing particle fraction (bubble-free)\\
    $\Phi_p^c$ & Critical particle fraction\\
    $r_p$ & Relative size of particles\\
    $\tau_{\rm crys}$ & Crystallization time\\
    \end{tabular}
    \caption{Notation}
    \label{tab:notation}
\end{table}

\section{The rheology of lava flows in the literature}\label{SectRheol}

Magma is a three-phase mixture composed of silicate melt, crystals and bubbles. The melt is considered as Newtonian fluid and the suspended bubbles or crystals modify the rheology of the mixture. As a result, lava rheology depends on its composition (melt chemical composition, crystals and bubbles concentrations) and temperature (crystallization and vesiculation). The resulting bulk rheology is non-Newtonian, which includes its dependence on strain-rate and yield stress. One of the most commonly used model is the Herschel-Bulkley rheology, in which the stress tensor is defined as
\begin{equation}\label{HB0}
\sigma = \left(\tau_c+ \eta |\dot\gamma|^{n}\right)\frac{\dot\gamma}{|\dot\gamma|}
\end{equation}
where $\dot\gamma$ is the shear rate, $\eta$ is the apparent viscosity of the mixture, $n$ is the flow index that indicates the degree of non-Newtonian behaviour and $\tau_c$ is the critical yield stress below which the velocity is constant or zero. The Bingham rheology is obtained when $n=1$. The Newtonian case is recovered when $\tau_c=0$ and $n=1$, in which case $\eta$ becomes the Newtonian viscosity. It is common in the literature to call the apparent viscosity $\eta$ a consistency factor, because when $n\neq 1$ the dimension of $\eta$ in \eqref{HB0} is not Pa.s (a viscosity sensu stricto) but Pa.s$^n$. Here we mostly use apparent viscosity instead of consistency factor to be consistent with the nomenclature of the original works from which stem the relationships, and often use simply viscosity as a shorthand.\\

The three parameters, $\eta, n, \tau_c$, depend on the composition of the lava and need constitutive relationships. Such relationships are scarce in the literature for three-phase flows, which is partly due to the difficulty of obtaining reliable and reproducible rheological measures in laboratory experiments. The most common analysis of these experimental results is semi-empirical, combining a theoretical approach of the mixture mechanical properties with some fitting parameters. To some extent, the approximations of such fitting is comparable to those introduced by the simplifications considered in the theoretical approximation and the uncertainty in laboratory measurements. In general, the uncertainty in rheological parameters increases at higher concentrations.\\

The {\bf presence of solids} influences the rheology of magma in complex ways, causing yield stress and strain rate dependence \cite{castruccio:2010, mueller:2010, mader:2013}. Unlike bubbles, particles causes the apparent viscosity to always increases with concentration.  The aspect ratio of the particles (maximum length over equivalent diameter), which we denote $r_p$, also controls the rheology. Both particle size and aspect ratio influence the maximum packing value $\Phi_m^*$. Experiments show that both viscosity and yield stress are independent of $r_p$ if they are expressed as a function of relative concentration, $\Phi_p/\Phi_m^*$. In reality, the dependence on $r_p$ is implicit when setting a specific value of $\Phi_m^*$. However, the flow index $n$ (indicating shear thinning behavior when $n<1$, and shear thickening when $n>1$) explicitly depends on $r_p$ \cite{mueller:2010}.
Most of the models in the literature are valid for monodisperse suspensions with $\Phi_p<\Phi_m^*$. Beyond this value, the solid suspension is subjected to plastic deformation, a regime not considered in this work.\\

The {\bf influence of bubbles} on lava rheology differs as a function of their concentration and of the shear conditions. Bubbles mostly influence the apparent viscosity of the suspension \cite{llewellin2005bubble, mader:2013,bokharaeian:2023}. They can also produce shear thinning, thereby modifying the flow index $n$  \cite{stein:1992,manga:1998, truby:2015, birnbaum:2021, soldati:2020}. Their influence on the yield stress only happens at high bubbles concentrations \cite{princen:1989, rouyer:2005, cohen:2014,namiki:2017,namiki:2022}. \\

The bulk viscosity is affected in different ways depending on the shape of bubbles, an effect captured by the dimensionless capillary number for steady (Ca) or unsteady (Cd) flows,
\begin{equation}
    {\rm Ca}=\lambda \dot\gamma;\qquad {\rm Cd}=\lambda\frac{\dot\gamma}{\ddot\gamma},
\end{equation}
where $\lambda=\frac{\eta_0 a}{\Gamma}$ is the bubble relaxation time\footnote{the timescale over which a bubble can respond to changes in its shear environment}. The capillary number Ca is the ratio of the viscous force (that deforms bubbles) to the surface tension force (that tends to keep bubbles spherical).
In either case (Ca or Cd), if the capillary number is small, the viscosity increases with bubble concentration because the bubbles remain spherical, hindering flow. Conversely, if the capillary number is large, then viscosity decreases as concentration increases because bubbles deform into an elongated shape that does not hinder flow. In the limit case ${\rm Ca}= 1$ or ${\rm Cd}= 1$, there is not clear relationship for the viscosity behaviour, which is often treated by simply matching the two functions for ${\rm Ca}<1$ and ${\rm Ca}>1$.\\

As for solids, most of the bubbly flow models proposed in the literature are valid for monodisperse suspensions with concentration $\Phi_b\lessapprox 0.5$. At higher concentrations, the system is in a foamy state, which has a different rheology because surface tension plays an important role \cite{namiki:2017,namiki:2022}.\\

The {\bf influence of the temperature} on lava flows is difficult to neglect. Temperature changes over space and time and strongly influences their rheology. During cooling, the temperature can change from that at the vent to the solidus temperature when the complete crystallization is achieved. In general, the temperature evolution of a lava flow over space and time is described by an advection equation. It controls crystallisation, which affects the bulk viscosity by changing the solid volume fraction but it is never easy to establish constitutive equations describing crystallization  \cite{marsh:1981,dragoni:1994,mader:2013}.\\

Depending on the phase considered, there are two ways to define the volume fraction of bubbles, $\Phi_b$, and crystals, $\Phi_p$.
The first definition is used when particles are progressively added to a suspension containing bubbles (i.e. crystallization in a bubble-bearing magma), 
\begin{align}\label{phibp1}
	\Phi_b = \frac{V_b}{V_l+V_b}, \qquad \Phi_p = \frac{V_p}{V_l+V_b+V_p}.
\end{align}
The second definition is for the case when bubbles are progressively added to a suspension with particles (i.e. vesiculation in an crystal-bearing magma), 
\begin{align}\label{phibp2}
\Phi_b^* = \frac{V_b}{V_l+V_b+V_p}, \qquad \Phi_p^* = \frac{V_p}{V_l+V_p}.
\end{align} 
These two definitions are related by the following formulas:
\begin{equation}\label{phibp_relation}
\begin{array}{lcl}
& &\displaystyle	\Phi_b = \frac{\Phi_b^*}{1-\Phi_p^*(1-\Phi_b^*)}, \qquad \Phi_p = \Phi_p^*(1-\Phi_b^*);
\\
\text{or inversely,} & &
\\
& &\Phi_b^*=\Phi_b(1-\Phi_p), \qquad \displaystyle \Phi_p^*=\frac{\Phi_p}{1-\Phi_b(1-\Phi_p)}.
\end{array}
\end{equation}

\subsection{Two-phase suspensions of solid particles}

When the lava is assumed to be composed of melt and crystals that occupy volumes $V_l$ and $V_p$ respectively, the particle volume fraction, $\Phi_p^*$ is then defined as in \eqref{phibp2}. We introduce $\Phi_m^*$ as the maximum solid packing fraction in the mixture, which is the densest possible concentration of solid particles. For mono-disperse spheres, $\Phi_m^* \sim 0.74$ but it is generally smaller for disordered suspensions \cite{mueller:2010}, for rougher particles, and for non-spherical particles \cite{mader:2013}. 
The viscosity must increase with the solid concentration and when $\Phi_p^*$ reaches $\Phi_m^*$ the viscosity must be infinite. This requirement was implemented first by \cite{roscoe:1952} and it was generalized by \cite{krieger:1959} by proposing the following definition of the viscosity
\begin{equation}\label{eq_eta_particles0}
    \eta=\eta_0\left(1-\frac{\Phi_p^*}{\Phi_m^*}\right)^{-B\Phi_m^*}
\end{equation}
where $B$ is the Einstein number. This relationship valid at low values of $\Phi_p^*$ is used as a basis for subsequent two-phase relationships that all try to best represent experimental data (\cite{castruccio:2010, mueller:2010} and references therein).  The original model proposed by \cite{roscoe:1952} considers $B\Phi_m^*=2.5$. One of the most popular relationship is that of \cite{maron1956application}, where the coefficient $B \Phi_m^* = 2$ is set.
\\\\
In \cite{mueller:2010}, experiments for several solid suspensions are presented together with relationships for the rheology coefficients $\eta, \tau_c, n$. They propose that
\begin{subequations}\label{rheo_mueller}
\begin{align}
	&\eta = \eta_0 \left(1-\frac{\Phi_p^*}{\Phi_m^*}\right)^{-B\Phi_m^*}, \label{eq_eta_particles}\\
	&\tau_c = \tau_c^* \left(\left(1-\frac{\Phi_p^*}{\Phi_m^*}\right)^{-2} -1\right), \label{eq_tau_mueller}\\
	& n = 1 - 0.2 r_p \left(\frac{\Phi_p^*}{\Phi_m^*}\right)^4;
	\label{eq_n_mueller}
\end{align} 
\end{subequations}
where $r_p$ is the aspect ratio of particles ($r_p=1$ for spheres).
In this case the coefficients $B$ (Einstein number) and $\tau_c^*$ (related to the particle size) are considered as fitting parameters. 
Relationship \eqref{eq_tau_mueller} is adopted from \cite{heymann:2002} on the consideration that the yield stress increases with increasing particle volume fraction. \cite{mueller:2010} found, indeed, that the yield stress is non negligible only at high particle volume fraction, $\Phi_p^*/\Phi_m^* \gtrapprox 0.8$. The value of $\tau_c^*$ is chosen to fit the experiments. For example, with spheres particles of mean radius 50~$\mu$m, they found $\tau_c^*=0.153$~Pa, for the maximum packing $\Phi_m^*=0.633$ (set to fit the viscosity in \eqref{eq_eta_particles} with $B \Phi_m^* = 2$).
The flow index $n$ is the parameter most sensitive to particle shape. Their experimental results reveal that shear thinning ($n<1$) becomes more significant as the solid volume fraction increases. This is encoded in the proposed relationship \eqref{eq_n_mueller} that ensures a decreasing $n$ with increasing $\Phi_p^*/\Phi_m^*$. The coefficient of $0.2$ and the exponent $4$ in this equation were chosen to fit the experiments, with a good fit for $\Phi_p^*/\Phi_m^* \leq 0.8$.
\\\\
A closely related model is presented in \cite{castruccio:2010}, which also considers a Herschel-Bulkley lava rheology. In this case, cubic crystals were used in the laboratory experiments. It was observed that the non-Newtonian behaviour ($\tau_c>0$ and shear thinning $n< 1$) emerges at a critical solid volume fraction, called $\Phi_c^*$, which depends on the nature of the mixture. The proposed relationships take this threshold into account:
\begin{subequations}\label{rheo_castruccio}
    \begin{align}
	&\eta = \eta_0 \left(1-\frac{\Phi_p^*}{\Phi_m^*}\right)^{-2.3},
	\label{eq_eta_castruccio}\\[4pt]
	&n =
	\left\{ \begin{array}{ll}
		\displaystyle 1 & \text{if } \Phi_p^* \leq \Phi_c^*\\
		\displaystyle 1 - 1.3 \, \frac{\Phi_p^*-\Phi_c^*}{\Phi_m^*} & \text{if } \Phi_p^* > \Phi_c^*,
	\end{array} \right.
	\label{eq_n_castruccio}\\[4pt]
	&\tau_c =
	\left\{ \begin{array}{ll}
		\displaystyle 0 & \text{if } \Phi_p^* \leq \Phi_c^*\\
		\displaystyle D \left(\Phi_p^*-\Phi_c^*\right)^8 & \text{if } \Phi_p^* > \Phi_c^*,
	\end{array} \right.
	\label{eq_tauc_castruccio}
\end{align}
\end{subequations}
with $D$ a fitting parameter. This rheology is also used in \cite{castruccio:2014} for experiments with different analogue crystal mixtures  and natural data using Etna samples. They set empirically that $\Phi_c^* = 0.44 \, \Phi_m^*$ and the obtained results are also in agreement with the rheology of \cite{mueller:2010} in \eqref{rheo_mueller}. They  note, however, that the yield stress relationship \eqref{eq_tauc_castruccio} does not adequately fit both analogue and natural samples, with the best-fit values of the parameter $D$ varying over three orders of magnitude.

\subsection{Two-phase suspensions of bubbles}
The influence of gas bubbles is analyzed considering that the lava is composed of melt and bubbles that occupy volumes $V_l$ and $V_b$, respectively. The bubble volume fraction, $\Phi_b$ is defined as in \eqref{phibp1}. \cite{llewellin2005bubble} gives a clear explanation focused on the viscosity dependence on how bubbles affect the rheology of the suspensions. A new relationship based on the concept of the relative viscosity and taking into account the capillary number is proposed and applied to a volcanic conduit flow model. The viscosity model is
\begin{equation}\label{eq_eta_bubble_Llewelin}
\eta = \left\{\begin{array}{ll}
\eta_0 \left(1-\Phi_b\right)^{-1} & \text{if } {\rm Ca} \ll 1 \\
\eta_0 \left(1-\Phi_b\right)^{5/3} & \text{if } {\rm Ca} \gg 1, 
\end{array}\right.
\end{equation}
which is valid for steady flows and for unsteady flows with ${\rm Cd}\ll 1$.
\\
This model was then generalized in \cite{mader:2013} to take both regimes into account:
\begin{equation}\label{eta_mader}
\eta =  \eta^\infty+\frac{\eta^0-\eta^\infty}{1+{\rm Cx}^s},
\end{equation}
where ${\rm Cx}=\sqrt{{\rm Ca}^2+{\rm Cd}^2}$, $s=2$ for monodisperse suspension, and $\eta^0, \eta^\infty$ are the viscosities in each regime ${\rm Ca}\ll 1$ and ${\rm Ca}\gg 1$, respectively, adopted from \eqref{eq_eta_bubble_Llewelin}.
\\\\
If the capillary number is large enough, the presence of bubbles can produce shear thinning behaviour. In the experimental results performed in \cite{truby:2015}, shear thinning is shown to occur even for the bubble suspensions without particles (see their figure 10a). This is in line with observations in \cite{stein:1992} with a fitting flow index $n\in[0.87,0.93]$, and the bubbles suspension simulations in \cite{manga:1998}, where a shear thinning behaviour is observed at small Ca ($<0.5$). \cite{truby:2015} thus propose a linear relationship for the flow index that decreases when the bubble concentration increases: 
\begin{equation}\label{eq_n_truby_bubbles}
n=1-0.334\Phi_b.    
\end{equation}
The influence of bubbles on the yield stress is only visible at high concentrations, $\Phi_b\gtrapprox 0.7$, when the suspension becomes a foam. In this case empirical laws follow the form $\tau_c=A(\Phi_b-\Phi_b^c)^2$, for some coefficient $A$ and $\Phi_b^c$ being the jamming packing fraction \cite{princen:1989, rouyer:2005, cohen:2014,namiki:2017,namiki:2022}. In this work we do not consider lavas with a high concentration of bubbles so the yield stress is independent of $\Phi_b$.

\subsection{Three phase suspensions}
The most complete model considers a Newtonian viscous fluid that contains both crystals and bubbles.
A (single) theoretical model is found in \cite{phan1997differential}, where the authors developed a differential scheme to calculate the effective viscosity of general multi-phase suspensions. In the three-phase case, it reads 
\begin{equation}\label{rheo_phan}
\eta = \eta_0 \left(1-\Phi_b\right)^{-1} (1-\Phi_p)^{-5/2},    
\end{equation}
where $\Phi_p$ can be replaced by the ratio $\Phi_p/\Phi_m^*$ in the case of a three-phase mixture. This model, however, has not been validated by experimentation.\\

In \cite{truby:2015}, a model for the rheology of a three-phase suspension is developed by combining the two-phase models of \cite{llewellin2005bubble, mueller:2010, mader:2013}. The model is validated using experimental data involving spherical glass beads. The authors focus on the low capillary number regime where the viscosity always increases with increasing bubble concentration. Furthermore, the yield stress measured in their experimental samples is negligible. This is consistent with the results of \cite{mueller:2010, mader:2013}, where the yield stress vanishes for solid concentration satisfying the condition $\Phi_p^*/\Phi_m^*<0.8$. Without yield stress, the Herschel-Bulkley rheology \eqref{HB0} reduces to a power law relationship.
The apparent viscosity of the suspension is given as a combination of the bubble model \eqref{eq_eta_bubble_Llewelin} and the particle model \eqref{eq_eta_particles}, assuming the Maron and Pierce coefficient $B \Phi_m^* = 2$,
\begin{subequations}\label{rheo_truby}
\begin{align}
	\eta = \eta_0 \left(1-\Phi_b\right)^{-1} \left(1-\frac{\Phi_p}{\Phi_m^*}\right)^{-2}.
	\label{eq_viscosity_truby}
\end{align}
The maximum packing fraction is the one proposed in \cite{mader:2013} written in terms of the particle aspect ratio $r_p$,
\begin{equation}\label{phipmax}
\Phi_m^*=\Phi_{m,s}^*\exp\left(-\frac{(\log_{10} r_p)^2}{2b^2}\right),
\end{equation}
where $\Phi_{m,s}^*$ is the maximum packing for spheres ($r_p=1$) and $b$ is a fitting parameter (for smooth particles $\Phi_{m,s}^*=0.66, b=1.08$ and for rough particles $\Phi_{m,s}^*=0.55, b=1$). The experimental data show that these relationships give promising results, except when both concentrations are large. This may be due to the inter-phase interactions, which are not accounted for in \eqref{rheo_truby}. 
This model follows the same structure as \eqref{rheo_phan} but has a different exponent in the solid concentration term.
The authors propose an empirical relationship for the flow index as a combination of the effects of both bubbles and particles from equations \eqref{eq_n_mueller} and \eqref{eq_n_truby_bubbles}, 
\begin{align}
	n = 1- 0.2 \left(\frac{\Phi_p}{\Phi_m^*}\right)^{4} - 0.334 \, \Phi_b,
	\label{eq_n_truby}
\end{align}
which is valid for a total concentration smaller than 0.5. Note that since their samples contain spherical particles, $r_p=1$ in \eqref{eq_n_mueller}.
\end{subequations}
\\\\
A purely empirical model for the rheology of a suspension with bubbles and particles is developed in \cite{birnbaum:2021}. The authors inverse the Herschel-Bulkley rheology to find relationships for $\eta, \tau_c, n$ using experimental flows for different volume fractions and a probabilistic approach. 
The empirical coefficients are:
\begin{align}
	\begin{array}{l}
		\displaystyle \eta = \eta_0 \left(1-\Phi_b^*\right)^{-B_g} \left(1-\frac{\Phi_p^*}{\Phi_m^*}\right)^{-B_s},\\[4.5mm]
		\displaystyle \tau_c = 10^{C_1(\Phi_p^* - \Phi_{c,\tau_c} )} + 10^{C_2 (\Phi_p^* + \Phi_b^* -\Phi_{c,\tau_c})};\\[1mm]
		\displaystyle n = \left\{ \begin{array}{ll}
			\displaystyle 1 & \text{if } \Phi_p^* (1 - \Phi_b^*) + \Phi_b^* \leq \Phi_{c,n},\\
			\displaystyle 1 + (C_3 - C_4 {\rm Ca}) (\Phi_{c,n} - \Phi_{p}^*(1 - \Phi_{b}^*) - \Phi_{b}^*), & \text{if } \Phi_p^* (1 - \Phi_b^*) + \Phi_b^* > \Phi_{c,n}.
		\end{array} \right.\\[5mm]
	\end{array}
	\label{rheo_birnbaum}
\end{align}
where the constants are $B_s=2.74\pm 1.56$, $B_g=1.98\pm 0.09$, $C_1=80\pm 10.9$, $C_2=1.98\pm 0.23$, $C_3=0.70\pm 0.25$, and $C_4=0.55\pm 0.31$. The viscosity follows the same structure as \eqref{rheo_truby} for the fitted constants, $\Phi_m^*=0.56\pm 0.20, B_g=1.98\pm 0.09, B_s=2.74\pm 1.56$. The yield stress appears beyond the threshold $\Phi_{c,\tau_c} = 0.35\pm 0.01$. The flow index $n$ is Newtonian below the threshold $\Phi_{c,n}=0.39\pm 0.12$ and shear-thinning or shear thickening above that threshold, depending on the capillary number. Experiments are mostly performed with bubble-rich fluids, with $0\leq \Phi_b^*\leq 0.82$ and $0\leq \Phi_p^*\leq 0.37$. In comparison with the Truby rheology in equation \eqref{rheo_truby}, the exponents in the viscosity differ by one. 
In the low capillary regime and for the two-phase case, with either bubbles or particles at a concentration $\Phi^*$, the flow index in \eqref{rheo_birnbaum} becomes
$$
\displaystyle n = \left\{ \begin{array}{ll}
			\displaystyle 1 & \text{if } \Phi^* \leq \Phi_{c,n},\\
			\displaystyle 1 - 0.7 (\Phi^*-\Phi_{c,n}), & \text{if }  \Phi^* > \Phi_{c,n}.
		\end{array} \right.
$$
In the case of a bubble mixture, the main difference with equation \eqref{eq_n_truby_bubbles} is the consideration of a threshold for the bubble volume fraction below which $n=1$. In a particle-laden suspension, this relationship is closer to that of \cite{castruccio:2010} in equation \eqref{eq_n_castruccio} than to the Mueller formula \eqref{eq_n_mueller}.

\subsection{Solid concentration, crystallization and temperature influence}

In \cite{marsh:1981} the crystallization in magma is analyzed as a function of temperature. Using natural data of different magmas, the author relates the eruption probability to the thermal probability and the rheological probability (which is related to the viscosity). From this analysis a constitutive equation for the solid volume fraction $\Phi_p^*$ is determined by the expression
\begin{equation}\label{eq_phip_marsh}
\Phi_p^* = \frac12\left(1-{\rm erf}\left(\frac{b_1}{\sigma_1}\Big(\frac{T-T_s}{T_l-T_s}- \theta_1\Big)\right)\right),
\end{equation}
with constant coefficients $b_1,\sigma_1, \theta_1$. In particular $\theta_1$ is the dimensionless temperature at which one half of the magma has crystallized, and $\sigma_1$ is related to the standard deviation of the temperature about the most probable state of crystallinity.
This equation is also used in lava flows models  \cite{tsepelev:2020,vetere:2021}.
\\\\
In \cite{dragoni:1994} a two-phase model with a Bingham rheology (that is, $n=1$) that takes into account the effect of temperature on crystallization is studied. The temperature is assumed to be constant in time and evolves according to:
\begin{align}\label{eq_temperature_dragoni}
	\partial_x T = -\frac{\varepsilon \sigma}{\rho c_p h u} T_e^4,
\end{align}
where $\rho$ is the lava density, $\varepsilon$ is the emissivity, $\sigma$ is the Stefan-Bolztmann constant, $c_p$ is the specific heat, $h$ is the lava thickness (or height) and $u$ its speed. The relationship of the solid fraction $\Phi_p$ as a function of temperature is assumed to be linear:
\begin{align}\label{eq_phip_dragoni}
	\Phi_p^* = \frac{T_l-T}{T_l-T_s} \Phi_p^{m,*}, \qquad T_s \leq T \leq T_l,
\end{align}
where $T_l$ is the liquidus temperature, $T_s$ the solidus temperature and $\Phi_p^{m,*}$ is the crystallization degree at $T=T_s$. When the lava is totally melted, $T=T_l$ and there is no crystals. Conversely, when $T=T_s$ the lava is totally solidified and $\Phi_p^*=\Phi_p^{m,*}$. 
The effective viscosity follows an exponential law depending on $T$ and $\Phi_p^*$:
\begin{align}\label{eq_eta_dragoni}
	\eta = \eta_a e^{a/T} + \eta_0 \left(e^{b\Phi_p^*}-1\right), 
\end{align}
with empirical parameters $\eta_a$, $a$ and $b$.
The first term describes the direct temperature dependence of the viscosity (which decrease exponentially with the temperature). The second one represents the effect of crystallization, which vanishes when there are no crystals ($\Phi_p^*=0$).
\\\\
In \cite{crisp1994influence}, the influence of crystallization on the evolution of lava flow temperature is studied in a theoretical model. The solid concentration is a function of time and it is not written in terms of the temperature but of the crystallization time $\tau_{\rm crys}$:
\begin{align}\label{eq_phip_crisp}
	\Phi_p^* = \Phi_{\rm crys}^* \left(1-e^{-t/\tau_{\rm crys}}\right),
\end{align}
where $\tau_{\rm crys}$ is a characteristic time for the lava to approach its maximum volume fraction of crystallization, $\Phi_{\rm crys}^*$. In other words, $\Phi_p^*$ increases with time until the maximum value $\Phi_{\rm crys}^*$ is reached. Values of $\tau_{\rm crys} \approx 0.5-10$ days are considered in \cite{crisp1994influence}. The temperature is assumed to depend only on time and its evolution equation reads
\begin{align}\label{eq_temperature_crisp}
	\partial_t T = \frac{L}{c_p} \partial_t \Phi_p^* - \frac{\varepsilon \sigma f}{\rho c_p h} T^4 + \frac{T-T_e}{\tau_e} + \frac{L}{c_p} \frac{\Phi_p^*-\Phi_e^*}{\tau_e}.
\end{align}
The terms on the right-hand side represent, in order, the heat added due to crystallization, the radiation cooling, the heat entrainment balance ($T_e$ is the temperature of the entrained material) and finally, the heat loss due to the melting of the entrained material ($\Phi_e^*$ is the crystal concentration in the entrained material). The constants appearing are the latent heat of crystallization $L$, the area fraction of exposed lava core at the surface $f$, the lava density $\rho$, the emissivity $\varepsilon$, the Stefan-Bolztmann constant $\sigma$, the heat capacity of lava in the core $c_p$, and the entrainment timescale $\tau_e$. 
\\\\
\cite{costa_macedonio_2005} propose a depth-averaged model for lava flows where crystallization is not explicitly considered. The evolution of temperature is given by:
\begin{equation} \label{eq_temperature_costa}
	\partial_t (h T)+ \partial_x (h \,T \, u)+ \partial_y (h\, T \, v)= Q_{rad}+Q_{conv}+Q_{cond}+Q_{visc},
\end{equation}
where $Q_{rad}$, $Q_{conv}$, and $Q_{cond}$ represent the radiative, convective and conductive temperature exchanges, respectively, and $Q_{visc}$ represents the viscous heating. These terms are given by:
\begin{align}
	& Q_{rad}= \frac{\epsilon \, \sigma \, f}{\rho c_p} (T_{env}^4-T^4) ,
	&&Q_{conv} = \frac{\lambda f}{\rho c_p} (T_{env}-T),\\
	& Q_{cond} = \frac{n_0 k}{\rho c_p h} (T_{c}-T) ,
	&& Q_{visc} = \frac{m \eta}{\rho c_p h} (u^2+v^2).
	\label{eq_heat_terms}
\end{align}
where $T_{env}$ is the external temperature, $T_c$ the temperature of the ground, $\epsilon$ is the emissivity coefficient (dimensionless), $f$ is the area of the boundary between lava and air, $\sigma$ is the Stefan-Boltzman constant, $\lambda$ is the atmospheric heat transfer coefficient, $k$ is the thermal conductivity and $n_0,m$ are dimensionless coefficients. As in \cite{dragoni:1994} the viscosity follows an exponential law, but it depends solely on temperature,
\begin{align}
	\eta = \eta_0 e^{-b(T-T_0)},
\end{align}
where $\eta_0$ is the viscosity at temperature $T_0$ and $b$ is a rheological parameter.
\\
The viscous heating term $Q_{visc}$ in \eqref{eq_heat_terms} is obtained as an approximation of the following term 
\begin{equation}\label{Qvisc_integral}
Q_{visc}\sim \frac{1}{\rho c_p} \int_{z_b}^{z_b+h} \eta |\partial_z (u,v)|^2 dz'.     
\end{equation}
In particular, they estimate the characteristic velocity boundary layer $\delta_v=\frac{h}{m}$ to write $\int_{z_b}^{z_b+h} \eta |\partial_z (u,v)|^2 dz' \sim \eta \frac{u^2+v^2}{\delta_v}= \frac{m\eta}{h}(u^2+v^2)$. The value $m=12$ is set for a parabolic profile of the velocity.

\section{The proposed lava flow model}

In this section we consider a 2-D model for a three-phase lava flow, with a Herschel-Bulkley rheology. We combine rheology relationships from the literature reviewed in the previous section. Then, we derive a depth-averaged version of the model in section \ref{sec_derivation_shallow} that is amenable to discretization.\\

The following system expresses the mass and momentum conservation in two dimensions for spatial variables $x$ and $z$ over an inclined plane with slope $\theta$. The velocity is ${\bf U}=(u,w)$, $\rho$ is the bulk density of the lava, $p$ is the pressure, $g$ is the gravity constant, and $\sigma$ is the shear stress tensor. The system is written in tilted coordinates with respect to the plane.
\begin{subequations}
  \begin{empheq}[left=\empheqlbrace]{align}
		& \displaystyle \partial_x u+\partial_z w =0, \label{eq_2D_uw0} \\
		& \displaystyle \rho \partial_t u + \rho \partial_x ( u^2) + \rho \partial_z (u w) + \partial_x p = \rho g \sin \theta + \partial_x \sigma_{xx} + \partial_z \sigma_{xz} \label{eq_2D_u0}	\\
		& \displaystyle \rho \partial_t (w) + \rho \partial_x ( uw) + \rho \partial_z (w^2) + \partial_z p = -\rho g \cos \theta + \partial_x \sigma_{xz} + \partial_z \sigma_{zz}. \label{eq_2D_w0}
	\end{empheq}
\label{eq_2D0}
\end{subequations}
The shear stress tensor $\sigma$ is given by the Herschel-Bulkley law \eqref{HB0}:
\begin{align}
	\left\{ 
		\begin{array}{ll}
			\displaystyle |\sigma| \leq \tau_c & \text{if } |D(\mathbf{U})| = 0\\
			\displaystyle \sigma = \mu D(\mathbf{U}) &  \text{if } |D(\mathbf{U})| \neq 0
		\end{array},
		\qquad \mu = \eta |D(\mathbf{U})|^{n-1} + \frac{\tau_c}{|D(\mathbf{U})|},
	\right.
	\label{eq_Herschel-Bulkley}
\end{align}
where $n>0$ is the flow index and $\tau_c$ is the yield stress (a threshold stress above which the flow behaves in a non-Newtonian fashion). For the case $n=1$,  $\eta$ is the apparent viscosity and in other cases it represents the consistency factor. Finally, $D(\mathbf{U})$ is the deformation tensor:
\begin{align*}
	D(\mathbf{U}) = \nabla \mathbf{U} + \nabla \mathbf{U}^t=
	\displaystyle \left( \begin{array}{cc}
		2 \partial_x u & \partial_x w + \partial_z u\\
		\partial_x w + \partial_z u & 2 \partial_z w
	\end{array} \right).
\end{align*}
The shear rate is then $|\dot \gamma| = |D({\bf U})|$.
Note that when $n=1$ and $\tau_c=0$ we recover a Newtonian behavior, with viscosity $\mu=\eta$. In the case $\tau_c=0$, the Herschel-Bulkley model reduces to a single power law.

\subsection{Closure for rheological parameters} \label{sec_closure_rheology}

We combine expressions for $\eta$, $n$ and $\tau_c$ coming from the literature and listed in Section \ref{SectRheol} so that they become suitable for our three-phase model. The volume fraction of bubbles $\Phi_b$ and crystals $\Phi_p$ are defined as in equation \eqref{phibp1}.\\

Following the approach outlined in \cite{truby:2015}, we propose a generalized rheological model that includes the effects of the aspect ratio of the solid particles and of the capillary number for bubbles. This generalization is based on the relationships presented in Section \ref{SectRheol} and is tested against some experimental data from the literature. As discussed before, the dependence of the viscosity and yield stress on particle aspect ratio is implicitly considered when setting the value of the maximum packing fraction. This dependence therefore only appears in the flow index. The capillary number clearly influences the apparent viscosity but its influence on the yield stress and the flow index has hardly been studied in magmas. To our knowledge, the only results are those presented in \cite{birnbaum:2021}, but the number of cases in their sampling is insufficient to establish a conclusive behaviour\footnote{There are only two samples with $\Phi_p=0$ and Ca$>1$. One yields shear thinning ($n=0.79$) and the other yields shear thickening ($n=1.32$).} (see their figure 3f).\\

\begin{figure}[hbtp]
	\centering      
	\includegraphics[width=0.5\textwidth]{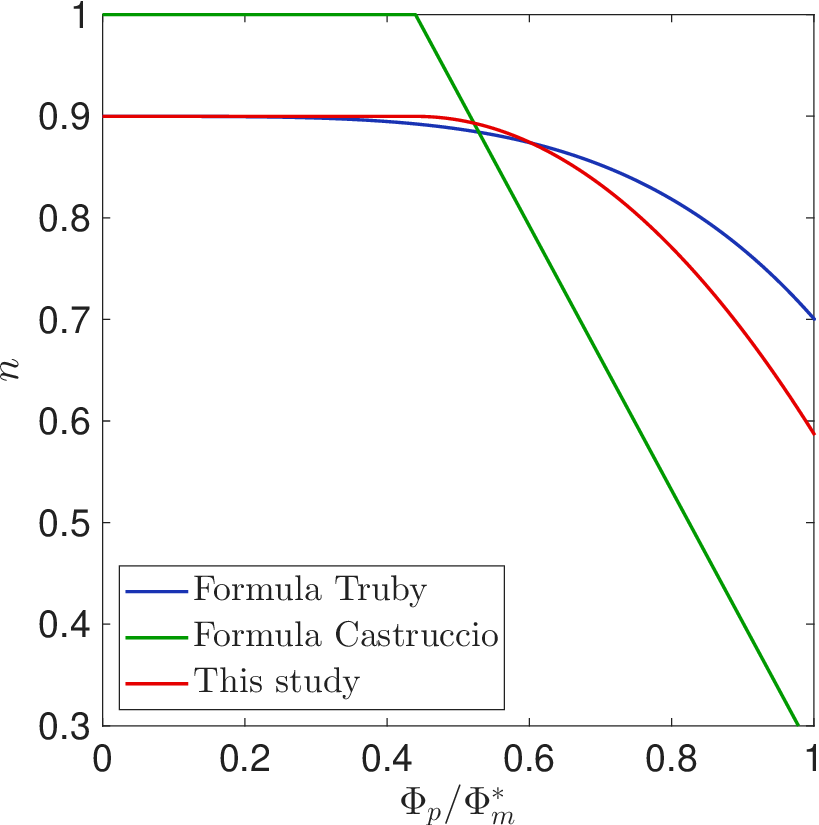}
	\caption{Comparison of flow index number $n$ depending on the solid volume fraction, with formulas in equations \eqref{eq_n_truby} (in blue), \eqref{eq_n_castruccio} (in green) and our new equation \eqref{n_prop} (in red). Values of parameters are $\Phi_m^*=0.73$, $\Phi_p^c=0.44 \Phi_m^*$, $\Phi_b=0.3$, $\alpha_n=1$.}
	\label{fig_comparison_n}
\end{figure}

In agreement the reasoning of \cite{truby:2015}, our apparent viscosity is a combination of the viscosity for solid suspensions in \cite{mueller:2010}\footnote{Note that we followed \cite{truby:2015} by using $\Phi_p$ instead of $\Phi_p^*$.} and that for bubbles suspensions in \cite{mader:2013}, given respectively in equations \eqref{rheo_mueller} and \eqref{eta_mader}: 
\begin{subequations}\label{rheo_prop_gral}
\begin{equation}
    \eta = \frac{\eta_0}{1+{\rm Cx}^2} \Big( (1-\Phi_b)^{-1}+{\rm Cx}^2 (1-\Phi_b)^{5/3}\Big) 
    \left(1-\frac{\Phi_p}{\Phi_m^*}\right)^{-2},
\end{equation}
where the capillary number is ${\rm Cx}=\sqrt{{\rm Ca}^2+{\rm Cd}^2}$ and $\Phi_m^*$ is given by \eqref{phipmax},
\begin{equation}
\Phi_m^*=\Phi_{m,s}^*\exp\left(-\frac{(\log_{10} r_p)^2}{2}\right),
\end{equation}
where $\Phi_{m,s}^*$ is the maximum packing for spheres ($r_p=1$). We choose $b=1$ (rough particle surface) in \eqref{phipmax}.\\
At low bubble concentrations, the yield stress is only affected by particle concentration. We consider the law proposed in \cite{mueller:2010} given in equation \eqref{eq_tau_mueller} and ensure that it yields a bounded value:
\begin{equation}\label{eq_tau_particles0}
    \tau_c = \min \left(\tau_c^{\max}, \tau_c^* \left(\left(1-\frac{\Phi_p^*}{\Phi_m^*}\right)^{-2} -1 \right) \right), 
\end{equation}
where $\tau_c^*$ is a fitting parameter and $\tau_c^{\max} \approx 10^9-10^{10}$~Pa is the Young modulus of crystals that bounds the maximum physically sound value for $\tau_c$ (cf. \cite{mao2015elasticity}). The crystal fraction definition in this formula is $\Phi_p^*$ (defined by equation \eqref{phibp2}), as the relationship in \cite{mueller:2010} is for a pure solid suspension to which we add bubbles. We do not consider the threshold for the solid concentration as in \cite{castruccio:2010} to find the limit case $\tau_c=0$ (see equation \eqref{eq_tauc_castruccio}). This is because equation \eqref{eq_tau_particles0} works well even at small concentrations when compared with experimental data (see figures 9c, 9d in \cite{mueller:2010}). \\
We propose a new relationship for the flow index $n$. \cite{truby:2015} show that this parameter is always affected by bubble concentration (see equation \eqref{eq_n_truby_bubbles}). Combining it with the Mueller model for particles gives the flow index in equation \eqref{eq_n_truby}, which yield values $n<1$ for any suspension. Figure 10a in \cite{truby:2015} shows that the experimental values of the flow index $n$ fit the line $1- 0.334 \Phi_b$ for $\Phi_p\leq 0.4$ but that above that particle concentration the flow index behaves differently. We take the fact that $n$ is independent of $\Phi_p$ when $\Phi_p\leq 0.4$ into account by following the suggestion of \cite{castruccio:2014}, which consists in introducing a threshold solid volume fraction $\Phi_p^c$ to cancel the particle influence on $n$ when $\Phi_p \leq \Phi_p^c$, see equation \eqref{eq_n_castruccio}. 
Adding the effect of the particle aspect ratio as in \eqref{eq_n_mueller}, we thus propose the following flow index
\begin{equation}
n = \left\{ \begin{array}{ll}
	\displaystyle 1- 0.334 \Phi_b & \text{if } \Phi_p \leq \Phi_p^c,\\
			\displaystyle 1- \alpha_n\, r_p \left(\frac{\Phi_p-\Phi_p^c}{\Phi_m^*}\right)^{4} - 0.334 \Phi_b & \text{if } \Phi_p > \Phi_p^c,
\end{array} \right.
\end{equation}
where $\alpha_n$ is a fitting parameter. In figure \ref{fig_comparison_n}, we represent the existing flow index relationships alongside ours for the case of spheres ($r_p=1$).

\end{subequations}

\subsubsection*{Particular case: spherical particles and low capillary number}
In the low capillary number regime (${\rm Cx}\ll 1$) and for spherical particles ($r_p=1$), the previous relationships become:
\begin{subequations}\label{rheo_prop}
\begin{equation}\label{eta_prop}
    \eta = \eta_0(1-\Phi_b)^{-1}
    \left(1-\frac{\Phi_p}{\Phi_m^*}\right)^{-2},
\end{equation}
\begin{equation}\label{tauc_prop}
    \tau_c = \min \left(\tau_c^{\max}, \tau_c^* \left(\left(1-\frac{\Phi_p^*}{\Phi_m^*}\right)^{-2} -1 \right) \right),
\end{equation}
\begin{equation}\label{n_prop}
n = \left\{ \begin{array}{ll}
	\displaystyle 1- 0.334 \Phi_b & \text{if } \Phi_p \leq \Phi_p^c,\\
			\displaystyle 1- \alpha_n \left(\frac{\Phi_p-\Phi_p^c}{\Phi_m^*}\right)^{4} - 0.334 \Phi_b & \text{if } \Phi_p > \Phi_p^c,
\end{array} \right.
\end{equation}
where $\Phi_m^*=\Phi_{m,s}^*$ is the maximum packing for spheres in presence of bubbles, and is equivalently $\Phi_m^*$ using \eqref{phibp_relation}. Coefficients $\tau_c^*,\alpha_n, \Phi_p^c$ are fitting parameters.
\end{subequations}

\subsection{Closures of bubbles and crystals concentrations}

Equations \eqref{eq_2D0} and \eqref{eq_Herschel-Bulkley} are closed using the expressions for $\eta$, $n$ and $\tau_c$ given in the previous section. These expressions depend on $\Phi_p$ and $\Phi_b$. We consider a constant concentration of bubbles $\Phi_b$ and we follow previously published models to set a constitutive equation for the crystal fraction $\Phi_p$.
\\
Our first approach is inspired that of \cite{crisp1994influence} where the crystals volume fraction $\Phi_p^*$ evolves in time following equation \eqref{eq_phip_crisp}.
If $\tau_{\rm crys}$ and $\Phi_{\rm crys}^*$ are constant, the equation \eqref{eq_phip_crisp} is the solution of the following ODE with initial condition $\Phi_p^*(t=0)=0$:
\begin{align*}
	\partial_t \Phi_p^* = - \frac{1}{\tau_{\rm crys}} \left(\Phi_p^*-\Phi_{\rm crys}^*\right).
\end{align*}
This formula can be build upon to take into account the advection caused by the lava flowing at a velocity ${\bf U}$:
\begin{align}\label{eq_phi_p_adv}
	\partial_t \Phi_p^* + \partial_x \left({\bf U} \Phi_p^*\right) = - \frac{1}{\tau_{\rm crys}} \left(\Phi_p^*-\Phi_{\rm crys}^*\right).
\end{align} 
The original equation \eqref{eq_phip_crisp} is also the analytic solution of \eqref{eq_phi_p_adv} for the static case and initial condition $\Phi_p^*(t=0)=0$. This is the same equation used in \cite{tsepelev:2020, zeinalova:2025} to solve the solid concentration in their lava flow model. 
The parameter $\Phi_{\rm crys}^*$ depends on the temperature and on the amount of dissolved water in the magma \cite{tsepelev:2020, laspina:2016}. The works of \cite{tsepelev:2020, zeinalova:2025} consider $\Phi_{\rm crys}^*$ as a constant with values of $0.8$ and $0.43, 0.63, 0.83$, respectively.
In our case we consider that the maximum volume percent crystallization $\Phi_{\rm crys}^*=1$ to represent crystallization towards total solidification.\\

\begin{figure}[hbtp]
	\centering      
	\includegraphics[width=0.7\textwidth]{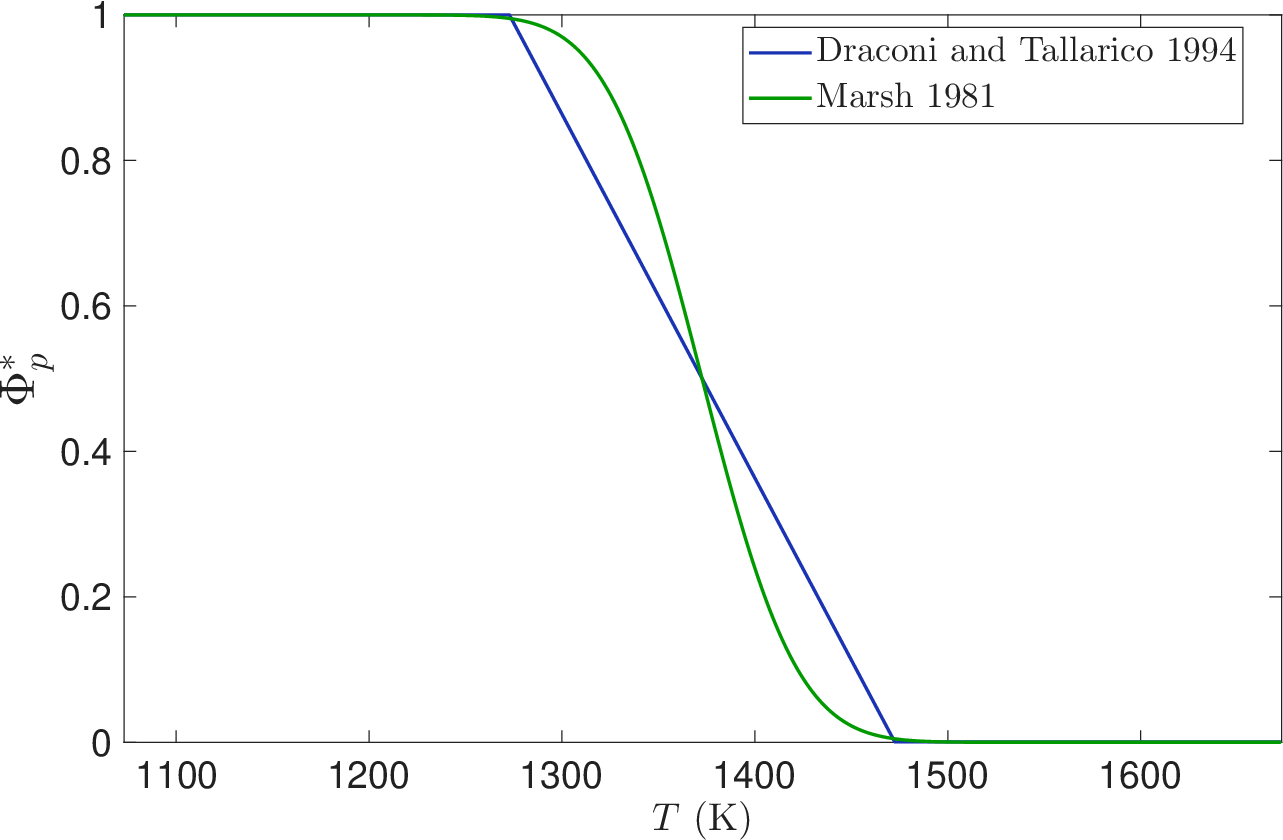}
	\caption{Comparison of solid fraction formulas in equations \eqref{eq_phip_marsh} and \eqref{eq_phip_T_lin0} for $\Phi_{\rm crys}^*=1$, $T_l=1473$ K, $T_s=1273$ K and coefficients $\theta_1=1/2, b_1=\sqrt{30}, \sigma_1=3/2$ from \cite{tsepelev:2020}.}
	\label{fig_comparison_phip}
\end{figure}

Equation \eqref{eq_phi_p_adv} makes the assumption that the crystallization time $\tau_{\rm crys}$ is constant. Its use is thus restricted to, for instance, small lava flows. To address more general flows, one should assume that crystallization depends on the lava temperature $T$. We considered such general flows in our second approach that involves an evolution equation of the temperature, and a constitutive equation relating $\Phi_p^*$ to $T$. We choose the very complete temperature equation from \cite{costa_macedonio_2005} given in equation \eqref{eq_temperature_costa} for the added advantage that it is already depth-averaged. Note that a similar approach was considered in \cite{tsepelev:2020} to model magma flow in a lava dome. One difference with our approach is that they neglected viscous heat dissipation because of the slow lava dome growth.\\
\\
In our model, the solid fraction is either that proposed by \cite{marsh:1981} and given in equation \eqref{eq_phip_marsh}, or that of \cite{dragoni:1994} given in equation \eqref{eq_phip_dragoni}. In the latter case, we consider here that the total crystallization degree is achieved at the solidus temperature, when $T=T_s$, so we set $\Phi_p^{m,*}=\Phi_{\rm crys}^*=1$, in line with our previous assumption in equation \eqref{eq_phi_p_adv}. We slightly modified this formula to ensure its validity for all positive temperature values:
\begin{align}\label{eq_phip_T_lin0}
	\Phi_p^* = \Phi_{\rm crys}^*\,\min \left(1, \frac{(T_l-T)_+}{T_l-T_s} \right), \quad T \geq 0.
\end{align}
When $T_s\leq T\leq T_l$, this equation is the same as the original. When $T>T_l$, $\Phi_p^*=0$ since the lava is totally molten, and if $T<T_s$ then $\Phi_p^*=1$, which is fully crystallized. Figure \ref{fig_comparison_phip} shows a comparison of the two relationships of $\Phi_p^*$ in equations \eqref{eq_phip_marsh} and \eqref{eq_phip_T_lin0}.\\

In summary, in both of our models, crystallization influences the rheology (viscosity, flow index, and yield stress) through the crystal concentration relationships. In the first approach, equation \eqref{eq_phi_p_adv} stipulates crystallization time, whereas in the second approach, the relationship \eqref{eq_phip_T_lin0} controls the temperature evolution, which in turn controls crystallization. 

\subsection{Validation of the choices made about rheology}

\begin{figure}[hbtp]
	\centering
	\includegraphics[width=\textwidth]{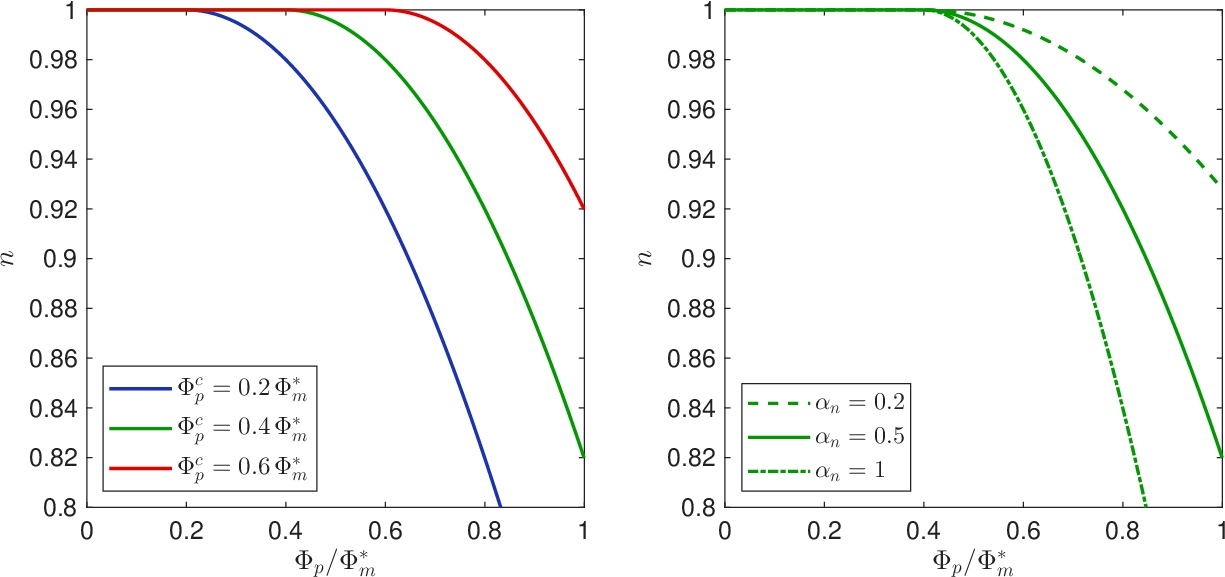}
	\caption{Influence of parameters $\Phi_p^c$ and $\alpha_n$ on flow index formula in equation \eqref{n_prop} when $\Phi_b=0$. On the left $\Phi_p^c$ varies and $\alpha_n=0.5$, on the right $\alpha_n$ varies and $\Phi_p^c=0.4 \Phi_m^*$, with $\Phi_m^*=0.56$.}
	\label{fig_compar_n}
\end{figure}

As we have combined different models to construct our rheology, it is important to check its validity by comparing it with previous studies and experimental data. Our rheology is given by equation \eqref{rheo_prop} where we only modified the flow index $n$ compared to previous models (the viscosity is that of \cite{truby:2015} for three-phase model and the yield stress is that of \cite{mueller:2010} for solid suspensions, applied here under the assumption of small bubble concentration). To calibrate the parameters $\alpha_n$ and $\Phi_p^c$ in the flow index relation \eqref{n_prop}, we compare the output of this formula with the experiments of \cite{truby:2015} for bubbles and solid suspensions, and with their respective models. This comparison is performed by keeping Ca$<1$ and considering spheres.
\\
In figure \ref{fig_compar_n} we represent the output of our relationship \eqref{n_prop} as a function of the normalized particle fraction $\Phi_p/\Phi_m^*$ without bubbles. The left panel shows $n$ for $\alpha_n=0.5$ and different values of $\Phi_p^c$, and the right panel shows $n$ for $\Phi_p^c=0.4\Phi_m^*$ and various $\alpha_n$. The parameter $\Phi_p^c$ sets the critical value above which shear thinning occurs ($n<1$). The coefficient $\alpha_n$ influences the rate of decrease of $n$, yielding larger rates at high $\alpha_n$ values.\\

\begin{figure}[hbtp]
	\centering
	\includegraphics[width=\textwidth]{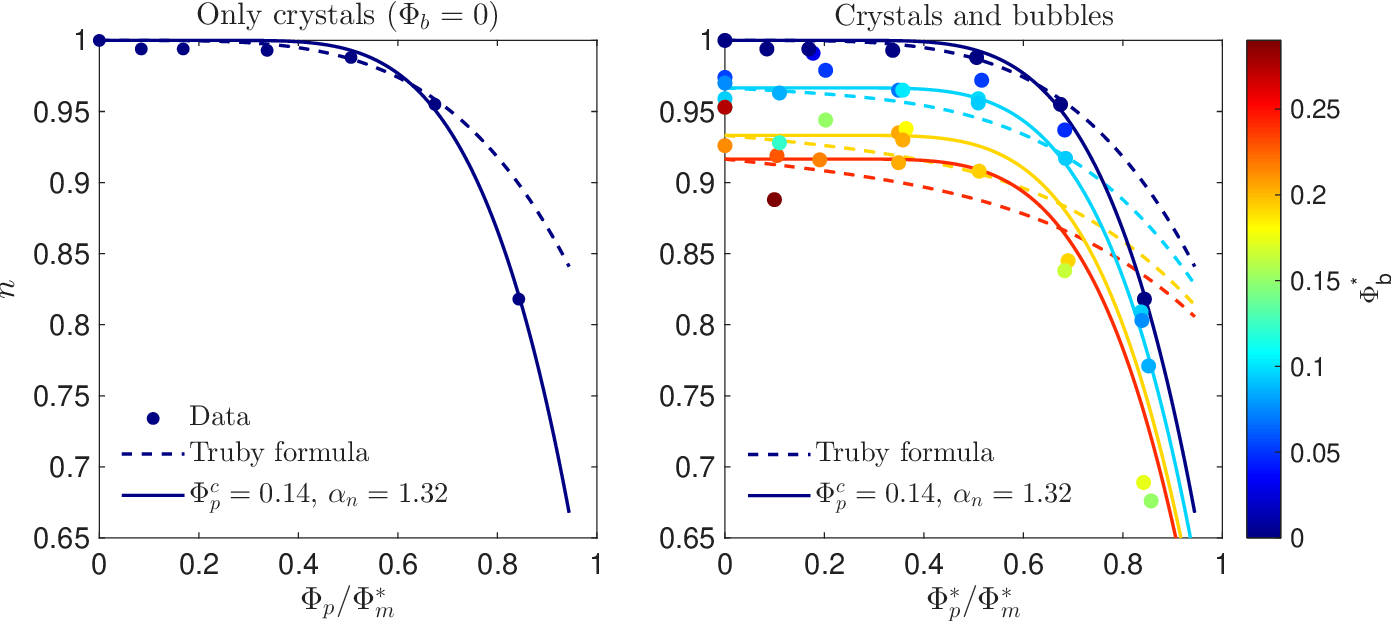}
	\caption{Comparison of data from \cite{truby:2015} with the flow index $n$ of their rheology (equation \eqref{eq_n_truby}, dashed lines) and this study best fit rheology (equation \eqref{n_prop} with $\Phi_p^c=0.14$, $\alpha_n=1.32$, $\Phi_m^*=0.593$, continuous lines). On the left, model with only crystals. On the right, model with crystals and bubbles. The bubble fraction is represented by the color bar that spans from $0$ (blue) to $0.3$ (red). The values of the bubble concentration for curves are: blue $\Phi_b=0$, cyan $\Phi_b=0.1$, yellow $\Phi_b=0.2$, red $\Phi_b=0.25$.}
	\label{fig_compar_Truby2}
\end{figure}

We used the experiments with crystals and bubbles of \cite{truby:2015} to calibrate our formula \eqref{n_prop}.
In figure \ref{fig_compar_Truby2} we compare the data and formulas of \cite{truby:2015} with our flow index relationship up to the maximum $\Phi_b$ value of 0.3 to remain within the validity domain of their relationship. We show the results for the two-phase experiments with only crystals, and for the three-phase flow with crystals and bubbles. The maximum crystal fraction is chosen as in \cite{truby:2015}, $\Phi_m^*=0.593$. Note that when only bubbles are present, our formula coincides with that of Truby. We use these data on solid suspensions to calibrate $\Phi_p^c$ and the coefficient $\alpha_n$ in the flow index equation ($\alpha_n = 0.2$ in Truby). We find a best fit for $\Phi_p^c=0.14$ and $\alpha_n = 1.32$ by using an orthogonal distance regression fitting procedure.
These modifications cause lower $n$ values at high solid volume fractions, which in \cite{truby:2015} are disregarded because they consider $\tau_c=0$ (as proposed by \cite{mueller:2010} for $\Phi_p^*/\Phi_m^*<0.8$). \\

\subsection{Derivation of the shallow water model} \label{sec_derivation_shallow}

We develop the depth-averaged model from the proposed 2D model, which we repeat here for clarity. 
\begin{subequations}\label{model_2D}
  \begin{empheq}[left=\empheqlbrace]{align}
		& \displaystyle \partial_x u+\partial_z w =0, \label{eq_2D_uw} \\
		& \displaystyle \rho \partial_t u + \rho \partial_x ( u^2) + \rho \partial_z (u w) + \partial_x p = \rho g \sin \theta + \partial_x \sigma_{xx} + \partial_z \sigma_{xz}, \label{eq_2D_u}	\\
		& \displaystyle \rho \partial_t w + \rho \partial_x ( uw) + \rho \partial_z (w^2) + \partial_z p = -\rho g \cos \theta + \partial_x \sigma_{xz} + \partial_z \sigma_{zz}. \label{eq_2D_w}
	\end{empheq}
\end{subequations}
The shear stress tensor $\sigma$ is given by the Herschel-Bulkley law \eqref{HB0}:
\begin{align}\label{HB1}
	\left\{ 
		\begin{array}{ll}
			\displaystyle |\sigma| \leq \tau_c & \text{if } |D(\mathbf{U})| = 0\\
			\displaystyle \sigma = \mu D(\mathbf{U}) &  \text{if } |D(\mathbf{U})| \neq 0
		\end{array},
		\qquad \mu =  \eta |D(\mathbf{U})|^{n-1} + \frac{\tau_c}{|D(\mathbf{U})|},
	\right.
\end{align}
where $D({\bf U})$ is the deformation tensor, $\dot\gamma=|D({\bf U})|$ is the shear rate, and the rheological coefficients are given in equations \eqref{eta_prop}, \eqref{tauc_prop} and \eqref{n_prop}, they depend on bubble and crystal concentrations.
The bubble concentration $\Phi_b$ is assumed constant and we consider two options for the solid concentration:\\
{\bf Model $1$.} The crystal fraction $\Phi_p^*$ evolves as equation \eqref{eq_phi_p_adv} at a rate controlled by the constant characteristic crystallization time $\tau_{\rm crys}$ until the maximum crystallization fraction $\Phi_{\rm crys}^*$ is reached:
\begin{align}\label{eq_phi_p_adv_bis}
	\partial_t \Phi_p^* + \partial_x \left({\bf U} \Phi_p^*\right) = - \frac{1}{\tau_{\rm crys}} \left(\Phi_p^*-\Phi_{\rm crys}^* \right).
\end{align}
{\bf Model $2$.} The crystal fraction $\Phi_p^*$ depends on the lava temperature following equation \eqref{eq_phip_T_lin0}. The temperature is computed using the depth-integrated equation \eqref{eq_temperature_costa}, where we choose a coherent approximation of the viscous heating term $Q_{visc}$ with our system. That is, 
\begin{align}\label{eq_phip_T_lin}
	\Phi_p^* = \Phi_{\rm crys}^*\,\min \left(1, \frac{(T_l-T)_+}{T_l-T_s} \right), \quad T \geq 0,
\end{align}
where 
\begin{equation} \label{eq_temperature_costa_adap}
	\partial_t (h T)+ \partial_x (h \,T \, u)+ \partial_y (h\, T \, v)= Q_{rad}+Q_{conv}+Q_{cond}+Q_{visc}.
\end{equation}
The terms $Q_{rad}$, $Q_{conv}$, and $Q_{cond}$ are taken as in \eqref{eq_heat_terms}, but we adopt the viscous heating term $Q_{visc}$ using the approximation introduced in equation \eqref{Qvisc_integral} as follows
\begin{equation}\label{Qvisc_integral2}
Q_{visc}= \frac{1}{\rho c_p} \int_{b}^{b+h} \eta |\partial_z (u,w)|^2 dz'.     
\end{equation}
\\\\
We need to introduce the boundary conditions to close the model.
At the free surface ($z=b+h$) we impose the kinematic condition and we consider no surface tension:
\begin{align}
	\partial_t (h+b) + \mathbf{U}\cdot n^h = 0, \qquad (\sigma - p I) \, n^h = 0, \qquad \text{at $z = b+h$}.
\end{align}
At the bottom ($z=b$), we impose the non-penetration condition and we consider a viscous shear for a coefficient $k_b$:
\begin{align}
	\mathbf{U}\cdot n^b = 0,\qquad ((\sigma - p I) \, n^b)_{\rm tg} = -k_b\,{\bf U}_{\rm tg}, 
	\quad \text{at $z = b$.}
\end{align}    \label{bc_bottom}
The normal vectors to the free surface and to the bottom surface, $n^h$ and $n^b$, are defined as follows:
\begin{align*}
	n^h = \frac{1}{\sqrt{1+\partial_x (h+b)^2}}
	\left( \begin{array}{c}
		 \partial_x (h+b)\\ -1
	\end{array} \right),
	\qquad
	n^b = \frac{1}{\sqrt{1+(\partial_x b)^2}}
	\left( \begin{array}{c}
		\partial_x b\\ -1
	\end{array} \right).
\end{align*}

\paragraph{Dimensional analysis.}
We introduce as usual the shallow hypothesis by considering that the ratio between the characteristic height $H$ and length $L$ is small, that is,  $\varepsilon = H/L \ll 1$. The dimensionless variables are denoted with tilde and given by:
\begin{align*}
	& (x,z,t) = \left( L \widetilde{x},H \widetilde{z}, \frac{L}{U} \widetilde{t} \right),
	\quad (u,w) = U\left( \widetilde{u}, \epsilon \widetilde{w} \right), \quad  h = H \widetilde{h},\quad b = H \widetilde{b},\\ & p = \rho g H \widetilde{p}, \quad \tau_{\rm crys}=\frac{U}{L}\widetilde{\tau}_{\rm crys}.
\end{align*}
We also introduce the dimensionless numbers: Froude number,  Reynolds number and Bingham number:
\begin{align*}
	Fr = \frac{U}{\sqrt{gH}}, \quad Re = \frac{U^{2-n}H^n}{\eta/\rho}, \quad Bi = \frac{\tau_c}{\eta} \left(\frac{H}{U}\right)^n.
\end{align*}
Omitting henceforth the tilde notation for simplicity, the dimensionless equations are:
\begin{align}
	\left\{\begin{array}{l}
		\displaystyle \partial_x u + \partial_z w =0, \\
		\displaystyle \partial_t u + \partial_x u^2 + \partial_z (u w) + \frac{\partial_x p}{\fr^2} - \frac{\sin \theta}{\varepsilon Fr^2} 
		= \frac{2\varepsilon}{Re} \partial_x \left(\mu \partial_x u \right) + \frac{\varepsilon}{Re} \partial_z \left(\mu \left(\partial_x w + \frac{1}{\varepsilon^2} \partial_z u \right)\right), \\
		\displaystyle \varepsilon^2 \left(\partial_t w + \partial_x ( uw) + \partial_z (w^2) \right) + \frac{\partial_z p}{Fr^2} + \frac{\cos \theta}{Fr^2} 
		= \frac{\varepsilon}{Re} \partial_x \left(\mu \left(\varepsilon^2 \partial_x w + \partial_z u \right)\right) + \frac{2\varepsilon}{Re} \partial_z \left(\mu \partial_z w \right), 
	\end{array} \right.	
	\label{eq_2d_epsilon}
\end{align}
where 
\begin{equation}\label{mu_adim}
    \mu=|D({\bf U})|^{n-1}+Bi\, \frac{1}{|D({\bf U})|}.
\end{equation}
Equation \eqref{eq_phi_p_adv_bis} for the model 1 does not change.
The dimensionless boundary conditions read:
\begin{subequations}
	\begin{empheq}{align}
		&\displaystyle \partial_t (b+h) + \frac{u \, \partial_x (b+h) - w}{\sqrt{1+\varepsilon^2 |\partial_x(b+h)|^2}}=0,\qquad \text{at } z=b+h, \label{eq_boundary_freesurface_1} \\
		&\displaystyle \partial_x (b+h) \left(\frac{1}{Fr^2}p - 2 \frac{\varepsilon }{Re} \mu \partial_x u \right) + \frac{1}{\varepsilon Re} \mu \left(\partial_z u + \varepsilon^2 \partial_x w \right) = 0, \qquad \text{at } z=b+h,
		\label{eq_boundary_freesurface_2}\\
		&\displaystyle \frac{1}{Fr^2} p + \frac{\varepsilon }{Re} \partial_x (b+h) \mu \left(\varepsilon^2 \partial_x w + \partial_z u \right) - 2 \frac{\varepsilon }{Re} \mu \partial_z w = 0, \qquad \text{at } z=b+h,
		\label{eq_boundary_freesurface_3}
	\end{empheq}
	\label{eq_boundary_freesurface}
\end{subequations}
and
\begin{subequations}
	\begin{eqnarray}
		&\displaystyle u \, \partial_x b = w,\qquad \text{at } z=b, \label{eq_boundary_bottomsurface_1}\\
		&\displaystyle -\frac{1}{Re}\mu (\partial_z u+\varepsilon^2\partial_x w)(1-\varepsilon^2|\partial_x b|^2)+2\frac{\varepsilon^2}{Re}\mu (\partial_x u+\partial_z w)\partial_x b=-\frac{1}{Re}k_b u , \quad \text{at } z=b. \label{eq_boundary_bottomsurface_2}
	\end{eqnarray}
	\label{eq_boundary_bottomsurface}
\end{subequations}

\paragraph{Depth averaging process.}
Now we integrate the system \eqref{eq_2d_epsilon} over the height $h$ using boundary conditions \eqref{eq_boundary_freesurface}-\eqref{eq_boundary_bottomsurface} to obtain a depth averaged model at the main order in $\varepsilon$. 
We introduce the depth-averaged velocity $\bar u=\frac{1}{h}\int_b^{b+h} u \, dz$.\\

Using the kinematic condition at the surface \eqref{eq_boundary_freesurface_1} and the no-penetration condition at the bottom \eqref{eq_boundary_bottomsurface_1}, we obtain the mass conservation equation:
\begin{equation}
\partial_t h + \partial_x (h \bar u) =\mathcal{O}(\varepsilon).
\end{equation}
We assume hydrostatic pressure hypothesis by neglecting the convective terms in the vertical momentum equation in \eqref{eq_2d_epsilon}. 
Then, by integration of this equation in $[z,b+h]$ and using boundary condition \eqref{eq_boundary_freesurface_3}, the pressure reads:
\begin{equation}
p=(b+h-z)\cos\theta+\mathcal{O}(\varepsilon).
\end{equation}
We embed this expression in the horizontal momentum equation and perform the vertical integration in $[b,b+h]$ using boundary conditions \eqref{eq_boundary_freesurface_1}, \eqref{eq_boundary_freesurface_2}, \eqref{eq_boundary_bottomsurface_1} and \eqref{eq_boundary_bottomsurface_2}. We obtain
\begin{equation}
\partial_t \left(h \bar u\right)  + \partial_x \left(h \bar u^2 \right) +\frac{1}{Fr^2}h\cos\theta  \partial_x (b+h)=\frac{1}{\varepsilon}\frac{1}{Fr^2}h\sin\theta -\frac{1}{\varepsilon Re}k_b \bar u+\mathcal{O}(\varepsilon).
\end{equation}
Note that the viscosity terms have been neglected and the viscosity effect only appears at the bottom through the friction term. Thus, we consider that the coefficient $k_b$ is an approximation of the viscous stress at the bottom in equation \eqref{mu_adim} (see \cite{costa_macedonio_2005, kelfoun:2015,biagioli:2023, bouchut:2025}). 
We analyze the boundary condition
\eqref{eq_boundary_bottomsurface_2}, which at the main order gives
$$
(\mu \partial_z u)_{|b}=k_b \bar u+\mathcal{O}(\varepsilon^2)
$$
where we identified $u_{|b}=\bar u$.
To set a value of $k_b$ in terms of the unknowns of the system, we use that $|D({\bf U})|=|\partial_z u|+\mathcal{O}(\varepsilon)$, and use the approximation 
\begin{equation}\label{dzu_aprox}
(\partial_z u)_{|b} \sim s_b\bar u/h,
\end{equation} 
for a constant coefficient $s_b$ (for granular flows in viscous regime, $s_b=3$ in \cite{cassar:2005}). Then, 
$$
k_b=\mu_b \frac{s_b}{h},
$$
where $\mu_b$ is an approximation of the apparent viscosity at the bottom.
Following \eqref{mu_adim} we introduce:
\begin{equation}\label{mu_approx}
    \mu_b=a_f\left(\left|s_b \frac{\bar u}{h}\right|^{n-1}+Bi\, \frac{1}{|s_b\frac{\bar u}{h}|}\right),
\end{equation}
where $a_f$ is a friction coefficient (see Supplementary Information \ref{app:friction}).
The dimensionless system then reads
\begin{subequations}
	\begin{empheq}[left=\empheqlbrace]{align}
		&\partial_t h + \partial_x (h \bar u) =0, \label{eq_syst1D_notclosed_h}\\
		&\displaystyle \partial_t \left(h \bar u\right)  + \partial_x \left(h \bar u^2 \right) + \frac{h \cos\theta}{Fr^2} \, \left( \partial_x (b+h) - \frac{\tan \theta}{\varepsilon} \right)
		= - \frac{a_f}{\varepsilon Re}\left( \left|s_b\frac{\bar u}{h} \right|^{n} + Bi \right) \frac{\bar u}{|\bar u|}. \label{eq_syst1D_notclosed_u}
	\end{empheq}
	\label{eq_syst1D_notclosed}
\end{subequations}
\\\\
To close the system we need to integrate the equation \eqref{eq_phi_p_adv_bis} for Model 1. As $\tau_{\rm crys}$ and $\Phi_{\rm crys}^*$ are assumed constant, the integrated equation simply reads
\begin{align}\label{eq_phi_p_model1}
	\partial_t (h\bar\Phi_p^*) + \partial_x \left(h\bar u \bar\Phi_p^*\right) = - \frac{h}{\tau_{\rm crys}} \left(\bar\Phi_p^*-\Phi_{\rm crys}^* \right),
\end{align}
where $\bar\Phi_p^*=\frac{1}{h}\int_b^{b+h} \Phi_p^* dz$. 
For Model 2, the crystal fraction is given in equation \eqref{eq_phip_T_lin}. The temperature equation \eqref{eq_temperature_costa_adap} is adopted for the one-dimensional case, and, in particular, we use the approximation in \eqref{dzu_aprox} to give the definition of the viscous heating term $Q_{visc}$ in \eqref{Qvisc_integral2} as
\begin{equation}\label{Qvisc_integral_fin}
Q_{visc}= \frac{s_b^2\eta }{\rho c_p h}  |\bar u|^2.  
\end{equation}
In other words, it is necessary to approximate $\partial_z u$ within the layer in \eqref{Qvisc_integral2}. Here we use as a simplification the same coefficient $s_b$ that is related to the shear rate at the flow base.

\subsection{Final proposed model}\label{sec:final_model}

We write the obtained depth-averaged system in dimensional form, dropping the bar notation for clarity.
\begin{subequations}\label{eq_syst1D_closed}
	\begin{empheq}[left=\empheqlbrace]{align}
		&\partial_t h + \partial_x (h u) =0, \label{eq_syst1D_closed_h}\\
		&\displaystyle \partial_t \left(h u\right)  + \partial_x \left(h u^2 \right) + g h \cos\theta \, \left( \partial_x (b+h) - \tan \theta \right) 
		= - \frac{a_f}{\rho} \left(\eta \left|s_b\frac{u}{h} \right|^{n} + \tau_c\right) \frac{u}{|u|}. \label{eq_syst1D_closed_u}
	\end{empheq}
\end{subequations}
This system is closed by setting the rheology with formulas \eqref{eta_prop}, \eqref{tauc_prop} and \eqref{n_prop} for $\eta$, $\tau_c$, $n$: $\forall \, 0 \leq \Phi_p^* < \Phi_m^*$,
\begin{subequations}\label{eq_rheology_final}
		\begin{equation} \displaystyle \eta = \eta_0 \left(1-\Phi_b\right)^{-1} \left(1-\frac{\Phi_p}{\Phi_m^*} \right)_+^{-2} \label{eq_rheology_final_eta}
		\end{equation}
		\begin{equation}
		\displaystyle \tau_c = \min \left(\tau_c^{\max}, \tau_c^* \left(\left(1-\frac{\Phi_p^*}{\Phi_m^*}\right)_+^{-2} -1 \right) \right)  \label{eq_rheology_final_tauc}
		\end{equation}
		\begin{equation}
		\displaystyle n = 1- \alpha_n \left(\frac{(\Phi_p-\Phi_p^c)_+}{\Phi_m^*}\right)^{4} - 0.334 \Phi_b \label{eq_rheology_final_n}
		\end{equation}
\end{subequations}
Note that thanks to the relation \eqref{phibp_relation} we write $\Phi_p$ in terms of $\Phi_b$ and $\Phi_p^*$ as 
$$
\Phi_p=\frac{\Phi_p^*(1-\Phi_b)}{1-\Phi_p^*\Phi_b}.
$$
For three-phase suspensions $\Phi_p<\Phi_p^*$ is satisfied because $0< \Phi_b< 1$. As a consequence, in the case $\Phi_p^*< \Phi_m^*$ the viscosity remains bounded.
In the limit case of solid suspensions, $\Phi_b=0$ and $\Phi_p=\Phi_p^*$, then both the yield stress and the viscosity reach their maximum respective values at the same time and the flow stops.\\

In the case $\Phi_m^*\leq \Phi_p^*\leq \Phi_{\rm crys}^*$ we have $\tau_c=\tau_c^{\max}$, where $\tau_c^{\max}\sim 10^9-10^{10}$ Pa, which is set to be large enough to obtain a stationary solution at rest. We thus obtain that $u=0$, $\partial_t h=0$ and it is not necessary to solve the rheology in \eqref{eq_rheology_final}.

\paragraph{Model $1$}
The crystal fraction $\Phi_p^*$ evolves as:
\begin{align}
	\partial_t (h \Phi_p^*) + \partial_x \left(hu \, \Phi_p^*\right) = - \frac{h}{\tau_{\rm crys}} \left(\Phi_p^*-\Phi_{\rm crys}^* \right),
	\label{eq_phi_p_final}
\end{align}
for $\tau_{crys}>0$  and $0< \Phi_{\rm crys}^*\leq 1$ given constants.

\paragraph{Model $2$}
The crystal fraction $\Phi_p^*$ depends on the lava temperature and their evolution are given by:
\begin{subequations}\label{model2_final}
    \begin{align}
	\left\{ \begin{array}{l}
		\displaystyle \Phi_p^* =\Phi_{\rm crys}^*\, \min \left( 1, \frac{(T_l-T)_+}{T_l-T_s} \right), 
		\\[4mm]
		\displaystyle \partial_t (h T)+ \partial_x (hu \,T) 
		=  C_{rad} (T_{env}^4-T^4) + C_{conv} (T_{env}-T) + \frac{C_{cond}}{h} (T_{c}-T) + \frac{C_{visc} \, \eta}{h} u^2.
	\end{array} \right.
	\label{eq_rheology_phi_version2}
\end{align}
The coefficients are
\begin{align}
	C_{rad}= \frac{\epsilon \, \sigma \, f}{\rho c_p},\quad
	C_{conv} = \frac{\lambda f}{\rho c_p},\quad
	C_{cond} = \frac{n_0 k}{\rho c_p } ,\quad
	C_{visc} = \frac{s_b^2}{\rho c_p },
\end{align}
\end{subequations}
where $T_{env}$ (K) is the external temperature, $T_c$ (K) the temperature of the ground, $\rho$ (kg.m$^{-3}$) is the lava density, $c_p$ (J.kg$^{-1}$.K$^{-1}$) is the heat capacity of lava, $\epsilon$ (dimensionless) is the emissivity coefficient, $f$ (dimensionless) is the fractional area of exposed lava to the air, $\sigma = 5.67 \times 10^{-8}$ W.m$^{-2}$.K$^{-4}$ is the Stefan-Boltzman constant, $\lambda$ (W.m$^{-2}$.K$^{-1}$) is the atmospheric heat transfer coefficient, $k$ (W.m$^{-1}$.K$^{-1}$) is the thermal conductivity and $n_0,s_b$ are dimensionless coefficients.

\subsection{Stationary solutions}

For Models $1$ and $2$ a stationary solution of system \eqref{eq_syst1D_closed} verifies $\partial_x q = 0$, that is, the discharge $q = hu$ is constant in space and time. We now establish stationary solutions for the uniform and non-uniform cases, defining $h$, $\Phi_p^*$ and $T$.

\paragraph{Uniform stationary solutions.}
\begin{figure}[hbtp]
	\centering
	\begin{subfigure}[b]{0.49\textwidth}
		\includegraphics[width=\textwidth]{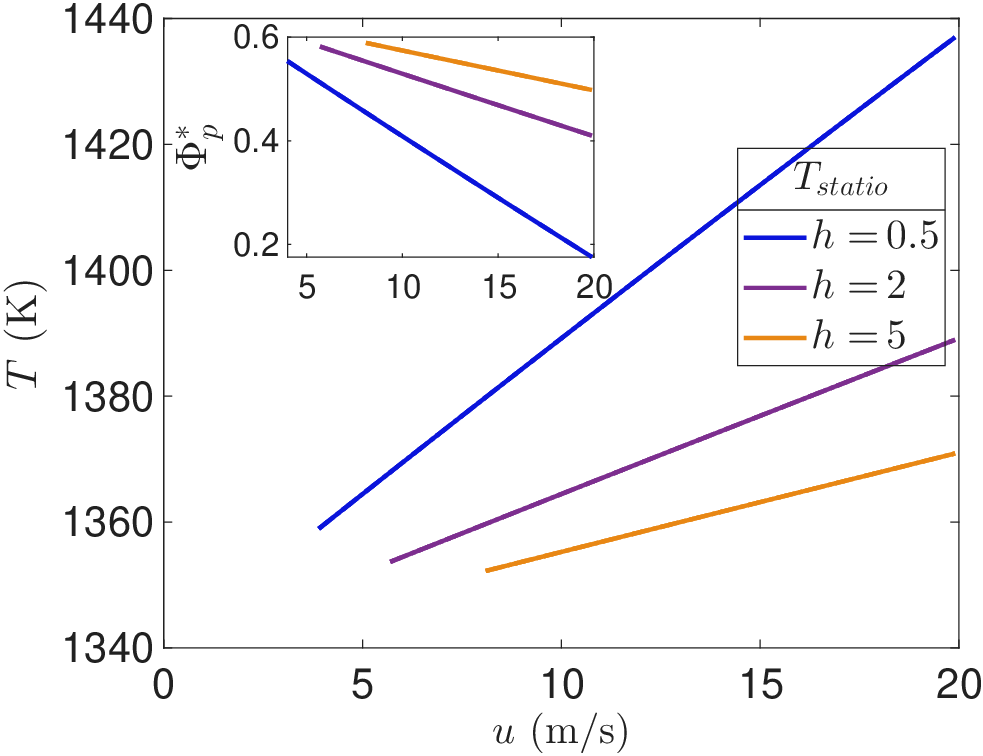}
		\caption{Stationary temperature $T$.}
		\label{fig_T_statio_T}
	\end{subfigure}\hspace{0mm}
	\begin{subfigure}[b]{0.49\textwidth}
		\includegraphics[width=\textwidth]{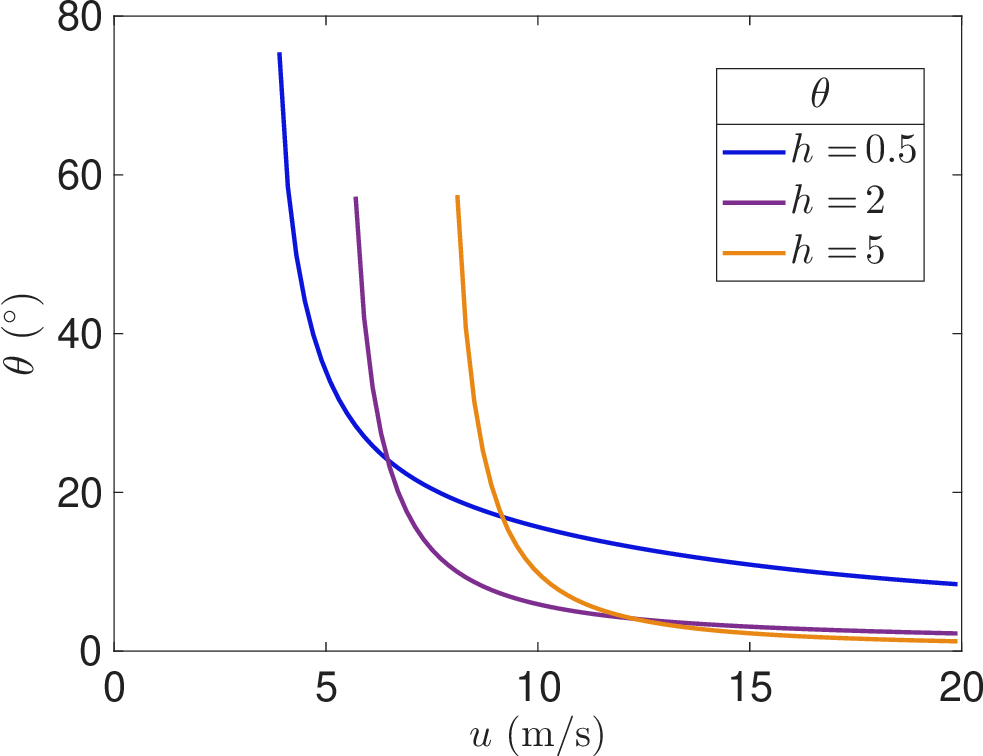}
		\caption{Minimal angle $\theta$.}
		\label{fig_T_statio_theta}
	\end{subfigure}\\[2mm]
	\begin{subfigure}[b]{0.49\textwidth}
		\includegraphics[width=\textwidth]{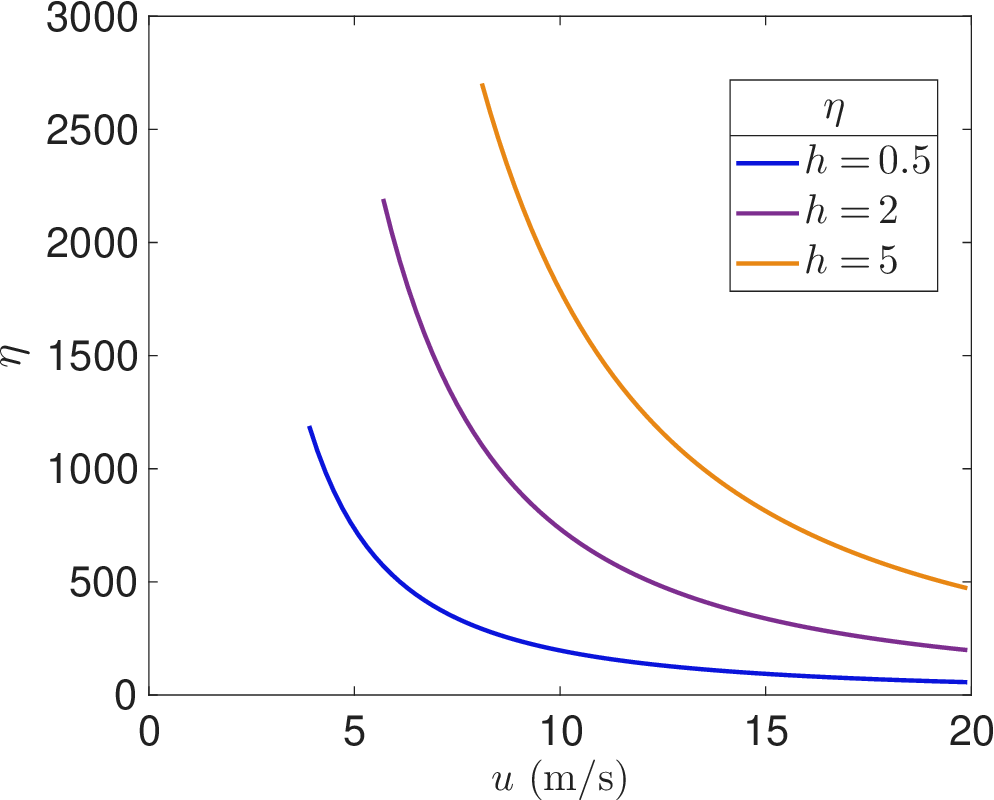}
		\caption{$\eta$.}
		\label{fig_T_statio_eta}
	\end{subfigure}
	\begin{subfigure}[b]{0.49\textwidth}
		\includegraphics[width=\textwidth]{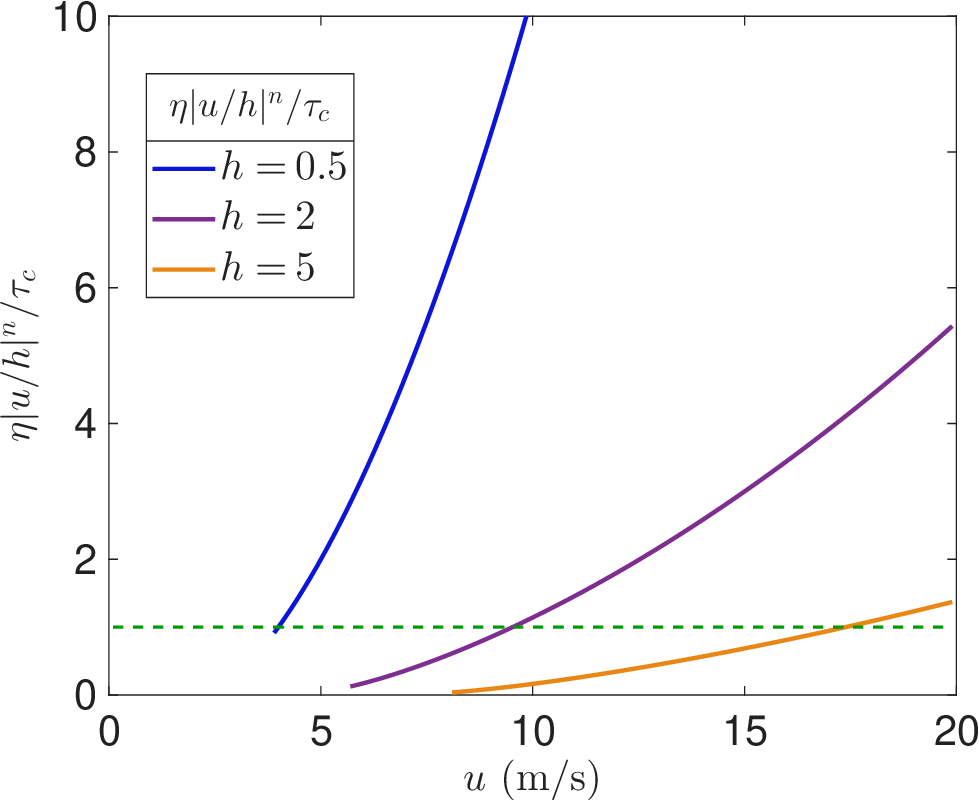}
		\caption{$\eta |u/h|^n/\tau_c$.}
		\label{fig_T_statio_etatauc}
	\end{subfigure}\hspace{2mm}
	
	\caption{Positive solutions of \eqref{eq_statio_unif_T} for different values of $u$, $h$. Parameters from table \ref{tab_param_test_scheme}.}
	\label{fig_T_statio}
\end{figure}
In the case of a uniform flow, with 
\begin{equation} \label{eq:def_zb}
z_b(x) = (L_x-x) \tan \theta + b(x),
\end{equation}
and with $h=h_0$ constant, for Model 1 there is only one possible definition:
$$
h(x,t)=h_0>0 \, \textrm{(constant)}, \quad q_0=0,
$$
because the uniform stationary solution of \eqref{eq_phi_p_final} is $\Phi_p^*=\Phi_{crys}^* \geq \Phi^*_m$. As a consequence $\tau_c=\tau_c^{\max}$, which we set to a large enough value to obtain stationary solutions at rest. In uniform stationary solutions, it is enough that $\rho g h_0 \sin \theta \leq \tau_c^{\max}$.

Model $2$ is, however, different because the evolution of $\Phi_p^*$ depends on that of the temperature. In the temperature evolution there is a mixing between radiative, convective, conductive and viscous effects, then it is possible to have cases where the velocity induces an increase of temperature that opposes the cooling by ambient temperature, thus allowing a constant lava flux. So for Model 2 we have $\Phi_p^*=\Phi_p^*(T)$, where $T$ is the solution of the following equation,
\begin{align}
C_{rad} (T_{env}^4-T^4) + C_{conv} (T_{env}-T) + \frac{C_{cond}}{h_0} (T_{c}-T) + \frac{C_{visc} \, \eta \left(\Phi_p(T)\right)}{h_0^3} q_0^2=0.
\label{eq_statio_unif_T}
\end{align}

If a solution of equation \eqref{eq_statio_unif_T} is such that $T\leq T_s$, by \eqref{eq_rheology_phi_version2} we obtain that $\Phi_p^*=\Phi_{crys}^*$ and $q_0=0$. If $T_s < T \leq T_l$, then the flux $q_0 = h u > 0$ if $\rho g h_0 \sin \theta \leq \tau_c(\Phi_p^*(T))$, and $q_0 = 0$ otherwise. 
In Figure \ref{fig_T_statio} the values of the positive stationary solutions (when they exists) are presented in terms of $u$ for several values of $h$.  In figure \ref{fig_T_statio_T} we can observe that a hotter temperature is reached for smaller height, and as a consequence smaller values of $\Phi_p^*$. We also observe that for a fixed value of $h$ the equilibrium temperature increases with velocity, which is due to viscous heating. In Figure \ref{fig_T_statio_theta} we can observe that smaller velocities require steeper slope angle so that a stationary solution exists. We can also observe a smaller viscosity when height or velocity increases in Figure \ref{fig_T_statio_eta}. Finally, Figure \ref{fig_T_statio_etatauc} displays the ratio between viscous and friction terms, which shows that friction effects are smaller than viscous ones when height or velocity increases.

\paragraph{Non-uniform stationary solutions.} In the general non-uniform case, solutions of \eqref{eq_syst1D_closed} verify, for $q_0 \in \mathbb{R}$:
\begin{align}
	\left\{ \begin{array}{ll}
		\displaystyle |g h \cos \theta \, \partial_x (z_b+h) | \leq \frac{\tau_c(\Phi_p^*)}{\rho} & \text{ if } q_0 = 0,\\
		\displaystyle - \frac{q_0^2 \partial_x h}{h^2} + g h \cos \theta \, \partial_x (z_b+h)
		= - \frac{a_f}{\rho} \left(\eta(\Phi_p(\Phi_p^*)) s_b^n\frac{q_0^n}{h^{2n}} + \tau_c(\Phi_p^*)\right) \frac{q_0}{|q_0|} & \text{ if } q_0 \neq 0.
	\end{array} \right.
    \label{eq_sol_statio_h}
\end{align}

For Model $1$ the crystal fraction $\Phi_p^*$ verifies a linear differential equation, the solution of which is given by:
\begin{align}
    \Phi_p^*(x) = 
    \left\{ \begin{array}{ll}
    	\displaystyle \, \Phi_{crys}^* & \text{ if } q_0 = 0,\\
    	\displaystyle \left(\Phi_p^*(0)-\Phi_{crys}^*\right) \exp \left(-\frac{\int_0^x h}{q_0 \tau_{crys}}\right) + \Phi_{crys}^* & \text{ if } q_0 \neq 0.
    \end{array} \right.
    \label{eq_sol_statio_phi} 
\end{align}

For Model $2$ we are not able to obtain an explicit expression of $\Phi_p^*$, because $\Phi_p^*=\Phi_p^*(T)$ (see equation \eqref{eq_rheology_phi_version2}), so $T$ is a solution of the following non-linear, initial value problem  
\begin{align}
    \left\{ \begin{array}{l}
    	\displaystyle q_0 \partial_x T  = C_{rad} (T_{env}^4-T^4) + C_{conv} (T_{env}-T) + \frac{C_{cond}}{h} (T_{c}-T) + \frac{C_{visc} \, \eta \left(\Phi_p^* \right)}{h^3} q_0^2, \\
    	T(0)=T_0.
    \end{array} \right.
    \label{eq_sol_statio_T}
\end{align}
In summary, semi-analytical stationary solutions can be computed solving equations \eqref{eq_sol_statio_h}-\eqref{eq_sol_statio_phi} for Model $1$, and equations \eqref{eq_sol_statio_h}-\eqref{eq_sol_statio_T} for Model $2$. See numerical results in Supplementary Information \ref{sec:test2}.

\section{Numerical tests}

A specific well-balanced numerical scheme is proposed, defined in terms of a hydrostatic reconstruction and a semi-implicit discretization.  Moreover, for the case of Model 2 a semi-implicit semi-linearized discretization is considered to approximate radiation terms. See Appendix \ref{app:num_scheme} for a complete description.

To show the convergence of the scheme, we carried out a comparison between the semi-analytical solutions and solutions of the general scheme given in section \ref{sec_model_num}. The details of the numerical tests are in Supplementary Information \ref{app:num}. All the numerical solutions behave well in the sense that there are no spurious oscillations caused by instabilities of the scheme, and the solutions converge to stationary states.\\

Briefly, the first series of verification tests of our numerical scheme was to compare uniform solutions of Models $1$ and $2$. Results show that, as expected, viscous effects are greater than the friction effects initially and this ratio decreases with time because of the deceleration. The second test series aimed at verifying that numerical solutions converge towards a known stationary solution. In this series, the stationary solution $h$ is a solution of the ODE \eqref{eq_sol_statio_h} was compared with the solution computed by the scheme with three different space steps. The third test series aimed at verifying that numerical solutions of both Models $1$ and $2$ converge towards unknown stationary solutions. Results show that the  differences  between  the  two  solutions  are  caused  by  the  friction, as the solid volume fraction $\Phi_p$ evolves differently in each model. In Model $1$, $\Phi_p$ is always increasing because the source term in \eqref{eq_phi_p_final} is positive (as $\Phi_p<1$). On the other hand, in Model 2, the sign of the source term in \eqref{eq_rheology_phi_version2} can change. The predominance of this viscous heating term for small $x$ values causes a positive total heat loss. Finally, the last series of tests aimed at checking the scheme accuracy in capturing dam break with wet and dry fronts. In both cases, the lava flows down the slope and its velocity decreases as expected when the crystal fraction increases. The lava velocity vanishes completely when the crystal fraction reaches the critical level. \\

\section{Model intercomparison: application to Mauna Ulu lava flow}

We use Model $2$ to reproduce a lava flow on Mauna Ulu in 1969–1974, on the Kīlauea volcano (Hawaï). Data from this eruption can be found in \cite{robert2014textural,harris2022anatomy}, where they study the morphology of this lava flow. The model PyFLOWGO is applied in \cite{harris2022anatomy} to compare numerical solutions to the data, and all the parameters that they use are given in the supplementary material of their paper. PyFLOWGO is an ideal candidate for a comparison because it is a 1D model applied to a channeled lava flow. Its approach is kinematic; it consists in a single differential equation that describes temperature evolution:
\begin{align}
    \partial_x T = \frac{Q_{rad}^{Harris}+Q_{conv}^{Harris}+Q_{cond}^{Harris}+Q_{rain}^{Harris}+Q_{visc}^{Harris}}{q \rho L_{cryst} \partial \phi / \partial T},
    \label{eq_T_pyflowgo}
\end{align}
where $L_{cryst}$ is the specific heat of crystallization.  
Our equation of $T$ is equivalent to this equation \eqref{eq_T_pyflowgo} if we suppose that the solution is stationary, that $q$ is constant, and that $c_p = L_{cryst} \partial \phi / \partial T$. The heat flux terms $Q^{Harris}$, in $W/m$, are given in the supplementary material of their paper, and they assume that $Q_{visc}^{Harris}=0$.\\

Unlike our model, equilibrium crystallization is calculated a posteriori in PyFLOWGO. The bubble-free crystal fraction $\phi$ is deduced from $T$:
\begin{align*}
    \partial_x \phi = \frac{\partial \phi}{\partial T} \frac{\partial T}{\partial x}.
\end{align*}
This led us to choose Model $2$ to carry out the comparison as it considers equilibrium crystallization implicitly, which is closer to the PyFLOWGO approach than the explicit crystallization of Model $1$. The quantity $\partial \phi / \partial T$ is given by a bimodal model (see \cite{robert2014textural}), the crystallization rate change after a given distance $x_{crit}$:
\begin{align*}
    \partial \phi / \partial T = 
    \left\{ \begin{array}{l}
        C_1 \quad \text{if } x \leq x_{crit}\\
        C_2 \quad \text{if } x > x_{crit}.
    \end{array} \right.
\end{align*}
It would have been possible to adapt the equation on $\Phi_p$ for Model $2$, where $\partial \phi / \partial T = -1/(T_l-T_s)$ is also a constant, replacing it by this bimodal model. Such a modification was not attempted to highlight the consequences of using an independent equilibrium crystallization path. \\

The kinematic approach of PyFLOWGO means that flow velocity is also calculated a posteriori. The mean velocity is computed with:
\begin{align*}
    V_{mean} = \frac{\tau_b d}{n_{shape} \eta^{Harris}} \left(1 - \frac{3}{2} \frac{\tau_c^{Harris}}{\tau_b} + \frac{1}{2} \left(\frac{\tau_c^{Harris}}{\tau_b}\right)^3 \right),
\end{align*}
where $n_{shape}= 3\left(1+\frac{h}{w}\right)^2$ is the channel shape factor, $\tau_b = \rho g d \sin \theta$ (Pa) is the basal shear stress and $d$ is the mean height of the channel. The viscosity formula supposes that crystals and bubbles are of the same size:
\begin{align}
    \eta^{Harris} = \eta_0 \left(1-\phi-\Phi_b \right)^{\frac{5\phi-2\Phi_b}{2\phi-\Phi_b}},
    \label{eq_viscosity_harris}
\end{align}
and the yield strength is given by
\begin{align}
    \tau_c^{Harris} = 0.01 \left(e^{0.08(T_0-T)}-1\right) + 6500 \phi^{2.85}.
    \label{eq_tauc_harris}
\end{align}
where $T_0$ is the temperature at the vent and $T$ is the lava temperature.

\subsection{Choice of parameters}

In \cite{harris2022anatomy} the input parameters used for running the simulation with PyFLOWGO are given in their Appendix A, Table A1. The velocity computed give them a transit time for lava to cross the domain of $12.8$ min. 
\\ \\
The parameters we use are given in Table \ref{tab_param_appl}. As our model differs from PyFLOWGO, a few parameters need to be discussed.
\begin{itemize}
    \item The conductivity coefficient $n_0$ in the term $Q_{cond}$ is the inverse of $\delta_T / h$, where $\delta_T$ is the thermal boundary layer between the lava and the bottom surface (see \cite{costa_macedonio_2005}). The fraction $\delta_T / h$ is approximated by a constant, set to $0.19$ in \cite{harris2022anatomy}.
    \item The viscous heat effect is neglected in \cite{harris2022anatomy}, here we take this effect into account in the term $Q_{visc}$ given in \eqref{Qvisc_integral_fin}. As discussed before, the relationship between the coefficient in \eqref{eq_heat_terms} of \cite{costa_macedonio_2005}  is $m=s_b^2$ and is set $m=12$, which corresponds to $s_b=3.4$. Considering that $s_b$ is a semi-empirical constant, here we follow the work of \cite{cassar:2005}, who proposes that $s_b=3$ in granular flows in the viscous regime.
    \item The thermal conductivity of lava $\kappa$ is given by $\kappa = \left(1.929-1.554 \Phi_b\right)^2$ (see \cite{harris2001flowgo}, \cite[p41]{peck1978cooling}).
    \item Comparing the temperature equations \eqref{eq_rheology_phi_version2} and \eqref{eq_T_pyflowgo} we can fix $c_p = -L_{heat} \, \partial_T \phi$ to obtain similar formulas. In the simulation we use $c_p = - L_{heat} \times 3.8 \times 10^{-3}$.
    \item The fluid viscosity $\eta_0$ depends on its temperature, and it is given by $\eta_0(T) = 10^{a+\frac{b}{T-c}}$ (with $T$ in Kelvin) in \cite{robert2014textural}. We chose the same coefficients as in \cite{harris2022anatomy}, which are $a=-4.7$, $b=5429.7$, and $c=595.5$.
    
    \item The bubble volume fraction was chosen to remain within the validity domain of \eqref{rheo_prop} ($\Phi_b<0.3$), which yields half the value chosen by \cite{harris2022anatomy}. 
    \item 
We modify the friction term in our model (right-hand side of \eqref{eq_syst1D_closed_u}) to mimic the lateral friction in the channel. For this aim, we follow \cite{wilson:1993} and we introduce the lateral friction together with the bottom friction in terms of the viscosity and the shape factor.  Thus, we set
$$
  a_f=3\left(1+\frac{2h}{w}\right)^2, 
$$ 
with $w$ the channel width and $h$ the lava thickness (see Supplementary Information \ref{app:friction} for more details). We set the value of $w$ as a constant based on the average of the real data.

\end{itemize}

\begin{table}[hbtp]
	\center 
	\begin{tabular}{| l| l | l |}
		\hline
		$g$ & $9.81$ m.s$^{-2}$ & gravity\\
        $\theta$ & $0^\circ$ & angle of the reference plane\\
		$\Phi_b$ & $ 0.2433$ & volume fraction of bubbles \\
		$\Phi_m^*$ & $0.6$ & maximum volume fraction of crystals \\
		$\Phi_p^c$ & $0.14$ & critical value for yield stress \\
        $\alpha_n$ & $1.32$ & coefficient of the flow index equation\\
		$\rho$ & $2900$ kg.m$^{-3}$ & fluid density $^a$\\
		$\eta_0$ & $26$ Pa.s$^n$ & consistency factor\\
		$\tau_c^*$ & $33$ Pa & characteristic critical yield stress\\
		$\tau_c^{\max}$ & $1.2 \times 10^{11}$ Pa & maximal yield stress\\
		$T_l$ & $1473$ K & liquidus temperature \\
		$T_s$ & $1268$ K & solidus temperature $^a$\\
		$T_{env}$ & $293$ K & external temperature $^a$\\
		$\tau_{crys}$ & $80000$ s & characteristic crystallization time \\
        $\Phi_{crys}^*$ & $1$ & maximum volume percent crystallization\\
		$\varepsilon$ & $0.98$ & emissivity $^a$\\
        $\sigma$ & $5.67 \times 10^{-8}$ W.m$^{-2}$.K$^{-1}$ & Stephan Boltzmann constant $^a$\\
        $f$ & $0.5$ & fraction area of exposed core $^a$\\
        $L_{heat}$ & $350000$ J/kg & latent heat of crystallization $^a$\\
        $c_p$ & $- L_{heat} \times 3.8 \times 10^{-3}$ J.kg$^{-1}$.K$^{-1}$ & heat capacity of the lava\\ 
        $\lambda$ & $70$ W.m$^{-2}$.s$^{-1}$ & atmospheric heat transfer coefficient $^b$\\
        $n_0$ & $100/19$ & coefficient for thermal conductivity $^a$\\
        $s_b$ & $3$ & coefficient for viscous heat and friction\\
        $k$ & $1.4$ W.m$^{-1}$.K$^{-1}$ & thermal conductivity of the lava\\
	    $w_0$ & $11.5$ m & averaged channel width\\
        \hline  
	\end{tabular}
	\caption{Model parameters for the application. $^a$Values of \cite{harris2022anatomy}, $^b$ values of \cite{costa_macedonio_2005}.}
	\label{tab_param_appl}
\end{table}

For clarity, we write our proposed model including the additional parameters for this application:
\begin{subequations}
\begin{equation}
\left\{\begin{array}{l}
\partial_t h + \partial_x (h u) =0,         \\
\displaystyle \partial_t \left(h u\right)  + \partial_x \left(h u^2 \right) + g h \cos\theta \, \left( \partial_x (b+h) - \tan \theta \right) 
		= -\frac{3}{\rho}\left(1+\frac{2h}{w_0}\right)^2\left(\eta \left|s_b\frac{u}{h} \right|^{n} +\tau_c\right)\frac{u}{|u|}. 
	\end{array}\right.
	\label{eq_syst1D_application}
\end{equation}
The rheology functions are defined $\forall \, 0 \leq \Phi_p^* < \Phi_m^*$ and assume that ${\rm Cx}\ll 1$,
		\begin{equation} \displaystyle \eta = \eta_0 \left(1-\Phi_b\right)^{-1} \left(1-\frac{\Phi_p}{\Phi_m^*} \right)_+^{-2} 
		\end{equation}
		\begin{equation}
		\displaystyle \tau_c = \min \left(\tau_c^{\max}, \tau_c^* \left(\left(1-\frac{\Phi_p^*}{\Phi_m^*}\right)_+^{-2} -1 \right) \right)  
		\end{equation}
		\begin{equation}
		\displaystyle n = 1- \alpha_n \left(\frac{(\Phi_p-\Phi_p^c)_+}{\Phi_m^*}\right)^{4} - 0.334 \Phi_b 
		\end{equation}
	\label{eq_rheology_final_app}
\end{subequations}
The crystal volume fraction is calculated following two different models. \\
{\bf Model $1$.}
The crystal fraction $\Phi_p^*$ evolves as:
\begin{align}
	\partial_t (h \Phi_p^*) + \partial_x \left(hu \, \Phi_p^*\right) = - \frac{h}{\tau_{\rm crys}} \left(\Phi_p^*-\Phi_{\rm crys}^* \right).
\end{align}
\\
{\bf Model $2$.}
The crystal fraction $\Phi_p^*$ depends on the lava temperature and their evolution are given by:
\begin{subequations}
    \begin{align}
	\left\{ \begin{array}{l}
		\displaystyle \Phi_p^* =\Phi_{\rm crys}^*\, \min \left( 1, \frac{(T_l-T)_+}{T_l-T_s} \right), 
		\\[4mm]
		\displaystyle \partial_t (h T)+ \partial_x (hu \,T) 
		=  C_{rad} (T_{env}^4-T^4) + C_{conv} (T_{env}-T) + \frac{C_{cond}}{h} (T_{c}-T) + \frac{C_{visc} \, \eta}{h} u^2.
	\end{array} \right.
	\label{eq_rheology_phi_version2a}
\end{align}
The coefficients are
\begin{align}
	C_{rad}= \frac{\epsilon \, \sigma \, f}{\rho c_p},\quad
	C_{conv} = \frac{\lambda f}{\rho c_p},\quad
	C_{cond} = \frac{n_0 k}{\rho c_p } ,\quad
	C_{visc} = \frac{s_b^2}{\rho c_p }.
\end{align}
\end{subequations}

The boundary conditions are of type Dirichlet at the upstream boundary, and of type no flux at the downstream. The bottom surface is taken from \cite{harris2022anatomy}.  The initial conditions involve four free parameters. To keep this comparison simple, we report the results obtained by only varying one of them, the effusion rate at the vent, $E_0$. These parameters are:
\begin{align*}
	& T_0 = 1438.15 \text{ K (temperature of the lava at the vent, data reported in\cite{harris2022anatomy})}\\
	& h_0 = 2.475 \text{ m (lava depth at the vent, \cite{robert2014textural})}\\
	& w_0 = 11.5  \text{ m (lava width at the vent, average channel width \cite{harris2022anatomy})}\\
	& E_0 = 90-420 \text{ m$^3$/s (effusion rate at the vent, data reported in \cite{harris_2009})}
\end{align*}

\subsection{Simulations results for Model 2}

The results of Model $2$ are presented first as its outputs are used to reduce the number of degrees of freedom explored with Model $1$. One important difference between Model $2$ and PyFLOWGO is that Model $2$ describes the transient propagation of the flow front. This propagation does not lead to a steady state but to a quasi-steady state that can only be understood by first studying the transient evolution of the lava flow.

\subsubsection{Transient evolution}

\begin{figure}[hbtp]
	\centering
	\includegraphics[width=0.85\textwidth]{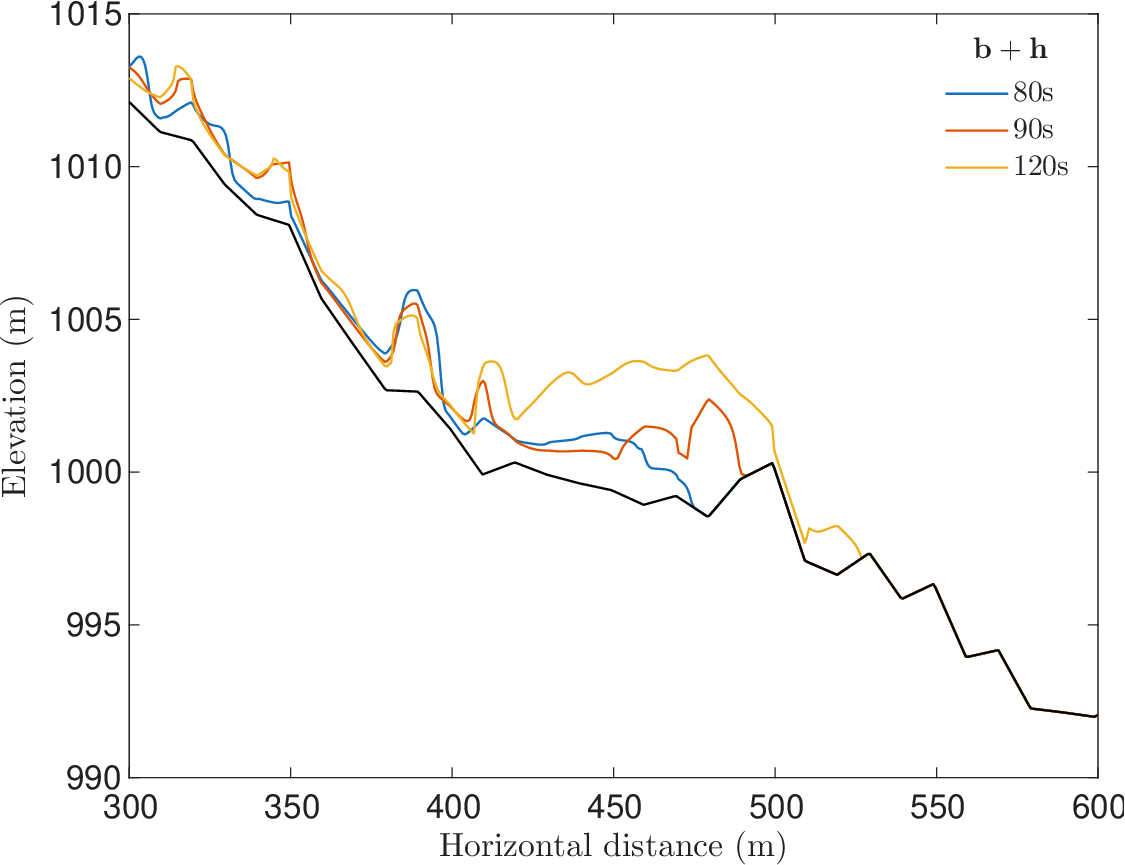}
	\caption{Model $2$ with $E_0=90$m$^3$/s.}
	\label{fig_harris_zoom1}
\end{figure}

Figure \ref{fig_harris_zoom1} shows the evolution of the lava thickness at three times in the proximal region with a vertical exaggeration (i.e. actual slopes are less steep than portrayed). At 90 s, the flow is blocked by a topographic high at $x=500$ m. The flow front ponds at that location for the next few seconds and, at $\approx 120$ s, it overcomes this small obstacle and a thinned front continues its travel downstream. Consistently with the results of the dam break with dry front, our model is thus able to accurately depict the effects of minute topographic changes.\\

\begin{figure}[hbtp]
	\centering
	\includegraphics[width=0.85\textwidth]{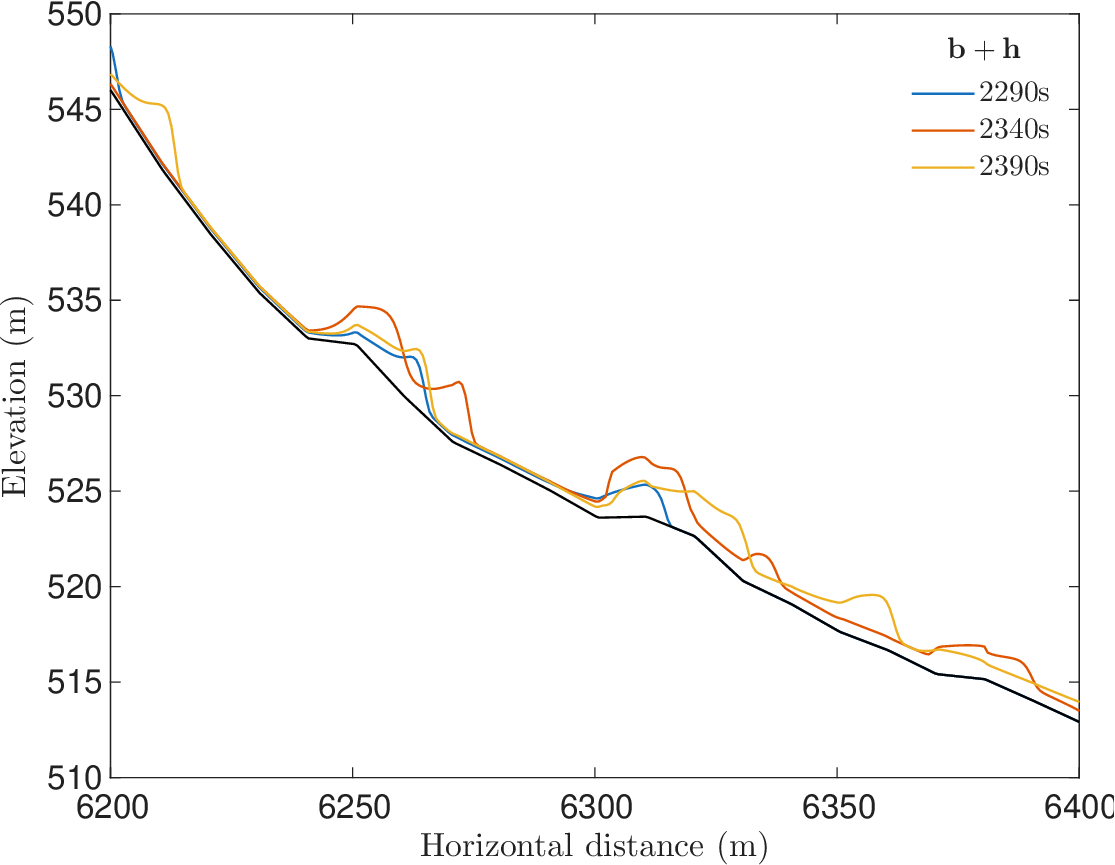}
	\caption{Model $2$ with $E_0=90$m$^3$/s.}
	\label{fig_harris_zoom2}
\end{figure}

Figure \ref{fig_harris_zoom2} shows the evolution of the lava thickness at three times in the distal region where the basal slope is the steepest (Fig. \ref{fig_harris_bh}). These three snapshots suggest a slightly complex picture. The flow seems to have a nearly consistently vanishing thickness in the smooth, steep part of the overall slope and thick lumps of lava clung to minor break-in-slope at $\approx 6250$ and $\approx 6320$ m. This distribution is caused by the rapid succession of the pond-then-breach behavior depicted proximally (Fig. \ref{fig_harris_zoom1}). The process is illustrated by a distinct lump of lava at 2390 s and $\approx 6200$ m that is traveling on the smooth, steep part of the slope. This lump will pond a few seconds later at the break-in-slope  around 6250 m. When repeated, this dynamics yields a pulsatory flow.\\

When the regional slope is weak and the lava is thick enough to create a continuous flow, such as proximally, the later evolution of the flow can reach a steady state. When, however, the regional slope is steep like in the distal region, the pulsatory nature of the flow precludes a steady state to be defined over a timescale comparable to the transit time of magma. We can thus expect, at best, a quasi-steady state to be established after that the flow front reaches the distal end of the computational domain.

\subsubsection{Quasi-steady state}

\begin{figure}[hbtp]
	\centering
	\includegraphics[width=0.75\textwidth]{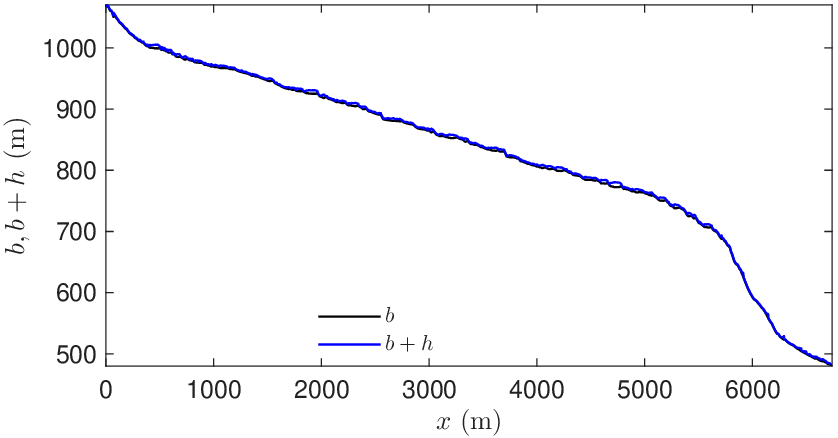}
	\caption{Model $2$ with $E_0=90$m$^3$/s, at time $3000$s. Bottom and free surface elevation.}
	\label{fig_harris_bh}
\end{figure}

\begin{figure}[hbtp]
	\centering
	\includegraphics[width=\textwidth]{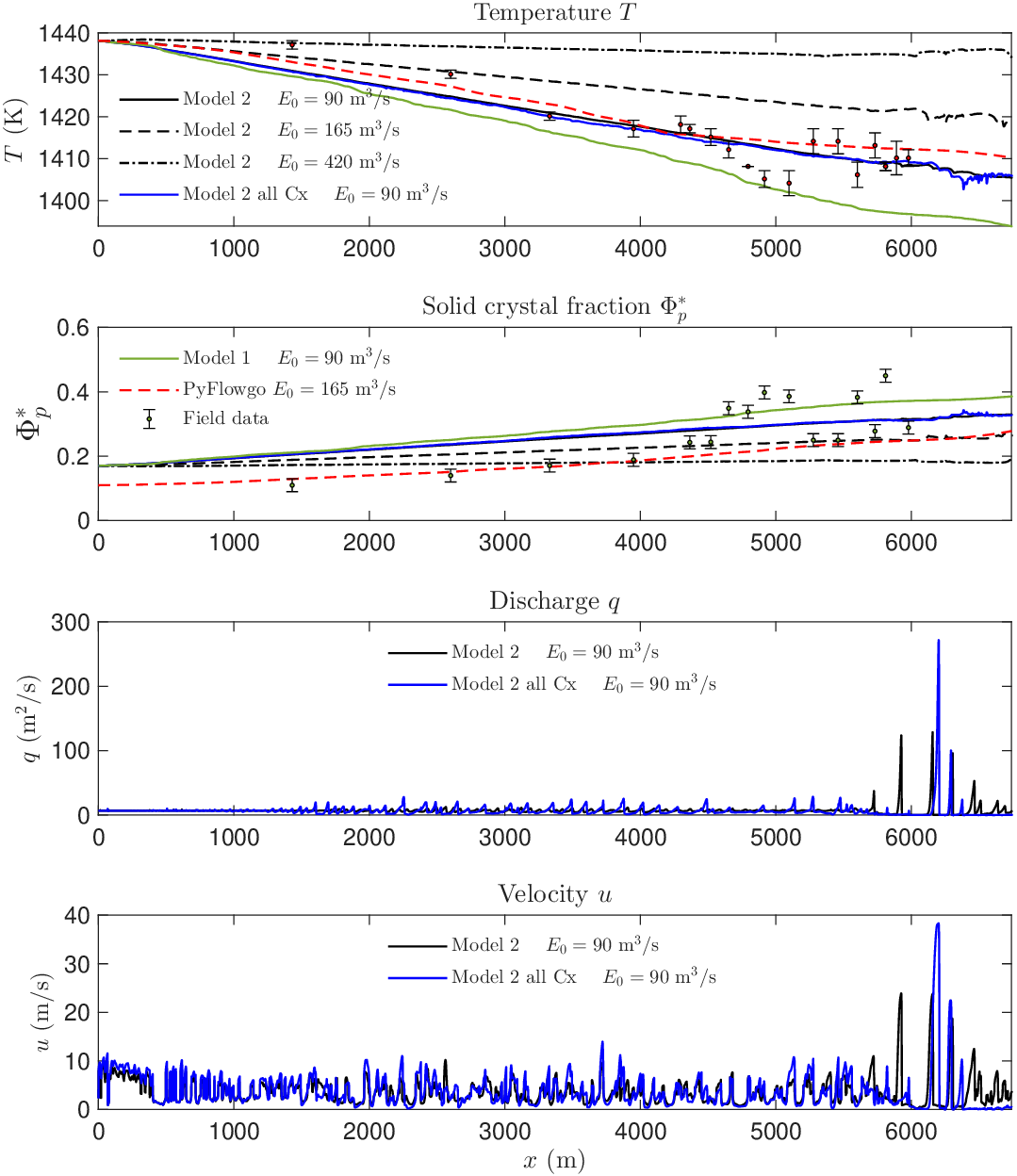}
	\caption{Evolution of temperature, crystal fraction, fluid discharge and velocity as a function of distance in quasi-steady state. Model $2$ runs were performed with $E_0=90$--$420$ m$^3$/s at $3000$ s. One Model $2$ run (all ${\rm Cx}$) was performed with the general viscosity \eqref{rheo_prop_gral}. The Model $1$ run was performed with $E_0=90$ m$^3$/s. The field data on $T$ and $\Phi^*_p$ are from \cite{harris2022anatomy}.}
	\label{fig_harris_u}
\end{figure}

The quasi-steady state simulations of Model $2$ with effusion rates in the range $90-420$ m$^3$/s are compared with the field data coming from the Mauna Ulu eruption. Figure \ref{fig_harris_bh} shows that the proximal region ($<5$ km), has a shallow regional slope, whereas the distal part of the domain ($>5$ km) is much steeper. We define the quasi-steady state as vanishing temporal variations of $h$ and $u$ proximally. Distally, the lava thickness, speed, and flux are highly irregular because of the pulsatory nature of the flow (Fig. \ref{fig_harris_zoom2}).\\

The following misfit on temperature is computed to choose the best fit:
\begin{align*}
    \text{err } = \frac{\sum \left(T_{data} - T_{simu}\right)^2}{\sum \left(T_{data} - T_{mean}\right)^2},
\end{align*}
$T_{mean}$ being the average computed temperature. The best fit is obtained for $E_0=90$ m$^3$/s, where err $=0.2278$ .\\

\begin{figure}[hbtp]
	\centering
	\includegraphics[width=\textwidth]{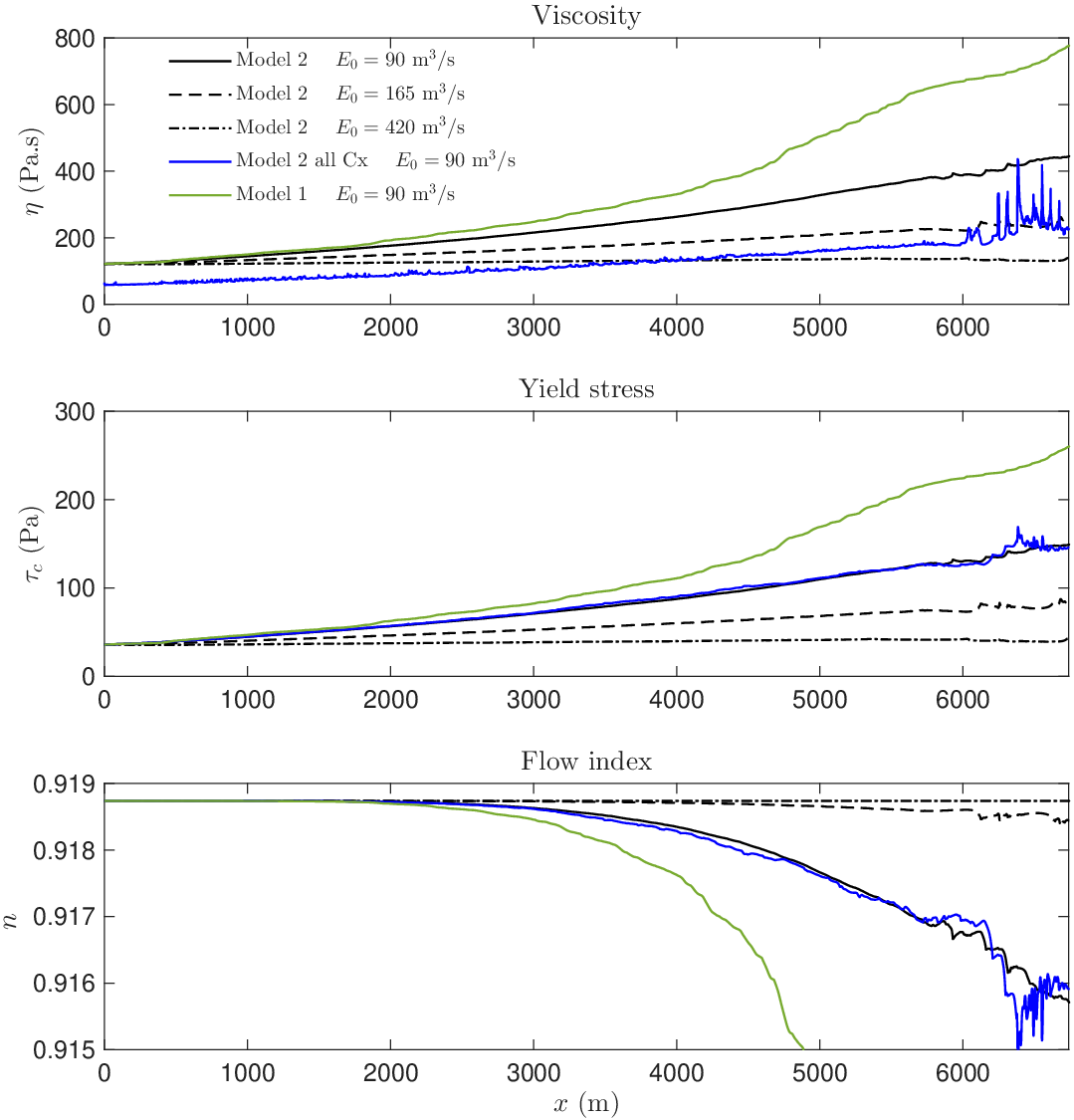}
	\caption{Evolution of rheological parameters ($\eta$, $\tau_c$, $n$) as a function of distance in quasi-steady state. Runs were performed with $E_0=90$--$420$ m$^3$/s at $3000$ s. One Model $2$ run (all ${\rm Cx}$) was performed with the general viscosity \eqref{rheo_prop_gral}. The Model $1$ run was performed with $E_0=90$ m$^3$/s.}
	\label{fig_harris_ntau}
\end{figure}

\begin{figure}[hbtp]
	\centering
	\includegraphics[width=\textwidth]{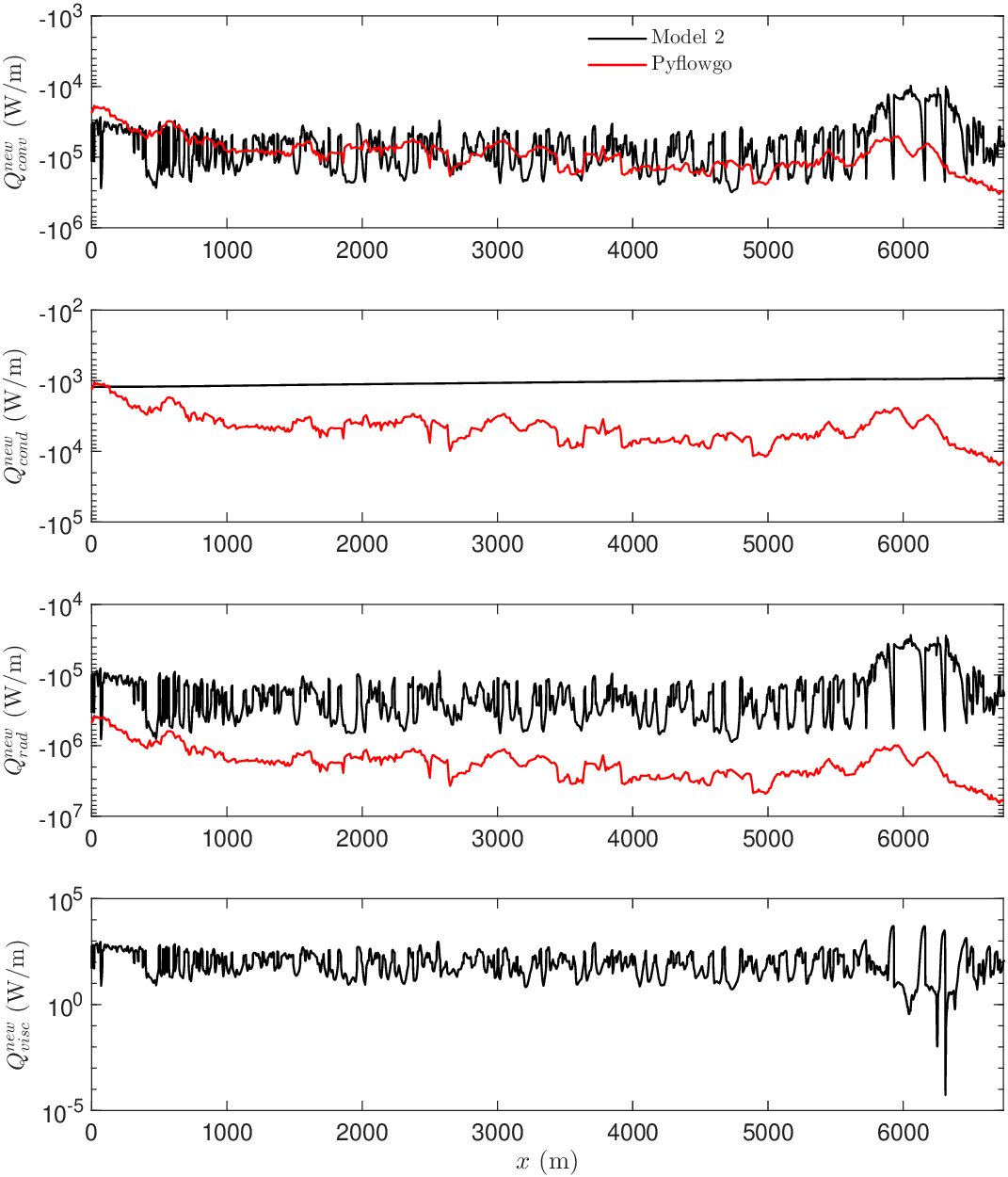}
	\caption{Heat exchanges terms as a function of distance at $4000$ s for Model $2$ with $E_0=90$ m$^3$/s and, when applicable, PyFLOWGO.}
	\label{fig_harris_model2_3}
\end{figure}

Figures \ref{fig_harris_u} and \ref{fig_harris_ntau} present model outputs for $E_0=90$, $165$, and $420$ m$^3$/s for $T$, $\Phi^*_p$, and $\eta$, whereas they present only the best fit output ($E_0=90$ m$^3$/s) for $q$, $u$, $\eta$, $\tau_c$ and $n$. The simulation times are chosen so that the lava flow is a quasi-steady state. The transit time for the best fit is $3290$ s, whereas it is $720$ s in the PyFLOWGO simulation.\\

Figure \ref{fig_harris_u} shows that the temperature is decreasing almost linearly with distance, regardless of the initial discharge rate. This cooling rate is fastest when the effusion rate is the smallest because the fluid velocity is then also the smallest. Conversely, the crystal fraction increases the most with distance at the smallest effusion rate. When compared with field data, the cool and high crystal fraction distal lava is best captured by the lowest effusion rates. The respective fits of the PyFLOWGO run for $T$ and $\phi$ are codependent because their crystallization rate was constrained by the temperature field data. As a result, the crystal fraction fit of PyFLOWGO is better than those of our model in the proximal region.\\

The fluid height, discharge rate, and velocity are represented only when $E_0=90$ m$^3$/s for clarity. The velocity range is $1-6$ m/s and both the fluid height and velocity vary because of the topography. An overall decrease of $h$ and an increase of $u$ is observed around $x=6$ km, because of the steeper slope there (see Fig. \ref{fig_harris_bh}) and the ensuing pulsatory flow that is clearly visible in the discharge rate output. The velocity range is $4-14$ m/s when $E_0=420$ m$^3$/s. At the same $E_0$ value as PyFLOWGO (165 m$^3$/s), the velocity range is $2-8$ m/s, which is much smaller than that of PyFLOWGO near the vent but similar to that of PyFLOWGO $>1$ km. \\

The viscosity, yields stress and flow index results are presented in Fig. \ref{fig_harris_ntau}. Both $\eta$ and $\tau_c$ increase with distance from the vent, and $n$ decreases with distance. Their rate of change is larger when the effusion rate is smaller. These results are similar to the ones computed with PyFLOWGO.\\ 

Figure \ref{fig_harris_model2_3} presents the evolution of heat terms when $E_0=90$ m$^3$/s, and a comparison with the results of PyFLOWGO simulations made in \cite{harris2022anatomy} (except for the viscous heating term that they neglect). These heat terms are given by the following formula, to match the units of $Q^{Harris}$:
\begin{align*}
    Q^{new} = \rho c_p h \, Q,
\end{align*}
where $Q$ are the heat terms given by equation \eqref{eq_heat_terms}. The two models have distinct formulas for the heat terms, so the results are quite different, especially for the conductive and radiative terms.\\

\subsection{Simulations results for Model 1}

We illustrate the difference of behavior between Models $1$ and $2$ by running Model $1$ with the same parameters as Model $2$ and $E_0=90$ m$^3$/s. The only additional free parameter that Model $1$ has compared to Model $2$ is the crystallization rate, $\tau_{\rm crys}$. To set its value, we  first define the crystal fraction ($\Phi_p^*$) of Model 1 as $\Phi_1$ to distinguish it from that of Model $2$, $\Phi_2$. We define similarly $u_1$, $u_2$, and $h_2$. In Model $1$, combining \eqref{eq_syst1D_closed_h} and \eqref{eq_phi_p_final} yields:
\begin{align}
	\partial_t \Phi_1 + u_1 \partial_x \Phi_1 = - \frac{\left(\Phi_1-\Phi_{\rm crys}^* \right)}{\tau_{\rm crys}} .
	\label{eq_phi_Lagr}
\end{align}

In Model $2$, calling $\mathcal{C}$ the right-hand side of the temperature equation in \eqref{eq_rheology_phi_version2a} and using \eqref{eq_syst1D_closed_h} yields:
\begin{align}
	\partial_t T + u_2 \partial_x T = \frac{\mathcal{C}}{h_2} .
	\label{eq_phi_Lagr2}
\end{align}
Using the definition \eqref{eq_rheology_phi_version2a} to substitute $T$ by $\Phi_2$ in \eqref{eq_phi_Lagr2} yields:
\begin{align}
	\partial_t \Phi_2 + u_2 \partial_x \Phi_2 = \frac{\Phi_{\rm crys}^* \mathcal{C}}{(T_l-T_s) h_2} ,
	\label{eq_phi_Lagr3}
\end{align}

Comparing equations \eqref{eq_phi_Lagr} and \eqref{eq_phi_Lagr3}, the evolution of crystal fraction in Model $1$ would thus roughly follow that of Model $2$ if
\begin{align}
	 \frac{\Phi_{\rm crys}^* \mathcal{C}}{(T_l-T_s) h_2} \approx   \frac{\left(\Phi_{\rm crys}^* -\Phi_1 \right)}{\tau_{\rm crys}}
	\label{eq_phi_Lagr4}
\end{align}
Simple equality between Models $1$ and $2$ using \eqref{eq_phi_Lagr4} is not possible as both $\mathcal{C}$ and $h_2$ on the left-hand-side are nonlinear functions of $t$ and $x$, whereas on the right-hand-side only $\Phi_1$ is a function of $t$ and $x$ \cite{crisp1994influence}. We can nevertheless define an equivalent crystallization rate, $\tau_{\rm equ}$, by assuming $\Phi_1=\Phi_2$ so that Model $2$ behaves as Model $1$ to the first order:
\begin{align}
	 \tau_{\rm equ} = \frac{h (T_l-T_s)(\Phi_2 -\Phi_{\rm crys}^* )}{\Phi_{\rm crys}^*\mathcal{C}}.
	\label{eq_tau_equ}
\end{align}

\begin{figure}[hbtp]
	\centering
	\includegraphics[width=0.9\textwidth]{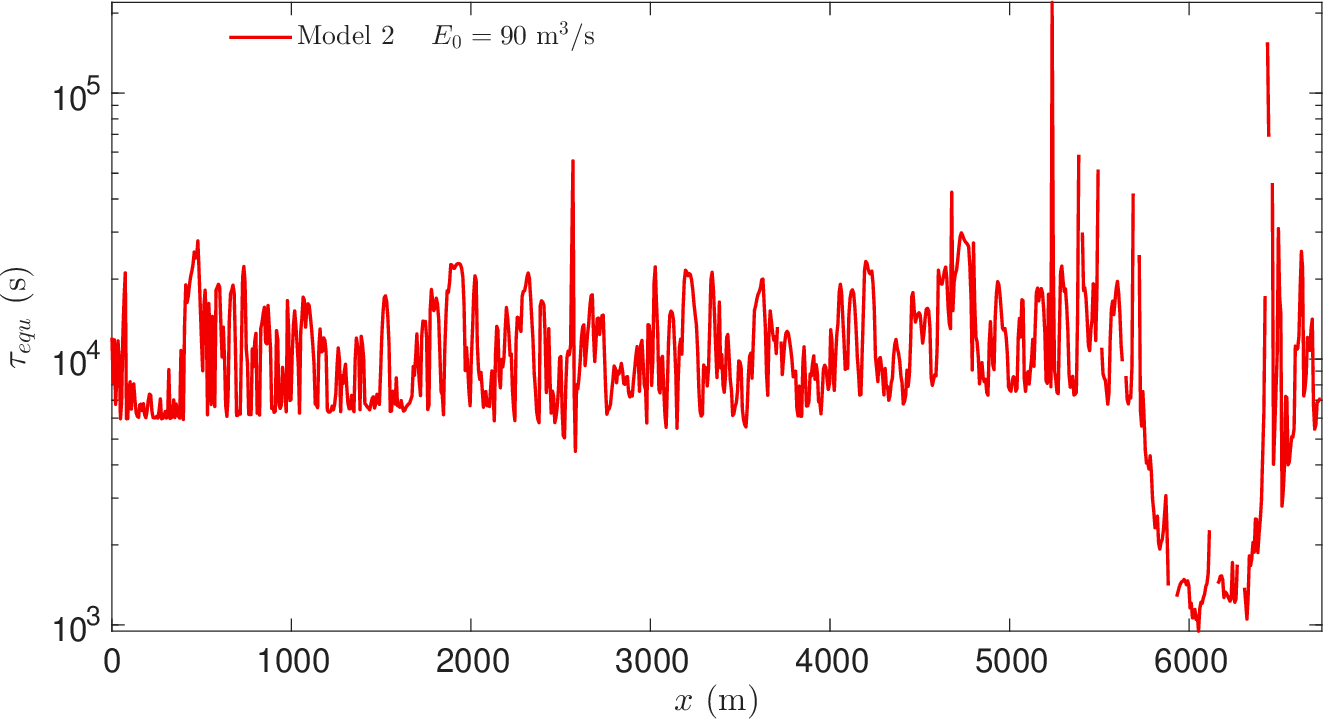}
	\caption{Equivalent crystallization rate of Model $2$ ($\tau_{equ}$) as a function of distance at $3000$ s.}
	\label{fig_harris_cryst}
\end{figure}

Figure \ref{fig_harris_cryst} shows the equivalent crystallization rate of Model $2$ as a function of distance. It suggests that a $\tau_{\rm equ}$ of $\approx 8\,000$ s is a reasonable first guess for Model $1$ to yield a crystal fraction evolution similar to that of Model $2$. We thus tried a few runs of Model $1$ with increasing values of $\tau_{\rm cryst}$ starting from $8\,000$ s. We show in Fig. \ref{fig_harris_u} a Model $1$ run with $\tau_{\rm cryst}=80\,000$ s ($\approx 1$ day, Table \ref{tab_param_appl}). It yields an evolution of crystal fraction slightly larger than Model $2$ and yet compatible with the field data. The cooling rate is faster than that of Model $2$, which results in a steeper increase of viscosity with distance (Fig. \ref{fig_harris_ntau}).\\

\cite{crisp1994influence} report and model a specific case in Hawaï of about 25 vol\% of microlites crystallized in $\approx 4$ hours and explore a range of $\tau_{\rm cryst}$ of 0.5--10 days. These timescales are consistent with our own range of exploration of 2 hours to 1 day. Model $1$ simplifies the complex evolution of cooling by parameterizing it with a simple crystallization rate modeled as an exponential decay. This simplification works because the dominant mechanism of deceleration is the viscosity increase due to crystallization (Fig. \ref{fig_harris_ntau}). The counterpart cost is the introduction of the additional free parameter $\tau_{\rm crys}$. Even using the outputs of Model $2$ as a guidance, we could not reduce the uncertainty on the range of $\tau_{\rm crys}$ below about one order of magnitude. \\

\section{Discussion on the model intercomparison}

Overall, our lava flow dynamics is strongly controlled by small-scale (hectometers) channel shape variation in width and slope, whereas PyFLOWGO predicts a large-scale (kilometers) deceleration. The transient nature of our model complicates a straightforward comparison because steep and rugged region cause pulsatory flow and thus only a quasi-steady state can be achieved. This strongly affects transit time calculations, slowing our transit times significantly compared to those of PyFLOWGO. Such effects are less important when focusing on the thermal evolution of both models, which can be quite similar when adjusting the initial discharge rate because this evolution is controlled by regional (kilometer) changes in slope and distance from the vent. \\

We compared Model $2$ to Model $1$, which replaces the multi-parametric evolution of temperature of Model $2$ by a much simpler crystallization relationship that adds one degree of freedom compared to Model $2$. Our results suggest that the tradeoff between the uncertainty associated with the additional degree of freedom of Model $1$ and the gain of eliminating temperature from Model $2$ (or PyFLOWGO) is uneven. In other words, replacing the apparently more complicated temperature equation \eqref{eq_phi_Lagr} by a simpler crystallization equation \eqref{eq_phi_Lagr2} introduces such large uncertainties on the input parameters (i.e. on $\tau_c$) that model results become unreliable. Model $2$ thus should be preferred over Model $1$.\\

One improvement brought by our model is to take gas into account. Our runs, however, used the viscosity formula that treats bubbles as hard spheres \eqref{rheo_prop}. The validity of this assumption that ${\rm Cx}\ll 1$ has to be tested, which can be done a posteriori by calculating a representative capillary number given by:
\begin{align*}
    \displaystyle {\rm Cx}=\sqrt{{\rm Ca}^2+{\rm Cd}^2}, \quad \text{with} \qquad 
    & {\rm Ca} = \frac{\eta_{solid} \, a}{\Gamma} \left|\frac{u}{h}\right|\\[1mm]
    & {\rm Cd} = \frac{\eta_{solid} \, a}{\Gamma} h,
\end{align*}
where $\eta_{solid} = \eta_0 \left(1-\Phi_p^*/\Phi_m^*\right)^{-2}$ is the viscosity of the mixture lava with crystals, $a=1$ mm is the bubbles radius and $\Gamma=0.24$ N/m \cite{COLUCCI2016113} is the surface tension at $1473$ K.\\
\begin{figure}[hbtp]
	\centering
	\includegraphics[width=0.9\textwidth]{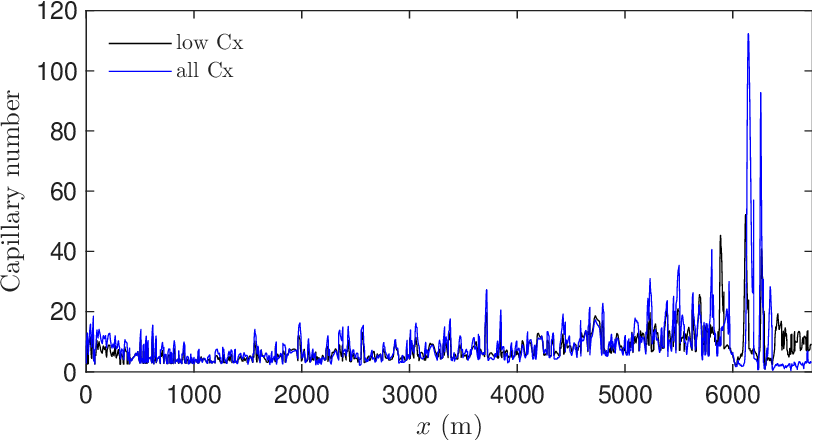}
	\caption{Capillary number evolution with distance for Model $2$ with $E_0=90$m$^3$/s at $3000$s for the hard sphere assumption (low ${\rm Cx}$) and the general viscosity (all ${\rm Cx}$).}
	\label{fig_harris_model2_4}
\end{figure}

The capillary number computed in the simulation at $E_0=90$ m$^3$/s is represented in Figure \ref{fig_harris_model2_4}. The number of peaks is higher when $x>5$ km because of the higher velocity and viscosity of the pulsatory flow. So, distally (or, more generally, whenever pulsatory behavior occurs), our formula \eqref{eq_rheology_final_app} for viscosity $\eta$ should be replaced by the general formula \eqref{rheo_prop_gral}. The temperature, crystal fraction, and speed in quasi-steady state of a run with the general viscosity relationship and $E_0=90$ m$^3$/s are shown in Fig. \ref{fig_harris_u}. The evolution of these three variables is hardly affected by the change in rheology until the final steep slope.\\

With the value of $\Phi_b$ used in our simulations, a change in the capillary number can decrease the value of the bulk viscosity by up to a factor two (for $\Phi_p^*=0.2$, the viscosity is divided by two when ${\rm Cx}$ varies from $0$ to $40$). This effect is quantified precisely in Fig. \ref{fig_harris_ntau}, which shows an almost proportional decrease in $\eta$ when the general viscosity is used instead of the low ${\rm Cx}$ approximation. As expected, the two other rheological parameters, $\tau_c$ and $n$, are mostly unaffected by this change in rheology. As a result, the a posteriori capillary numbers of both rheologies are very similar (Fig. \ref{fig_harris_model2_4}).

\section{Conclusion}

In this paper we developed a new three-phase lava flow model with a Herschel–Bulkley rheology. The rheology is derived combining semi-empirical formulas from the literature and comparing with data, and depends on the solid crystal fraction. We use a depth average approach of the two dimensional equations after performing a dimensional analysis, which yields a system of hyperbolic partial differential equations. Two versions of the model are described, Model $1$ involving a transport equation on the solid crystal fraction with a relaxation term, whereas Model $2$ involves a temperature dependence of this crystal fraction. The temperature equation is a transport equation with convective, conductive, radiative and viscous heating exchanges terms. \\

A numerical finite volume scheme was developed to solve this model. It is a semi-implicit HLL scheme with an hydrostatic reconstruction to capture the wet/dry transitions. We performed classical tests (Supplementary Information \ref{app:num}) to validate this scheme and verify its convergence. They show that the scheme is robust in different situations. It does not have instabilities, it converges to stationary solutions, and behaves well with wet/dry fronts. \\

Finally, we apply both closure variations of our model (Models $1$ and $2$) to the case of the Mauna Ulu lava flow (1969-1974). This channeled lava flow has already been studied with the PyFLOWGO kinematic model \cite{harris2022anatomy}, so we could carry out a model inter-comparison using field data constraints. Our model is the most complex of the two as it involves equations on the lava height $h$ and its velocity $u$ whereas in PyFLOWGO only a mean velocity is computed using an explicit formula.\\

The best-fit volumetric flux at the vent differs by a factor $\approx ~2$ between PyFLOWGO and our model; the smallest value being predicted by our model. Our lava flow dynamics is strongly controlled by small-scale (hectometers) channel shape variations in width and slope, whereas PyFLOWGO predicts a large-scale (kilometers) deceleration. The thermal evolution of both models is similar and is controlled by the distance from the vent and by large-scale changes in substratum slope. The transit time (i.e. the travel time of a lava aliquot from vent to some distal position), however, differs strongly between the two models. This is due to the transient nature of our model, which correctly predicts that confined lava traveling down a steep slope with irregular, small-scale breaks will yield a series of cascading fill-then-breakout lumps that cause the overall flow to be pulsatory. These pulses, when present, dominate the local dynamics and thus the local transit time. They also preclude a strict steady state to be reached.\\

We defined a model closure (Model $1$) with a simple crystallization relationship that had the potential to replace the multi-parametric evolution of temperature of the other variation (Model $2$). Broadly, these two closures help us to explore possible system simplification prior to extend our model in a future work by adding processes such as disequilibrium crystallization. Unfortunately, our results show that this simplification adds one very poorly constrained parameter (the characteristic crystallization time) that is a major control of flow dynamics. The thermal evolution of the flow done in Model $2$ thus yields much more accurate results than those obtained with the poorly constrained Model $1$.\\

One improvement brought by our model is the possible presence of up to 30 vol\% of gas. Theoretically, taking gas bubbles into account is best done with a general rheological relationship \eqref{rheo_prop_gral} that is valid for the whole range of capillary number to properly account for possible bubble deformation under shear. In the conditions explored herein, bubbles modulate viscosity within a factor 2 with a shear thinning behavior, decelerating slow flows and accelerating fast flows. We used a simplified rheology by treating bubbles as hard spheres (i.e. assuming low capillary number). Our results suggest that dynamic parameters other than bulk viscosity are hardly affected by this assumption.

\begin{appendices}
\section{Numerical scheme}\label{app:num_scheme}

In this section, we present a numerical scheme to solve the two different versions of the model \eqref{eq_syst1D_closed}-\eqref{eq_rheology_final}, with \eqref{eq_phi_p_final} for model 1, and \eqref{model2_final} for model 2.

\subsection{Uniform equations}

First we consider the uniform version of model \eqref{eq_syst1D_closed}, that is with constant quantities in the $x$ direction. 
A semi-implicit scheme is used for the fluid evolution; the gravity and yield stress terms are treated together while the viscosity term is treated in a semi-implicit way. Denoting by $f^k$ the approximation of a function $f$ at time $k \Delta t$, with $k\in \mathbb{N}$, the scheme is:
\begin{align*}
	\left\{
		\begin{array}{l}
			\displaystyle h^{k+1} = h^k=h^0,\\
			\displaystyle u^{k+1} = u^{k+1/2} - \Delta t\, \frac{a_f \, \eta^k \, s_b^{(n^k)}}{\rho} \frac{|u^{k+1/2}|^{(n^k)-1}}{(h^k)^{(n^k)+1}} u^{k+1},
		\end{array}
	\right.
\end{align*}
where 
\begin{align*}
	& \displaystyle u^{k+1/2} = u^k + \Delta t g \sin \theta + \Delta t \, f^k,\\
	& \displaystyle \text{with } f^k = - \underset{\tau_c^k / (\rho h^k)}{\text{proj}} \left( \frac{u^k}{\Delta t} + g \sin \theta \right),
\end{align*}
and $\tau_c^k=\tau_c(\Phi_p^k)$, $n^k=n(\Phi_p^k)$,  $\eta^k=\eta(\Phi_p^k)$ defined by (\ref{eq_rheology_final_eta})-(\ref{eq_rheology_final_n}).
The function proj is given by:
\begin{align}
	\underset{a}{\text{proj}} \, b = 
	\left\{ \begin{array}{ll}
		b &\text{if } |b| < a,\\
		a \, \text{sgn}(b) & \text{if } |b| \geq a,
	\end{array} \right.
	\label{eq_function_proj}
\end{align}
so that the definition of $f$ is such that the gravity and yield stress terms do not change the sign of $u$.
The uniform version of the linear equation on $\Phi_p^*$ \eqref{eq_phi_p_final} is discretized using an implicit method:
\begin{align*}
	\displaystyle \Phi_p^{*k+1} = \frac{\Phi_p^{*k} + \Delta t \, \Phi_{crys}^*/\tau_{crys}}{1+\Delta t / \tau_{crys}}.
\end{align*}
The uniform equation on $T$ \eqref{eq_rheology_phi_version2} is treated in a semi-implicit fashion by using a semi-linearization around $T^{k}$ for the non-linear term $T^4_{env}-T^4$:
\begin{align*}
	T_{env}^4-\left(T^{k+1}\right)^4 
	& = T_{env}^4-\left(T^k\right)^4 + (T^k)^4-(T^{k+1})^4\\
	& = T_{env}^4-\left(T^k\right)^4 + 4 \left(T^k\right)^3 \left(T^k-T^{k+1}\right) + O\left(\left(T^k-T^{k+1}\right)^2\right)\\
	& = T_{env}^4 + 3 \left(T^k\right)^4 - 4 \left(T^k\right)^3 T^{k+1} + O\left(\left(T^k-T^{k+1}\right)^2\right).
\end{align*}
Consequently the scheme is written as
\begin{align*}
	\displaystyle T^{k+1} = \frac{h^0 T^k + \Delta t \left(C_{rad} \left(T_{env}^4 + 3\left(T^k\right)^4 \right) + C_{conv} T_{env} + C_{cond} T_c / h_0 + C_{visc} \eta^{k} \left(u^{k+1}\right)^2 / h_0 \right)}{h^0 + \Delta t \left(4 C_{rad} \left(T^k\right)^3 + C_{conv} + C_{cond}/h^0\right)}.
\end{align*}

\subsection{Non-uniform equations} \label{sec_model_num}

An HLL finite volume scheme is used for the fluid evolution. As for the uniform scheme, the gravity and yield stress are treated with an hydrostatic reconstruction and the viscosity term is treated semi-implicitly. We introduce:
\begin{align*}
	\displaystyle 
	W = \left( \begin{array}{c} h \\ q \\[1mm] w_3 \end{array} \right), \quad
	F(W) = \left( \begin{array}{c} q \\ \displaystyle  \frac{q^2}{h} \\[3mm] \displaystyle \frac{q w_3}{h}\end{array} \right), \quad
	S(W) = \left( \begin{array}{c} 0 \\ \displaystyle  g h \cos \theta \\ 0\end{array} \right),
\end{align*}
where $q=hu$, $w_3=h\Phi_p^*$ for Model $1$, and $w_3=hT$ for Model $2$.
Then the system \eqref{eq_syst1D_closed} is written as
\begin{align*}
	\displaystyle 
	\partial_t W + \partial_x F(W) + S(W) \partial_x (h+z_b) =
	 \left( \begin{array}{c} 0 \\[1mm]
	 \displaystyle  - \frac{a_f}{\rho} \left(\eta \left|s_b\frac{u}{h} \right|^{n} + \tau_c\right) \frac{u}{|u|}\\[3mm]
	 \displaystyle S_{w_3}
	 \end{array} \right),
\end{align*}
where the source term $S_{w_3}$ depends on the model:
\begin{equation*} \label{eq:def_zb}
    S_{w_3} = 
    \left\{ \begin{array}{ll}
        \displaystyle - \frac{h}{\tau_{\rm crys}} \left(\Phi_p^*-\Phi_{\rm crys}^* \right) & \text{for Model $1$,} \\
        \displaystyle C_{rad} (T_{env}^4-T^4) + C_{conv} (T_{env}-T) + \frac{C_{cond}}{h} (T_{c}-T) + \frac{C_{visc} \, \eta}{h} u^2  & \text{for Model $2$.}
    \end{array} \right.
\end{equation*}
The numerical scheme is given by:
\begin{align*}
	\left\{
	\begin{array}{l}
		\displaystyle h_i^{k+1} = h_i^{k+1/2},\\
		\displaystyle u_i^{k+1} = u_i^{k+1/2} - \Delta t a_f \frac{\eta_i^{k}s_b^{(n_i^k)}}{\rho} \frac{|u_i^{k}|^{(n_i^{k})-1}}{(h_i^{k})^{(n_i^{k})+1}} u_i^{k+1},\\
		w_{3,i}^{k+1} = f_{w_3}(W_{i}^{k},W_{i}^{k+1}),
	\end{array}
	\right.
\end{align*}
where $i$ and $k$ denote respectively the space and time discretization, and $f_{w_3}$ is described at the end of this section. The intermediate step $k+1/2$ consists in the HLL scheme with hydrostatic reconstruction:
\begin{align*}
	& \displaystyle W_i^{k+\frac{1}{2}} = W_i^k + 
	\frac{\Delta t}{\Delta x} \left(F_{i+\frac{1}{2}}^k-F_{i-\frac{1}{2}}^k + 
	\frac{1}{2}S_{i+\frac{1}{2}}^k \left((h_{i+\frac{1}{2}}^k)^+ - (h_{i+\frac{1}{2}}^k)^-\right) 
	+ \frac{1}{2}S_{i-\frac{1}{2}}^k \left((h_{i-\frac{1}{2}}^k)^+ - (h_{i-\frac{1}{2}}^k)^-\right)\right),
\end{align*}	
where
\begin{align*}
	\displaystyle 
	& F_{i+\frac{1}{2}}^k = \frac{F\left(\left(W_{i+\frac{1}{2}}^k\right)^+\right) + F\left(\left(W_{i+\frac{1}{2}}^k\right)^-\right)}{2} - \frac{1}{2} V_{i+\frac{1}{2}}^k,
	& S_{i+\frac{1}{2}}^k = \frac{S\left(\left(W_{i+\frac{1}{2}}^k\right)^-\right) + S\left(\left(W_{i+\frac{1}{2}}^k\right)^-\right)}{2},
\end{align*}
where $V_{i+\frac{1}{2}}^k$ is the numerical viscosity given by:
\begin{align*}
	 \displaystyle & V_{i+\frac{1}{2}}^k = \alpha_{i+\frac{1}{2}}^k \left(\left(W_{i+\frac{1}{2}}^k\right)^+ - \left(W_{i+\frac{1}{2}}^k\right)^-\right) + \beta_{i+\frac{1}{2}}^k \left(F\left(W_{i+\frac{1}{2}}^k\right)^+ - F\left(W_{i+\frac{1}{2}}^k\right)^-\right),\\[3mm]
	 \displaystyle & \alpha = \frac{S_R|S_L| - S_L|S_R|}{S_R-S_L}, \quad
	 \beta_{i+\frac{1}{2}}^k = \frac{|S_R| - |S_L|}{S_R-S_L},\\[1mm]
	 \displaystyle & S_R = \max \left(\lambda^+,\lambda^- \right), \quad  S_L = \min \left(\lambda^+,\lambda^- \right).
\end{align*}
The quantities $\lambda^{\pm} = q/h \pm \sqrt{gh \cos \theta}$ are the eigenvalue of $A$, which is the Jacobian matrix of $F$. 
The reconstructed heights $h^+$ and $h^-$ include the bottom surface height and the yield stress term, and are given by:
\begin{align*}
	\left\{
	\begin{array}{l}
		\displaystyle (h_{i+\frac{1}{2}}^k)^+ = \left(h_{i+1}^k - \left(-z_{b,i+1}+z_{b,i} + \frac{\Delta x f_{i+\frac{1}{2}}^k}{g \cos \theta }\right)_+\right)_+,\\
		\displaystyle (h_{i+\frac{1}{2}}^k)^- = \left(h_{i}^k - \left(z_{b,i+1}-z_{b,i} - \frac{\Delta x f_{i+\frac{1}{2}}^k}{g \cos \theta }\right)_+\right)_+,
	\end{array}
	\right.
\end{align*}
where 
$$
z_{b,i}=b_i+\tan \theta \, x_i	 \quad \displaystyle \text{and} \quad  f^k_{i+\frac{1}{2}} = - \underset{\tau_{c,i+\frac{1}{2}}^k / (\rho h_{i+\frac{1}{2}}^k )}{\text{proj}} \left(\frac{u_{i+\frac{1}{2}}^k}{\Delta t} - \frac{g \cos \theta}{\Delta x} \left(z_{b,i+1}-z_{b,i} + h_{i+1}^k-h_i^k) \right) \right).
$$
The function proj is defined by equation \eqref{eq_function_proj}. The reconstructed discharge $q$ and the crystals--temperature quantity $w_3$ are given by:
\begin{align*}
	\left\{
	\begin{array}{l}
		\displaystyle (q_{i+\frac{1}{2}}^k)^+ = \frac{q_{i+1}^k}{h_{i+1}^k} (h_{i+\frac{1}{2}}^k)^+\\[4mm]
		\displaystyle (q_{i+\frac{1}{2}}^k)^- = \frac{q_{i}^k}{h_{i}^k} (h_{i+\frac{1}{2}}^k)^-
	\end{array}
	\right., \qquad
	\left\{
	\begin{array}{l}
		\displaystyle (w_{3,i+\frac{1}{2}}^k)^+ = \frac{w_{3,i+1}^k}{h_{i+1}^k} (h_{i+\frac{1}{2}}^k)^+\\[4mm]
		\displaystyle (w_{3,i+\frac{1}{2}}^k)^- = \frac{w_{3,i}^k}{h_{i}^k} (h_{i+\frac{1}{2}}^k)^-.
	\end{array}
	\right.
\end{align*}
Finally, a CFL condition is used to compute the time step:
\begin{align*}
	\Delta t^k = \frac{C \Delta x}{\underset{i}{\max} \, \lambda_i^k},
\end{align*}
where $\lambda_i^k$ are the eigenvalues of the Jacobian matrix $A_i^k$. \\ \\
The discretization of the source term for $w_3$ is made with an Euler-implicit discretization in time for Model $1$ and a semi-implicit one for Model $2$. That is, for $w_3 = h\Phi_p$ we have:
\begin{align*}
	\displaystyle \left(h\Phi_p^*\right)_i^{k+1} = \frac{\left(h\Phi_p^*\right)_i^{k+1/2} + \Delta t \, h_i^{k} \Phi_{crys}^*/\tau_{crys}}{1+\Delta t / \tau_{crys}}.
\end{align*}
For $w_3 = h  T$ we use the same semi-implicit discretization of the source term as in the uniform scheme:
\begin{align*}
	\displaystyle (hT)_i^{k+1}
	& =  \Bigg( (hT)_i^{k+1/2} + \Delta t \bigg( C_{rad} \left(T_{env}^4 + 3\left(T_i^{k}\right)^4 \right) 
	 + C_{conv} T_{env} + C_{cond} T_c / h_i^{k} \\ & + C_{visc} \eta_i^{k} \left(u_i^{k}\right)^2 / h_i^{k} \bigg) \Bigg) / \Bigg( 1+ \Delta t \left( 4 C_{rad} \left(T_i^{k}\right)^3/h_i^{k} + C_{conv}/h_i^{k} + C_{cond}/\left(h_i^{k}\right)^2 \right) \Bigg).
\end{align*}

\bigskip

Let us finally describe the treatment of wet/dry fronts, where some regularization is necessary to approximate the velocity and temperature, calculated by the terms $q/h$ and $(hT)/h$. The first term is usually approximated in Saint-Venant models as $u=q/h$ in wet/dry fronts, where we assume that if $h=0$ then $u=0$. From the numerical point of view a regularization technique must be considered to treat the division by zero, we set the following one,
$$
u \approx  \frac{ q \,  h \, \sqrt{2}}{\sqrt{h^4+\max(h\, , \, 10^{-5})^4}}.
$$
In Model 2 we have an extra difficulty, related to the computation of $T$ in wet/dry fronts that is calculated as $(hT)/h$. The source term of the equation associated with $(hT)$ is defined in terms of $T$, see equation \eqref{model2_final}, so it must be solved. If we use a similar regularization as the previous one, then we get $T=0$ in wet/dry fronts, implying that the lava is solid. As a consequence, the lava flow stops, leading to a non realistic situation.
For this reason, we use a new regularization, which set $T=T_l$ in wet/dry fronts, defined as follows,
$$
T \approx \frac{1}{2} \bigg( (\, T_l-\widehat{T} \, ) \tanh(10- 10^3\, h)+T_l+\widehat{T} \bigg), 
$$
where $\widehat{T}=(hT) \,h\,  \sqrt{2} /\sqrt{h^4+\max(h,10^{-5})^4}$. Note that in Model 1 we do not find this problem because the source term of the evolution equation of $(h \phi_p^*)$ is multiplied by $h$, see equation \eqref{eq_phi_p_final}, and then we do not need to calculate $\phi_p^*$.

\end{appendices}

\section*{Acknowledgments}
This research has been partially supported by the European Union - NextGenerationEU program and by grant PID2022-137637NB-C22 funded by  MCIN/AEI/10.13039/501100011033 and ``ERDF A way of making Europe''. Authors thank IMUS-Maria de Maeztu grant CEX2024-001517-M - Apoyo a Unidades de Excelencia María de Maeztu for supporting this research, funded by MICIU/AEI/10.13039/501100011033 and the applied research and innovation project "Development of mathematical tools for transfer" (Mathware) at IMUS, co-financed by the EU - Ministry of Finance and Public Administration - European Funds - Junta de Andalucía -Consejería de Universidad, Investigación e Innovación and the projects “Programa operatIvo PPIT-FEDER Andalucía 2021-2027 SOL2024-31596 and SOL2024-31708 for partially supporting this research.

\bibliographystyle{apalike}
\bibliography{biblio} 

\begin{thebibliography}{99}
\bibitem[Chevrel et al., 2018]{chevrel:2018} Chevrel, M. O., Labroquère, J., Harris, A. J., and Rowland, S. K.
(2018). Pyflowgo: An open-source platform for simulation of channelized lava thermo-
rheological properties. Computers \& Geosciences, 111:167–180.
\bibitem[Dobran, 2001]{dobran:2001} Dobran, F. (2001). Volcanic Processes. Springer New York, NY.
\bibitem[Wilson and Parfitt, 1993]{wilson:1993} Wilson, L. and Parfitt, E. A. (1993). The formation of perched
lava ponds on basaltic volcanoes: the influence of flow geometry on cooling-limited lava
flow lengths. Journal of Volcanology and Geothermal Research, 56(1):113–123.
\end{thebibliography}

\unappendix

\clearpage
\setcounter{page}{1}

\beginsupplement

\addcontentsline{toc}{section}{Supplementary Information}

\begin{center}
{\huge Supplementary Information to:\\ A transient depth-averaged lava flow model with a Herschel–Bulkley rheology accounting for three phases} 

\vspace{1cm}

{\large J. Binard$^1$, A. Burgisser$^2$, E.D. Fernández-Nieto$^1$, G. Narbona-Reina$^1$}
\end{center}

{\scriptsize $^1$ Dpto. Matemática Aplicada I E.T.S Arquitectura, Universidad de Sevilla, Sevilla, Spain}

{\scriptsize $^2$ Univ. Grenoble Alpes, Univ. Savoie Mont Blanc, CNRS, IRD, Univ. Gustave Eiffel, ISTerre, Grenoble, France}

\vspace{1.1cm}

\begin{abstract}
    This Supplementary Information document contains a section listing the various numerical tests (Tests 1--5) performed with our lava flow model and a section describing how to recover viscous shear effects in 1D models.
\end{abstract}

\section{Numerical tests}\label{app:num}

We perform some test of the numerical schemes to ensure that they behave well in different situations.
The rheology parameters $\alpha_n = 1.32$ and $\Phi_p^c = 0.14$ are chosen according to the data comparison done in figure \ref{fig_compar_Truby2}.

\begin{table}[hbtp]
	\center 
	\begin{tabular}{| l| l | l |}
		\hline
		$L$ & $3.5$ m & Domain length \\
		$N$ & $400$ & Number of discretisation points \\
		CFL & $0.8$ & Constant for CFL condition\\
		\hline  
	\end{tabular}
	\caption{Numerical parameters of the simulations}
	\label{tab_paramnum_simu_nonuniform}
\end{table}

\begin{table}[hbtp]
	\center 
	\begin{tabular}{| l| l | l |}
		\hline
		$\alpha_n$ & $1.32$ & coefficient of the flow index equation\\
		$\Phi_b$ & $0.2$ & volume fraction of bubbles \\
		$\Phi_m^*$ & $0.6$ & maximum volume fraction of crystals \\
		$\Phi_p^c$ & $0.14$ & critical value for yield stress \\		
		$\theta$ & $15^\circ$ & angle of the tilted plane \\
		$\rho$ & $2500$ kg.m$^{-3}$ & fluid density \\
		$\eta_0$ & $26$ Pa.s$^n$ & consistency factor\\
		$\tau_c^*$ & $33$ Pa & characteristic critical yield stress\\
		$\tau_c^{\max}$ & $1.2 \times 10^{11}$ Pa & maximum yield stress\\
        $a_f$ & 1 & friction coefficient\\
        $s_b$ & $1$ & coefficient for shear rate approximation\\
		$T_l$ & $1474$ K & liquidus temperature \\
		$T_s$ & $1268$ K & solidus temperature \\
		$T_{env}$ & $293$ K & external (air) temperature \\
		$\tau_{crys}$ & $10$ s & characteristic crystallisation time\\
		$C_{rad}$ & $10^{-12}$ m.s$^{-1}$.K$^{-3}$ & constant for radiative loss \\
		$C_{conv}$ & $2 \times 10^{-4}$ m.s$^{-1}$.K$^{-3}$ & constant for convective loss \\
		$C_{cond}$ & $ 10^{-5}$ m.s$^{-1}$.K$^{-3}$ & constant for conductive loss \\
		$C_{visc}$ & $10^{-4}$ m.s$^{-1}$.K$^{-3}$ & constant for viscous loss \\
		\hline  
	\end{tabular}
	\caption{Model parameters for tests $1-5$ of the numerical results. }
	\label{tab_param_test_scheme}
\end{table}

\subsection{Test 1: comparison of uniform solutions of Models 1 and 2}

\begin{figure}[hbtp]
	\centering
		\begin{subfigure}[t]{0.46\textwidth}
		\centering     
		\includegraphics[width=\textwidth]{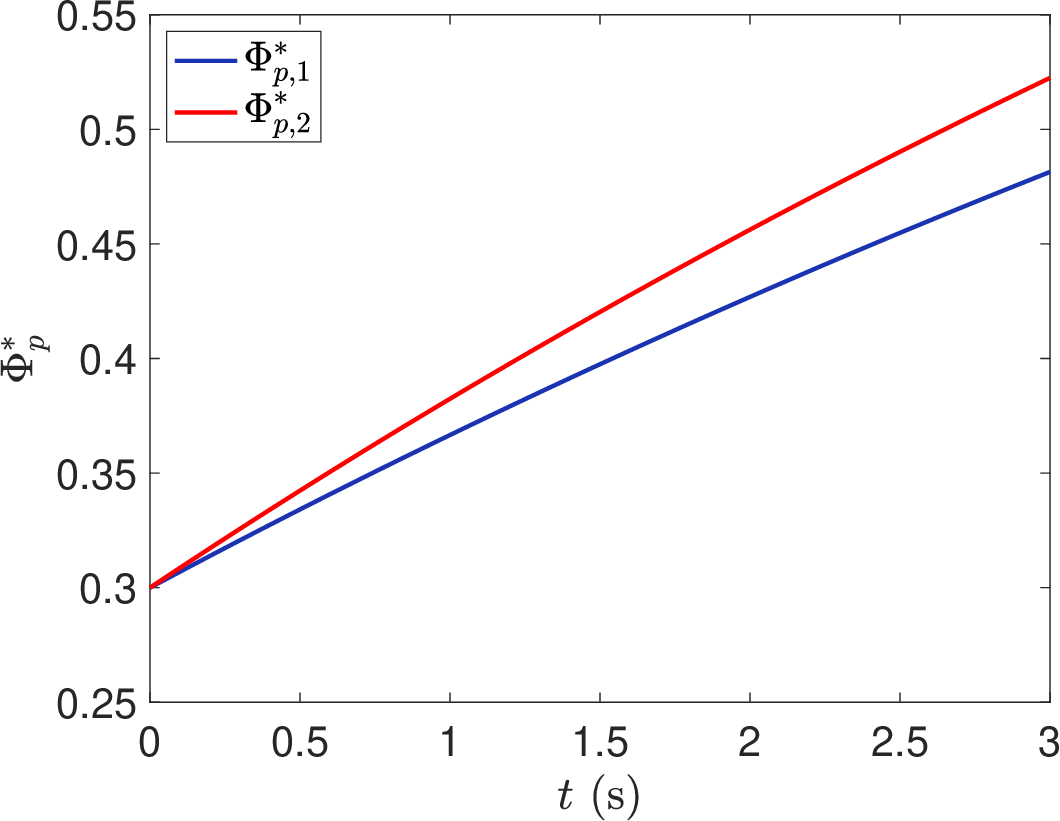}
		\caption{Evolution in time of $\Phi_p^*$.}
		\label{fig_uniform_phi}
	\end{subfigure}
\quad 
    \begin{subfigure}[t]{0.5\textwidth}
		\centering      
		\includegraphics[width=\textwidth]{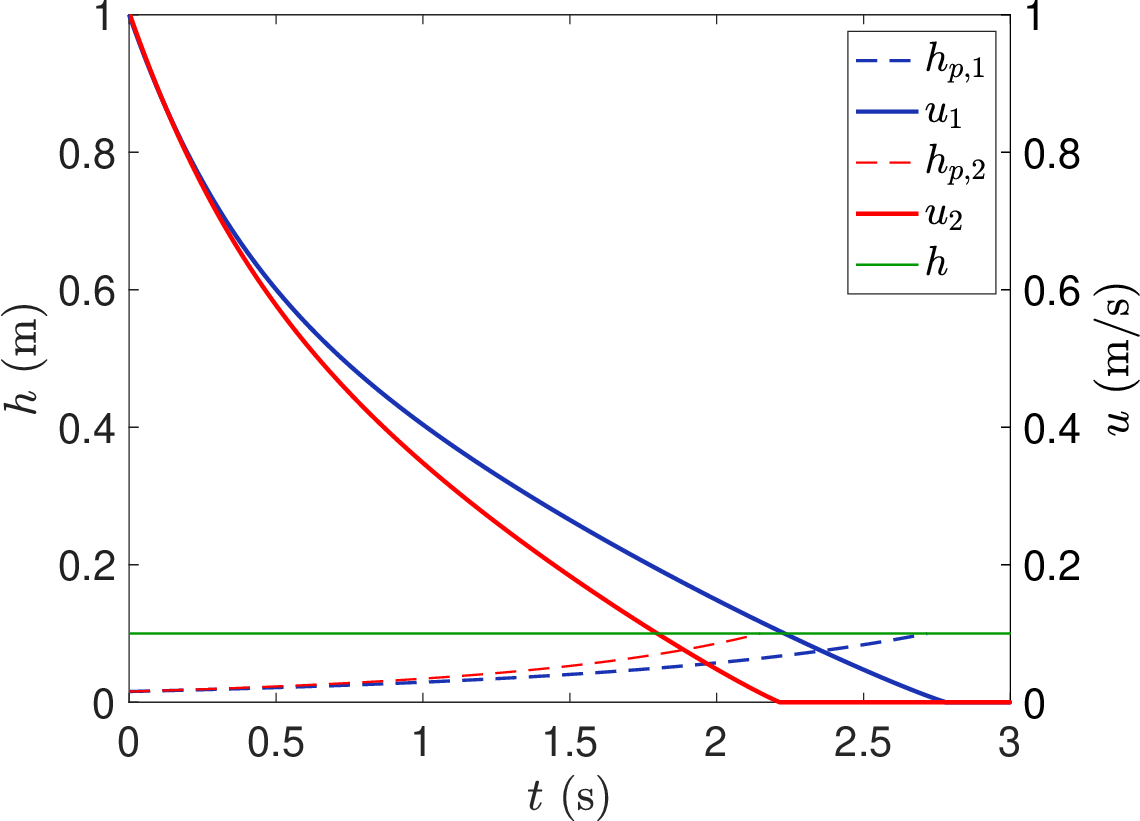}
		\caption{Evolution in time of $h_p$, $u$.}
		\label{fig_uniform_h_u}
	\end{subfigure}\\[2mm]
	\begin{subfigure}[t]{0.5\textwidth}
		\centering     
		\includegraphics[width=\textwidth]{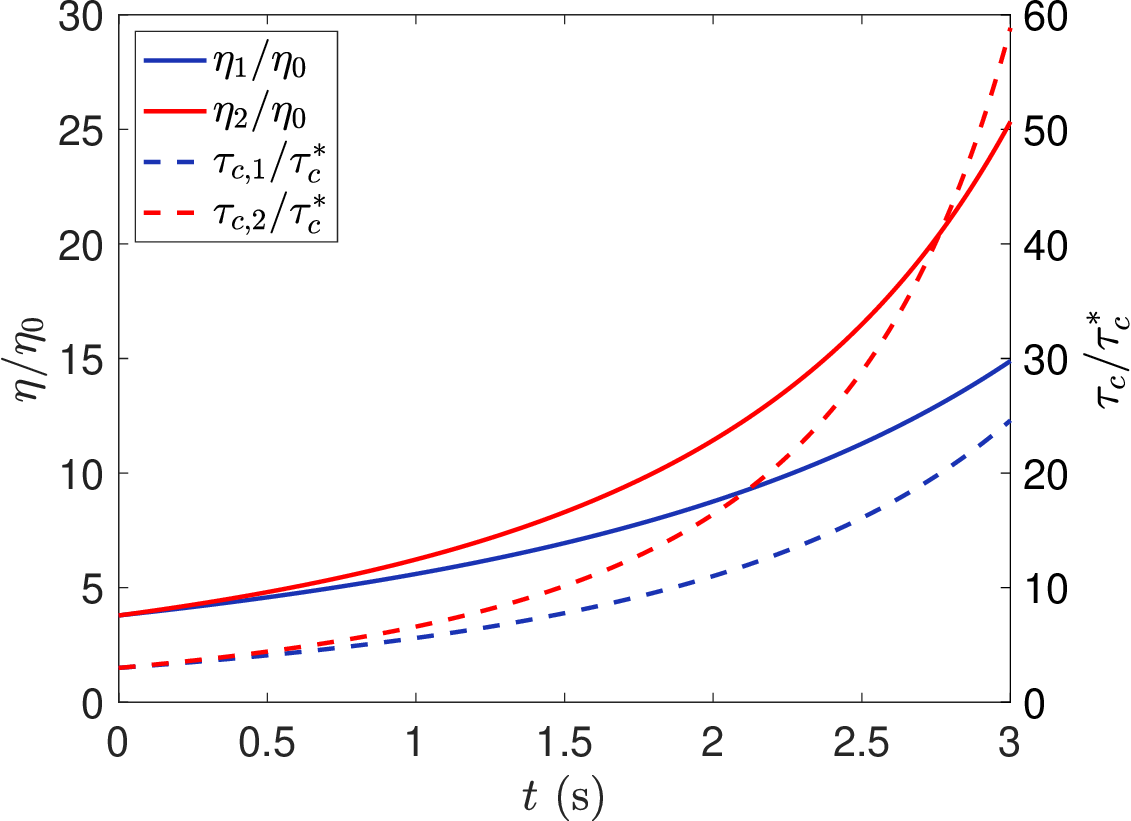}
		\caption{Evolution in time of $\eta$, $\tau_c$.}
		\label{fig_uniform_eta}
	\end{subfigure}
	\quad
	\begin{subfigure}[t]{0.46\textwidth}
		\centering     
		\includegraphics[width=\textwidth]{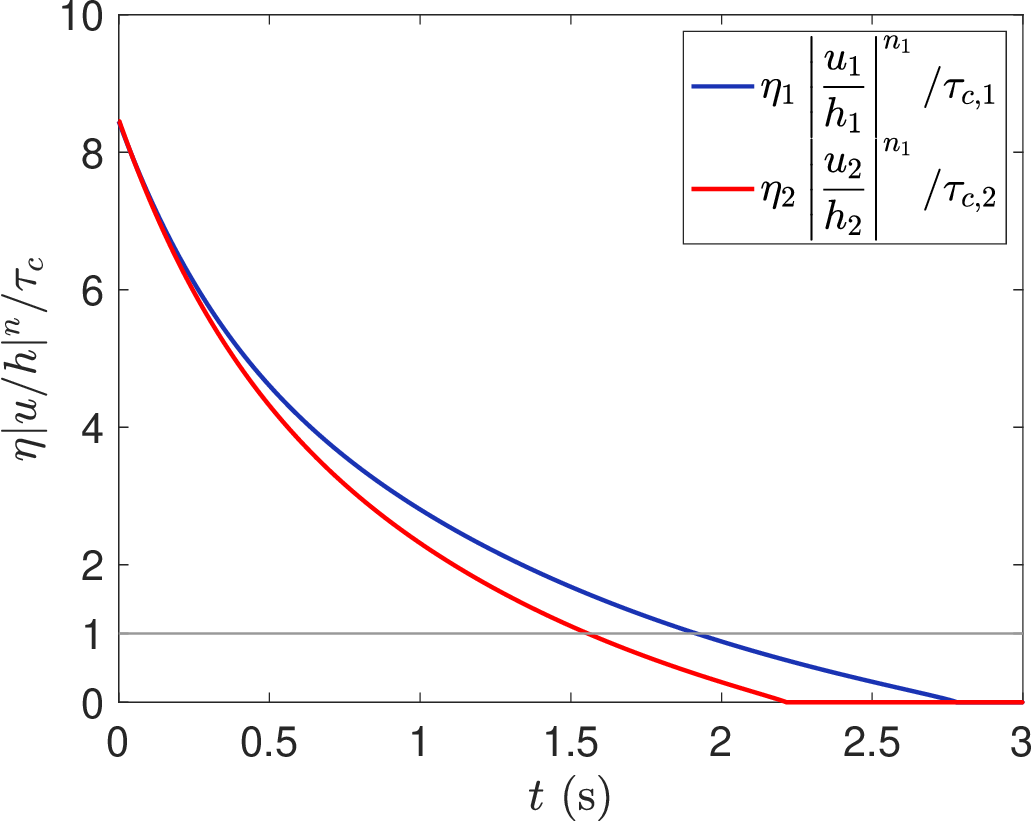}
		\caption{Evolution in time of $\eta \left| \frac{u}{h} \right|^{n}/\tau_{c}$.}
		\label{fig_uniform_etatauc}
	\end{subfigure}
	\caption{Test $1$, with $h=0.1$ m, $u|_{t=0} = 1\,$m/s, $\Phi_p|_{t=0} = 0.3$, $T|_{t=0} = 1412$ K. Model $1$ in blue, Model $2$ in red.}
	\label{fig_uniform}
\end{figure}

The parameters used for the two models are listed in Table \ref{tab_param_test_scheme}. The flow height $h$ is constant, equal to $0.1$ m and the initial discharge $u_0=1$ m/s. The initial temperature is set to $1412$ K and the initial crystal volume fraction $\Phi_p^*(0)=0.5 \Phi_m^*$. Figure \ref{fig_uniform} shows the temporal evolution of the quantities $h$, $q$, $\Phi_p^*$, $\eta$, $\tau_c$.
First, $\Phi_p^*$ as a function of time for Models $1$ and $2$ is showed in figure \ref{fig_uniform_phi}. The evolution of $u$ and $h_p = \frac{\tau_c}{\rho g \sin \theta}$ are showed in figure \ref{fig_uniform_h_u}. In accordance with the equations, the discharge $q$ vanishes when $h_p \geq h = 0.1$ m, so when $t\geq 2.2$ s for Model $1$ and $t\geq 2.7$ s for Model $2$. For this test the ratios $\eta/\eta_0$ and $\tau_c/\tau_c^*$ are very similar in both models, albeit they are mostly smaller with Model 1 (Figure \ref{fig_uniform_eta}). In Figure \ref{fig_uniform_etatauc} the ratio between viscous and drag terms are presented. In both models the viscous effects are greater than the friction effects initially and this ratio decreases with time because of the deceleration.
\clearpage

\subsection{Test 2: Convergence test to stationary solutions}\label{sec:test2} 

\begin{figure}[hbtp]
	\centering
	\includegraphics[width=0.5\textwidth]{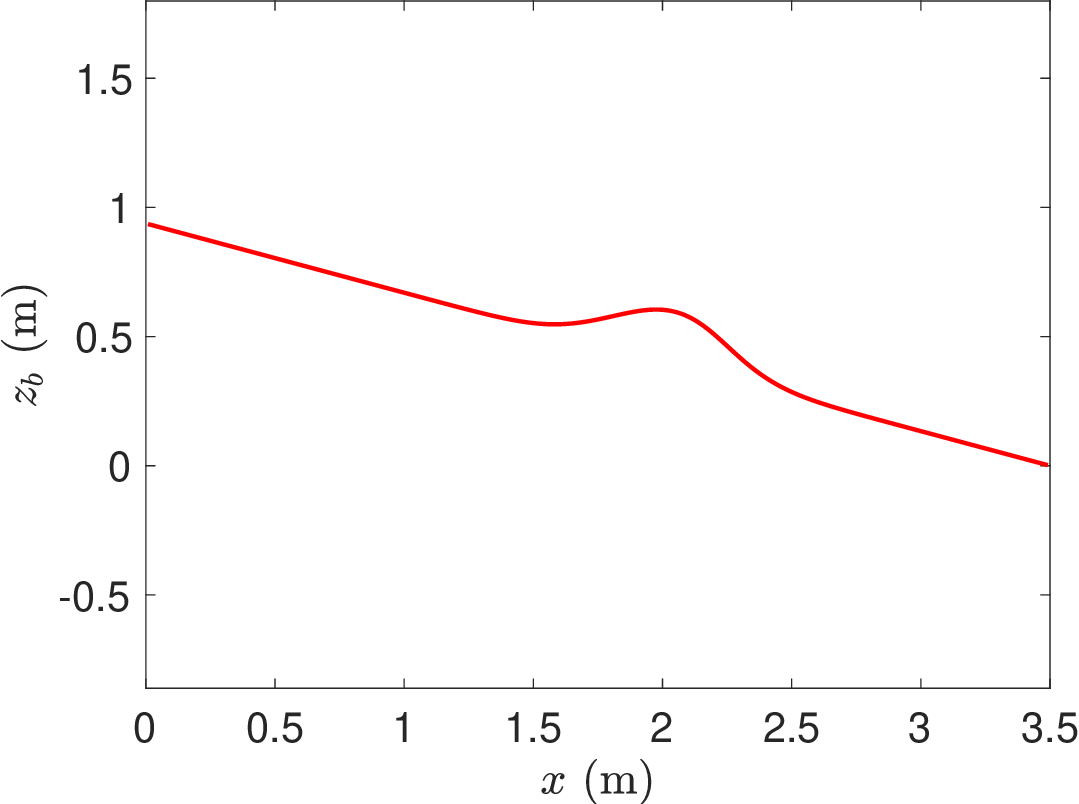}
	\caption{Bottom surface with a bump.}
	\label{fig_statio_bottom_z}
\end{figure}

In this test, we compare the solutions obtained after sufficient time with the corresponding stationary solutions. Parameters are given in Tables \ref{tab_paramnum_simu_nonuniform} and \ref{tab_param_test_scheme}. The stationary solution $h$ is a solution of the ODE \eqref{eq_sol_statio_h}, computed using an explicit Euler method. A Dirichlet boundary condition is used, at the upstream boundary in the supercritical case and the downstream one for the subcritical case. In Model $1$ the fraction $\Phi_p^*$ is computed with the explicit equation \eqref{eq_sol_statio_phi} for a supercritical solution, and with the following equation for a subcritical one:
\begin{align}
	&\displaystyle  \Phi_p^*(x) = 
		\displaystyle \left(\Phi_p^*(L)-\Phi_{crys}^*\right) \exp \left(\frac{\int_{L}^x h}{q_0 \tau_{crys}}\right) + \Phi_{crys}^* & \text{ if } q_0 \neq 0.
	\label{eq_sol_statio_phi_sub}
\end{align}
In Model $2$ the lava temperature $T$ is solution of equation \eqref{eq_sol_statio_T} and computed numerically. \\
To compute the solution with the general scheme we use the same boundary conditions for $h$, $q$, $\Phi_p^*$ as in test $2$. The bottom surface is a tilted plane with a bump, given by 
\begin{align}
	b(x) = 0.2 \, e^{-10(x-2)^2},
	\label{eq_bump}
\end{align}
and represented in figure \ref{fig_statio_bottom_z}. \\ \\

\begin{figure}[hbtp]
	\centering
	\begin{subfigure}[t]{0.49\textwidth}
		\centering      
		\textbf{Model $1$}\par\medskip
		\includegraphics[width=\textwidth]{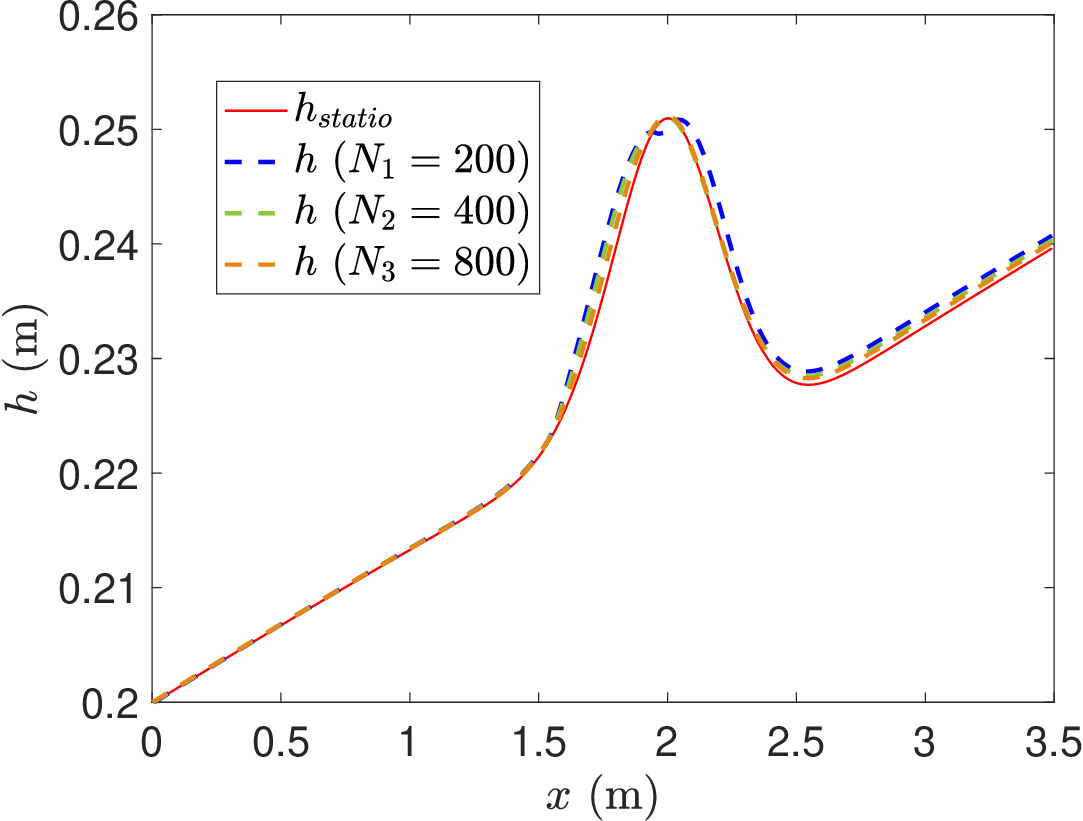}
		\label{fig_statio_friction_smaller_dx_h}
	\end{subfigure}
	\begin{subfigure}[t]{0.49\textwidth}
		\centering
	    \textbf{Model $2$}\par\medskip
		\includegraphics[width=\textwidth]{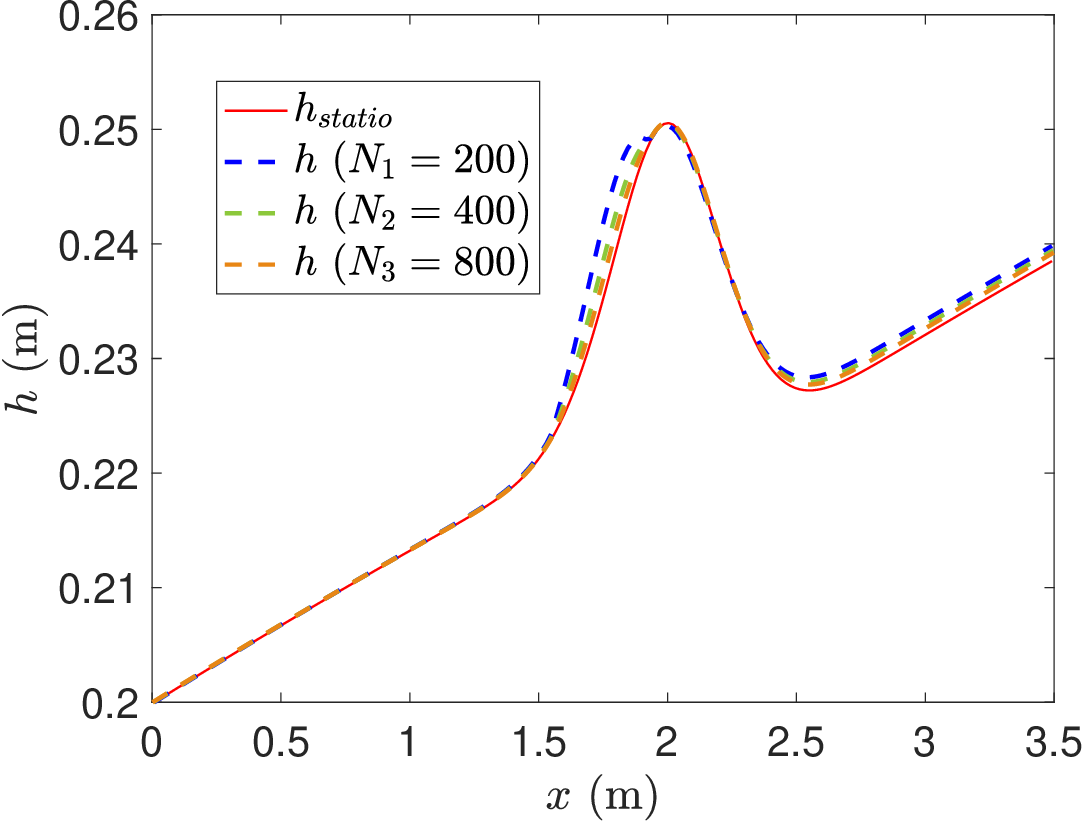}
		\label{fig_statio_friction_smaller_dx_model2_h}
	\end{subfigure}\\
	
	\begin{subfigure}[t]{0.49\textwidth}
		\centering      
		\includegraphics[width=\textwidth]{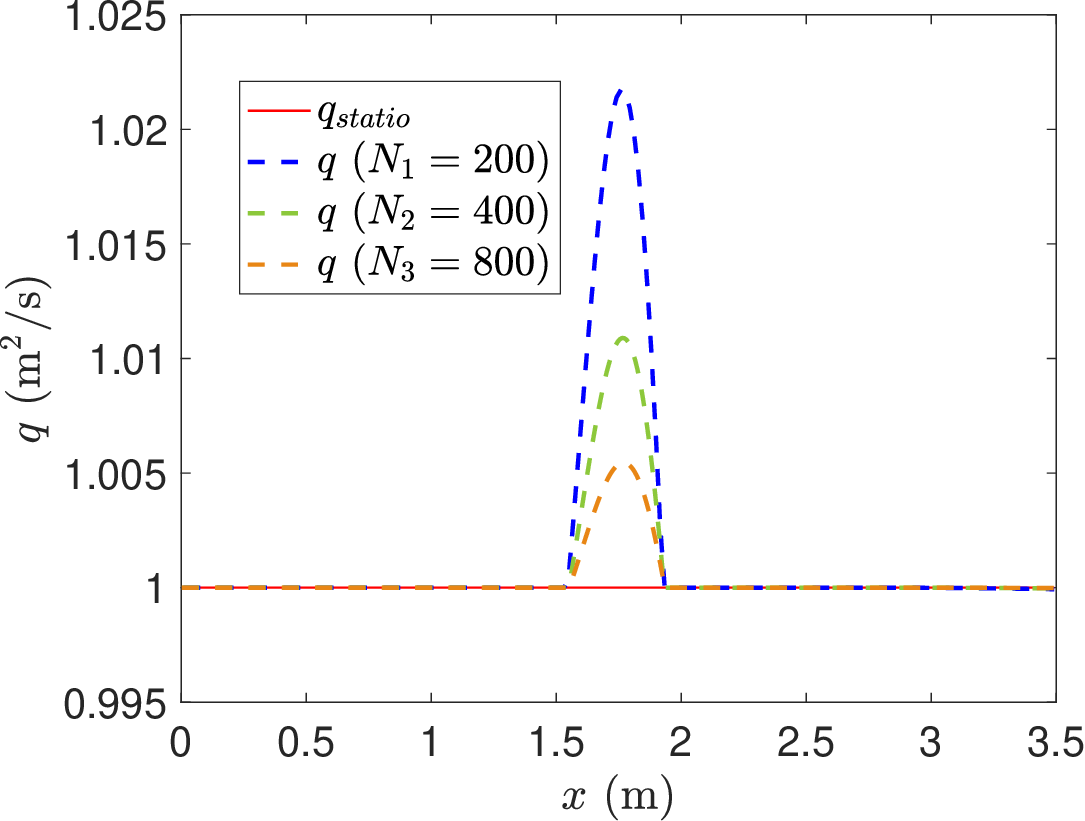}
		\label{fig_statio_friction_smaller_dx_q}
	\end{subfigure}
	\begin{subfigure}[t]{0.49\textwidth}
		\centering      
		\includegraphics[width=\textwidth]{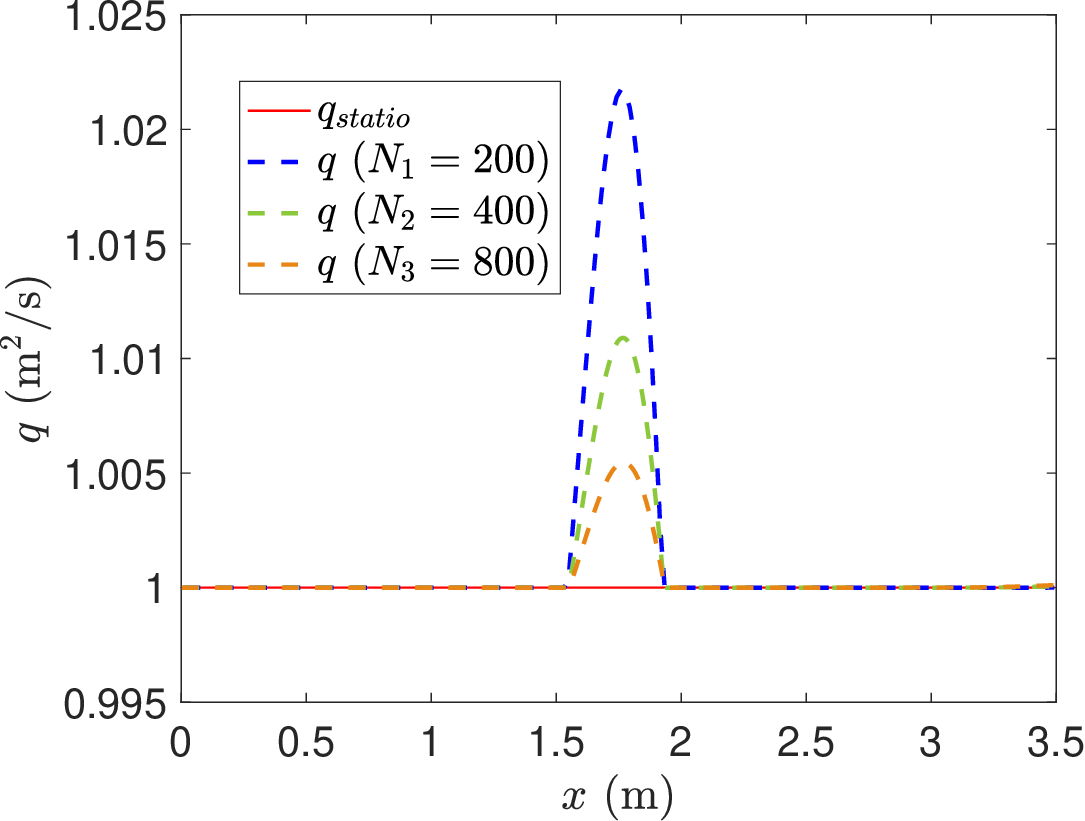}
		\label{fig_statio_friction_smaller_dx_model2_q}
	\end{subfigure}\\
	
	\begin{subfigure}[t]{0.49\textwidth}
		\centering      
		\includegraphics[width=\textwidth]{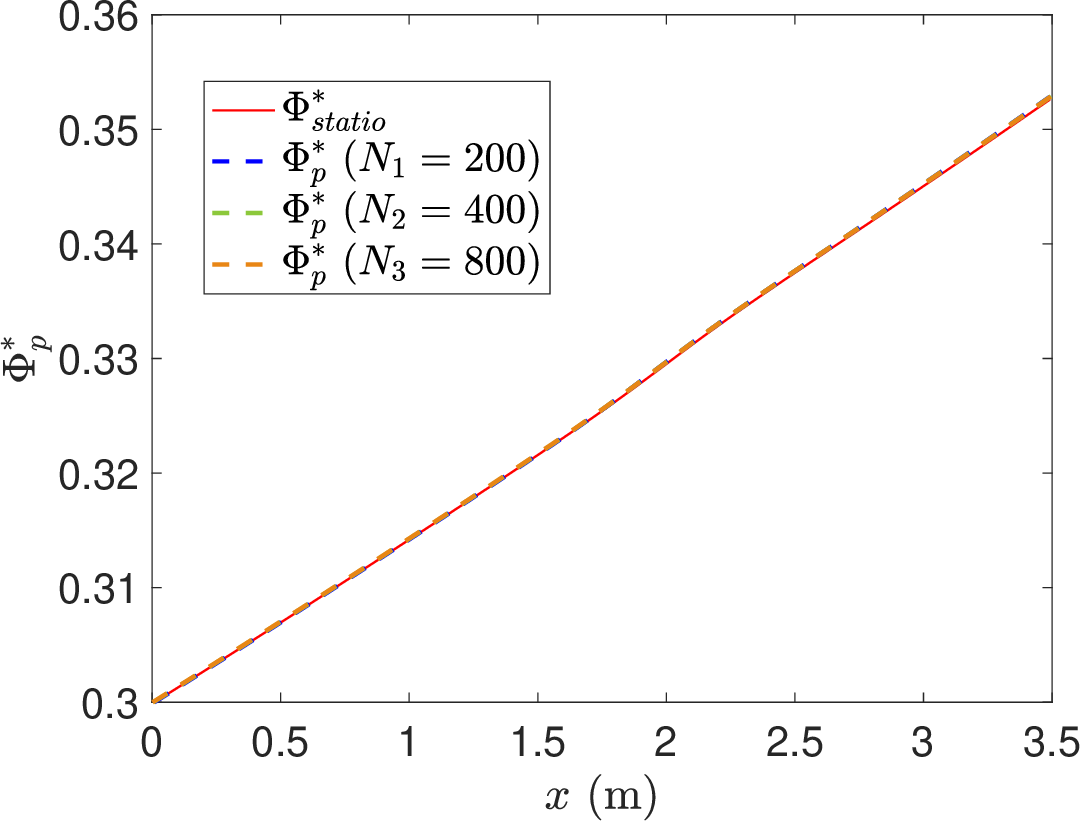}
		\label{fig_statio_friction_smaller_dx_phi}
	\end{subfigure}
	\begin{subfigure}[t]{0.49\textwidth}
		\centering      
		\includegraphics[width=\textwidth]{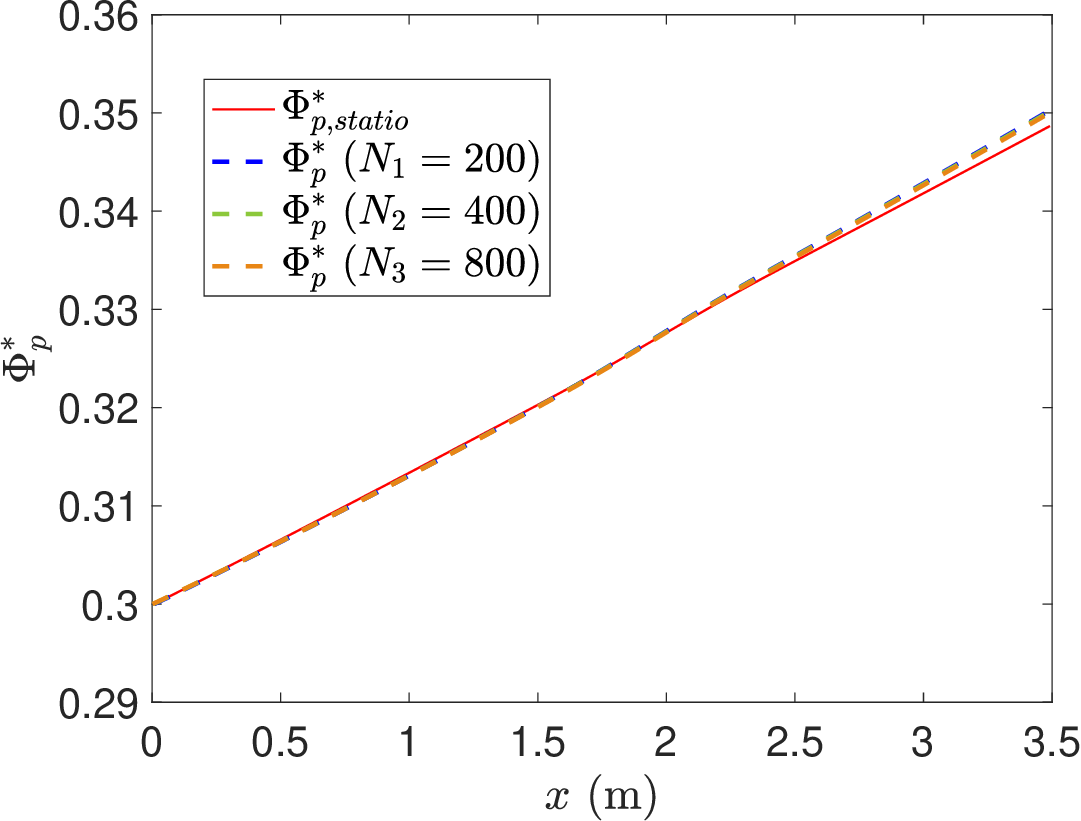}
		\label{fig_statio_friction_smaller_dx_model2_phi}
	\end{subfigure}
	\caption{Test $2$, supercritical regime with bump, $\Phi|_{x=0}=0.3$, $q_0=1$ m$^2$/s. Comparison of solution at time $1\,$s for $N_1=200$ (blue dashed lines), $N_2=400$ (green dashed lines), $N_3=800$ (orange dashed lines), with the stationary solution (continuous red lines).}
	\label{fig_statio_friction_smaller_dx}
\end{figure}

Figure \ref{fig_statio_friction_smaller_dx} shows a stationary solution in a supercritical regime, with three different space step $\Delta_x$ that correspond to discretisations with respectively $N_1=200$, $N_2=400$ and $N_3=800$ points. The red continuous line represents the analytic solution, the three others are the solution computed by the scheme with three different space steps. The error decrease when the space step is smaller (the green and orange curves), which shows the convergence of the scheme for this solution.

\subsection{Test 3: Comparison of stationary solutions}

We compare stationary solutions 
of the two models for different crystals fractions with the same parameters as in test $2$. For the discharge $q$ a Dirichlet condition is set at the upstream boundary: $q|_{x=0} = 1\,$m$^2$/s, and a Neumann condition at the downstream boundary: $\partial_x q|_{x=L} = 0$. Flow height $h$ is fixed at the upstream boundary in the supercritical regime: $h|_{x=0} = 0.1\,$m, and at the downstream boundary in the subcritical regime: $h|_{x=L} = 1.65$ m. Bottom surface $z_b$ is either a tilted plane ($b=0$), or a tilted plane with a bump given by equation \eqref{eq_bump}  (Figure \ref{fig_statio_bottom_z}). \\ \\
\begin{figure}[hbtp]
	\centering
	\begin{subfigure}[t]{0.49\textwidth}
		\centering      
		\includegraphics[width=\textwidth]{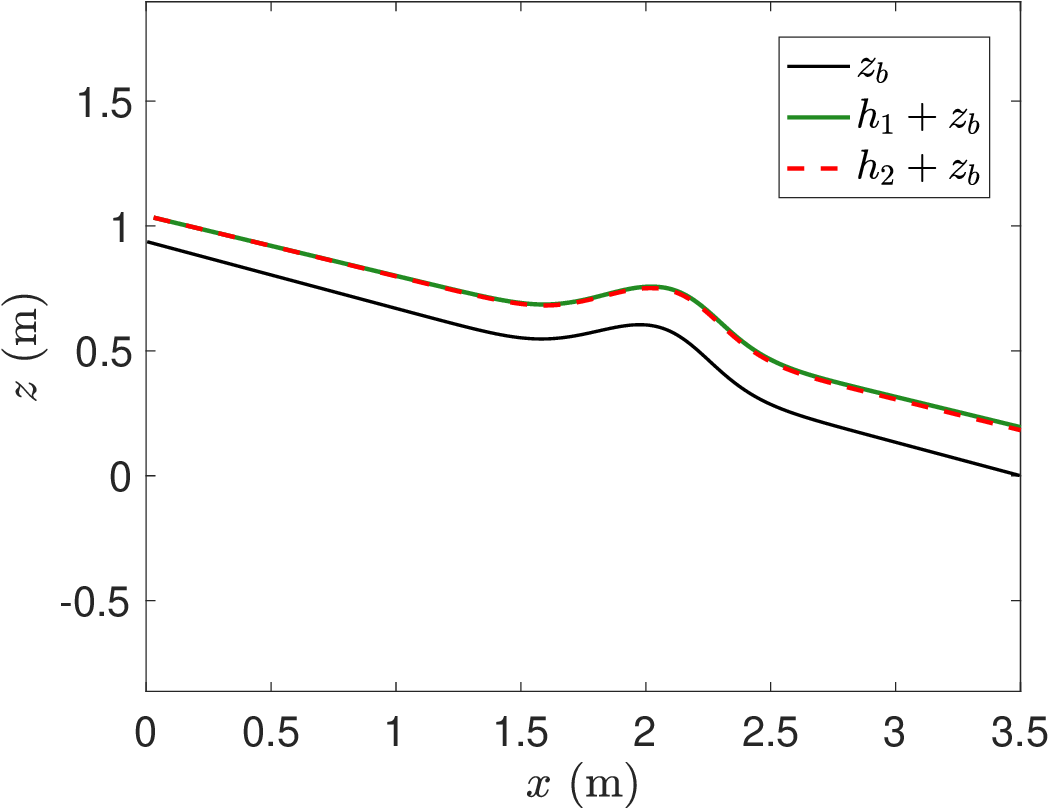}
		\label{fig_statio_compar_heat_lost_h}
	\end{subfigure}
	\begin{subfigure}[t]{0.49\textwidth}
		\centering      
		\includegraphics[width=\textwidth]{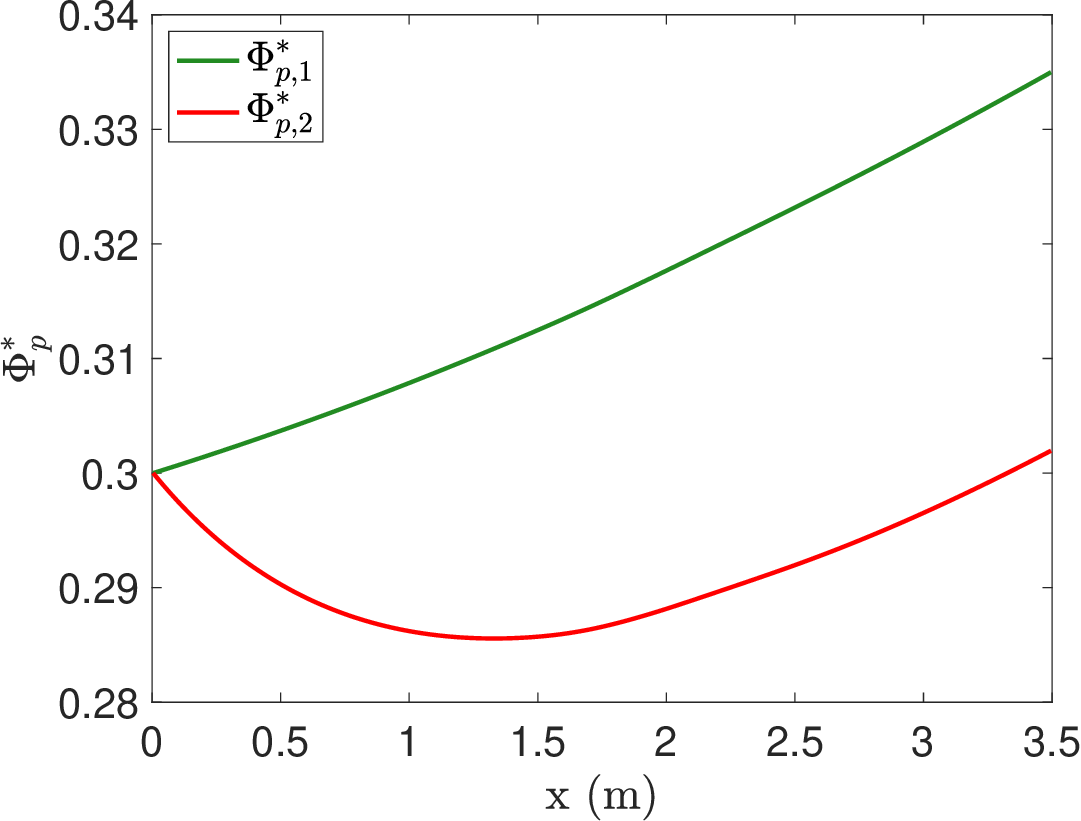}
		\label{fig_statio_compar_heat_lost_phi}
	\end{subfigure}\\
	\begin{subfigure}[t]{0.49\textwidth}
		\centering      
		\includegraphics[width=\textwidth]{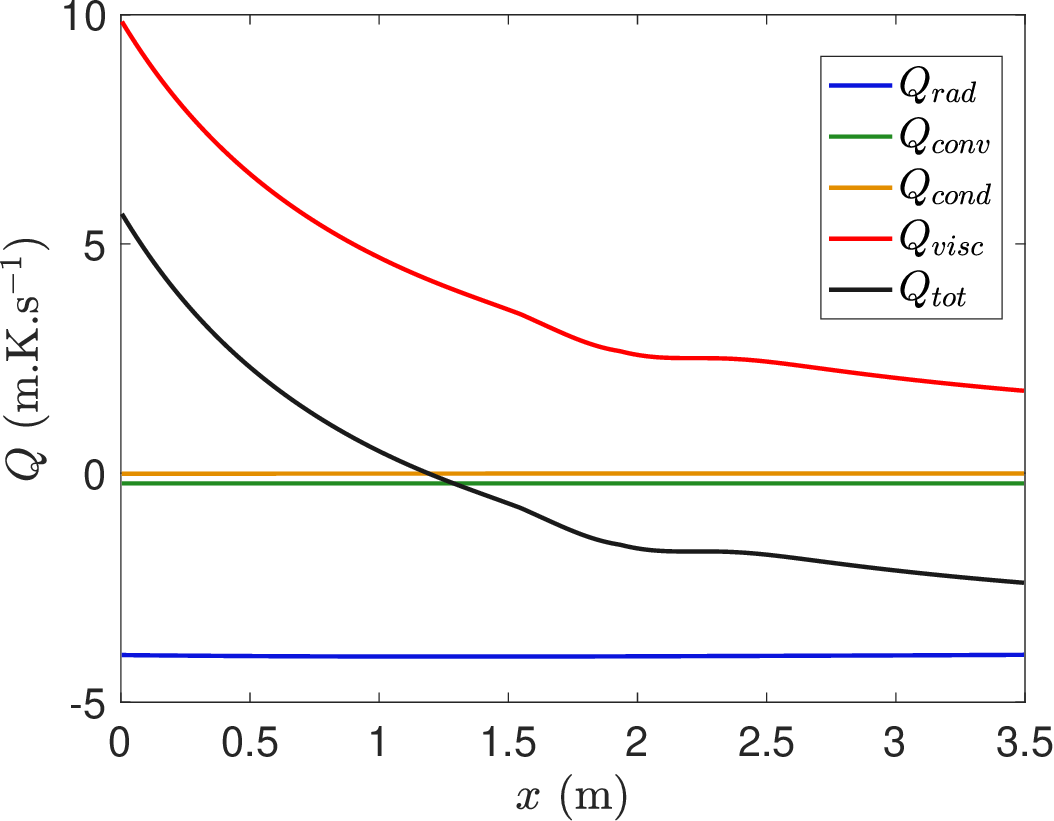}
		\label{fig_statio_compar_heat_lost_Q}
	\end{subfigure}
	\begin{subfigure}[t]{0.49\textwidth}
		\centering
		\includegraphics[width=\textwidth]{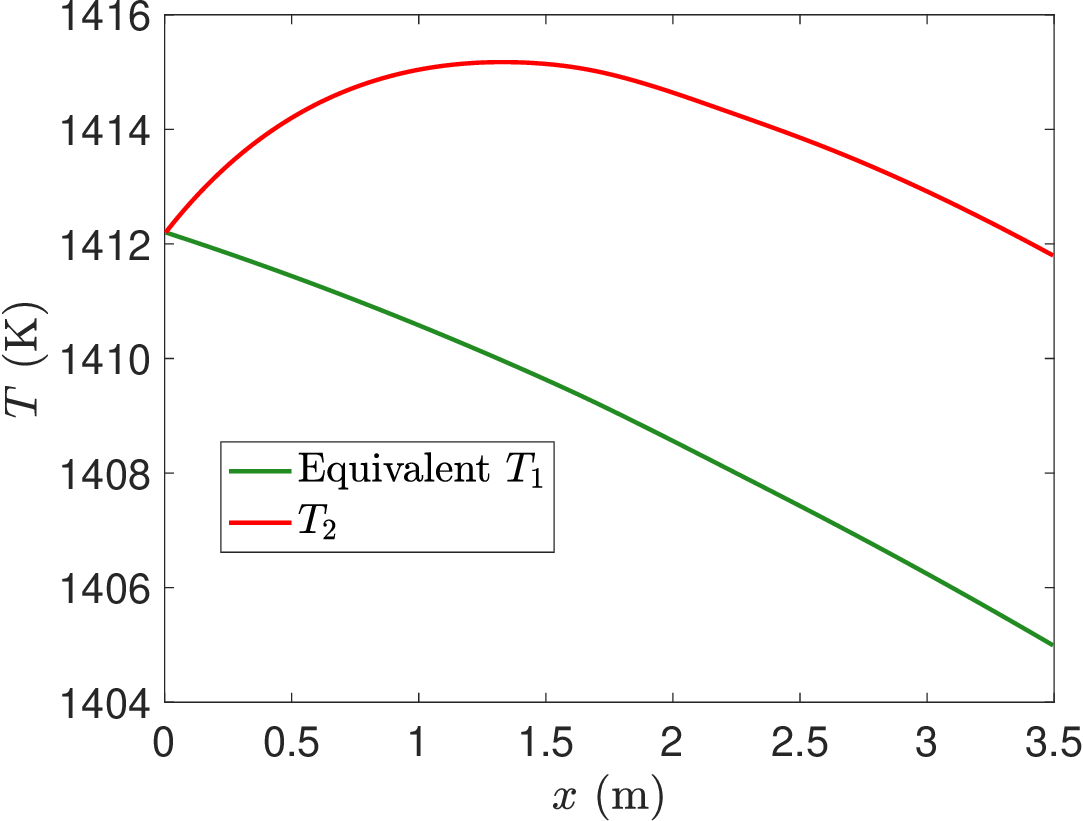}
		\label{fig_statio_compar_heat_lost_T}
	\end{subfigure}
	\caption{Test $3$, supercritical regime with bump, $h|_{x=0}=0.1$, $q_0=1$ m$^2$/s, $\Phi|_{x=0}=0.3$, $N=400$. Comparison of solution of Models $1$ and $2$, at time $1\,$s. In the panel with the comparison of heat lost terms, $Q_{tot} = Q_{rad}+Q_{cond}+Q_{cond}+Q_{visc}$.}
	\label{fig_statio_compar_heat_lost}
\end{figure}

Figure \ref{fig_statio_compar_heat_lost} shows stationary solutions of Model $1$ (green lines), and Model $2$ (red lines). The differences between the two solutions are caused by the friction, as the solid volume fraction $\Phi_p$ evolves differently in each model. With Model $1$, $\Phi_p$ is always increasing because the source term in equation \eqref{eq_phi_p_final} is positive (as $\Phi_p < 1$). On the other hand, in Model $2$, the sign of source term in equation on $T$ \eqref{eq_rheology_phi_version2} can change. In figure \ref{fig_statio_compar_heat_lost} we can see that this term ($Q_{tot}$) is positive when $x<1$, and negative when $x>1.5$, which causes the change of variation of $\Phi_p$ and $T$. The radiative and convective term are negative (as long as $T>T_{env}$), and the viscous term is always positive. The predominance of this viscous heating term for small $x$ causes a positive total heat loss. \\ \\
\begin{figure}[hbtp]
	\centering
    \begin{subfigure}[t]{0.49\textwidth}
		\centering      
		\includegraphics[width=\textwidth]{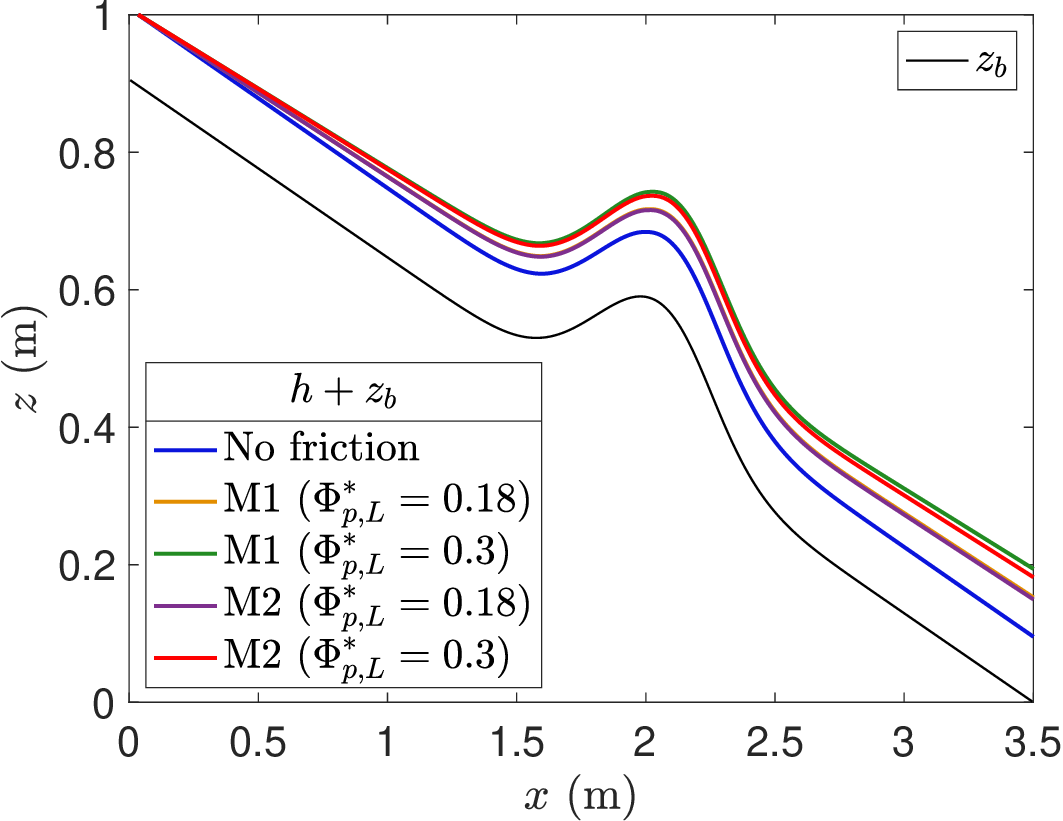}
		\label{fig_statio_comparfriction_super_h2}
	\end{subfigure}
	\begin{subfigure}[t]{0.49\textwidth}
		\centering      
		\includegraphics[width=\textwidth]{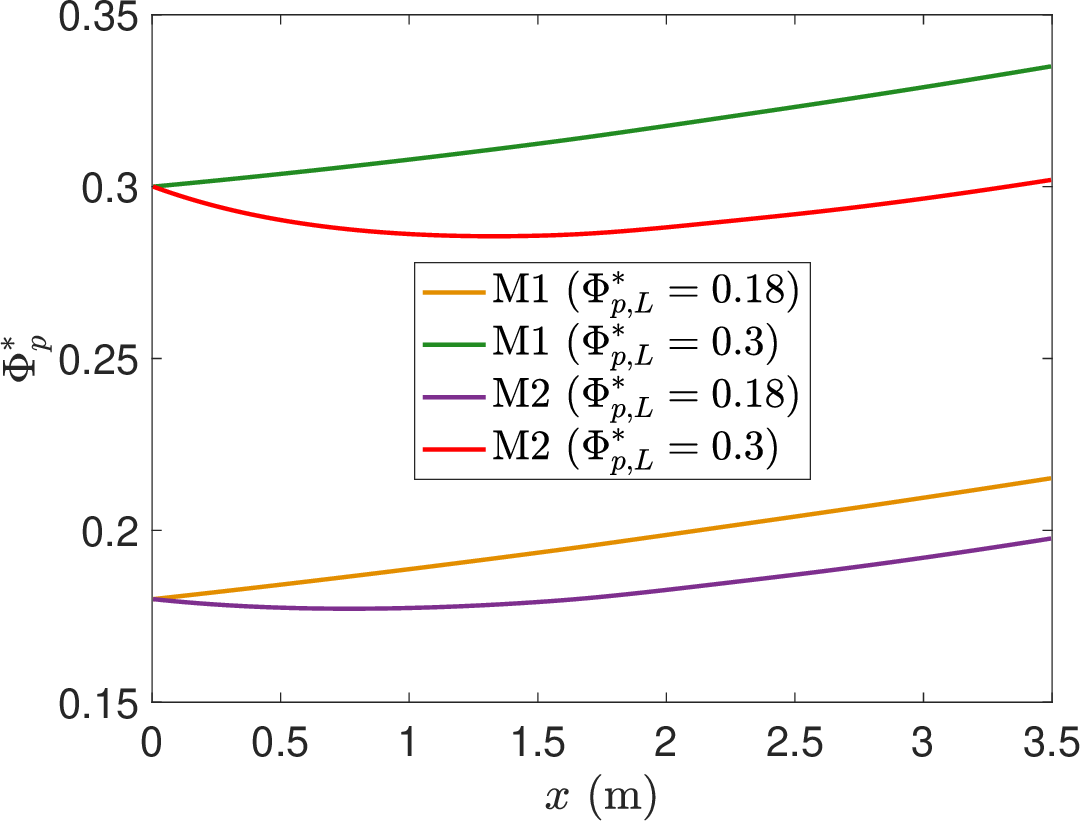}
		\label{fig_statio_comparfriction_super_phi2}
	\end{subfigure}
	\caption{Test $3$, Models $1$ and $2$, supercritical regime, $h|_{x=0}=0.1$, $q_0=1$, $N=400$. Comparison of solutions at time $1\,$s for different friction: no friction, $\Phi_p^*|_{x=0} = \Phi_{p,L}^* = 0.18$, $\Phi_p^*|_{x=0} = \Phi_{p,L}^* = 0.3$.}
	\label{fig_statio_comparfriction_super}
\end{figure}
\begin{figure}[hbtp]
	\centering
  	\begin{subfigure}[t]{0.49\textwidth}
		\centering      
		\includegraphics[width=\textwidth]{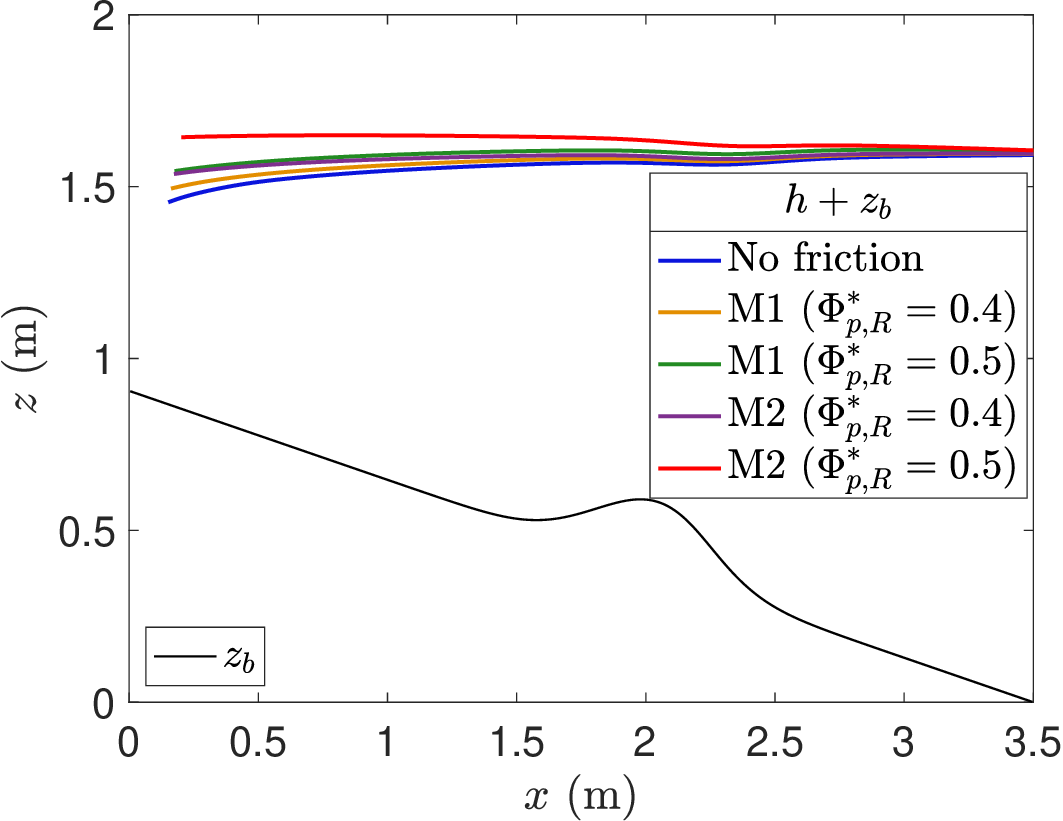}
		\label{fig_statio_comparfriction_sub_h2}
	\end{subfigure}
	\begin{subfigure}[t]{0.49\textwidth}
		\centering      
		\includegraphics[width=\textwidth]{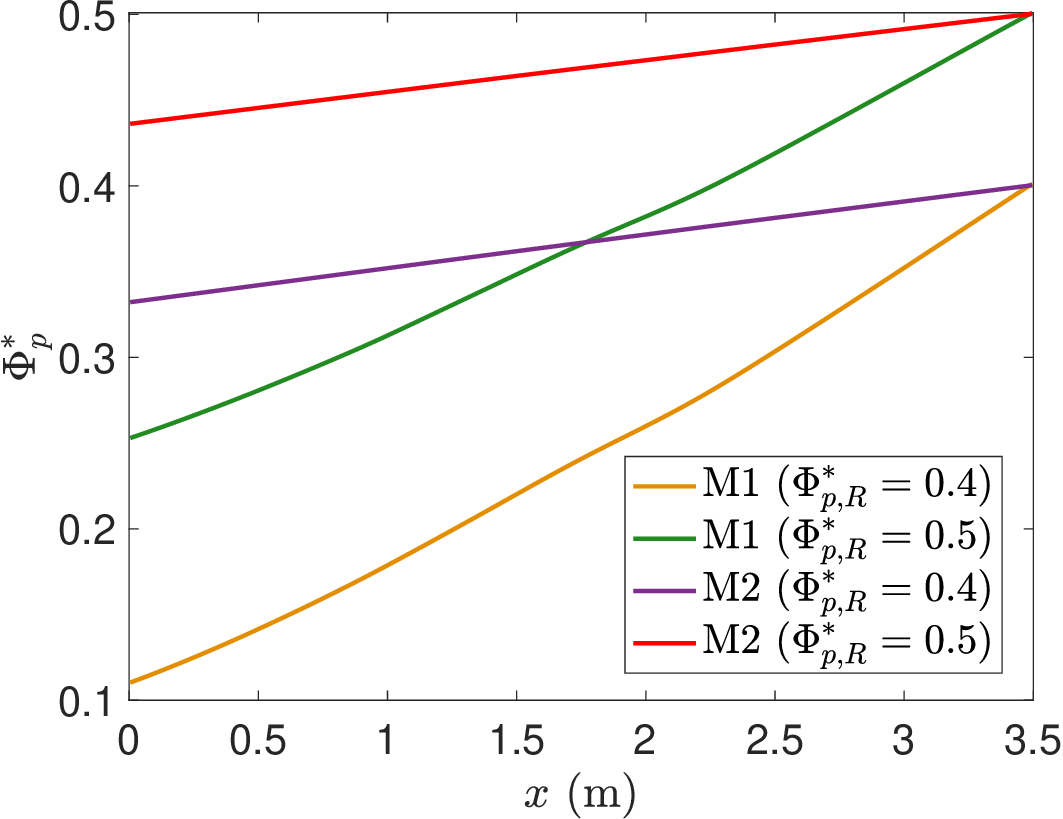}
		\label{fig_statio_comparfriction_sub_phi2}
	\end{subfigure}
	\caption{Test $3$, Models $1$ and $2$, subcritical regime, $h|_{x=L} = 1.65$, $q_0=1$ m$^2$/s, $N=400$. Comparison of solutions at time $1\,$s for different friction: no friction, $\Phi_p^*|_{x=L} = \Phi_{p,R}^* = 0.4$, $\Phi_p^*|_{x=L} = \Phi_{p,R}^* = 0.5$.}
	\label{fig_statio_comparfriction_sub}
\end{figure}

\noindent Next we compare solutions for different initial conditions on $\Phi_p^*$. In figures \ref{fig_statio_comparfriction_super}-\ref{fig_statio_comparfriction_sub} the blue curves are the solutions for the model without friction (classical shallow water equations) at time $5\,$s. This solution is compared with those of the complete model \eqref{eq_syst1D_closed}-\eqref{eq_rheology_phi_version2}, with two different boundary conditions for $\Phi_p$. The orange and purple curves are respectively the solutions of Models $1$ and $2$ when $\Phi_0 = 0.18$. The green and red curves are respective solutions of the models when $\Phi_p|_{x=0} = 0.3$.
In figure \ref{fig_statio_comparfriction_super} the solutions are supercritical and the minimum Froude number is $10.3$ without friction and $2.65$ with friction. The increase of friction leads to a downstream thickening of the flow. Figure \ref{fig_statio_comparfriction_sub} presents subcritical solutions. The Froude number has a maximum value of $0.7$ without friction and $0.56$ with friction. The fluid height $h$ is fixed at the downstream boundary, and the free surface $z_b+h$ tends to flatten at higher friction.\\ \\
All these numerical solutions behave well, there are no spurious oscillations caused by instabilities of the scheme, and the solutions converge to stationary states.


\subsection{Test 4: Dam break without wet/dry front}

\begin{figure}[hbtp]
	\centering
	\begin{subfigure}[t]{0.49\textwidth}
		\centering      
	    \textbf{Model $1$}\par\medskip
		\includegraphics[width=\textwidth]{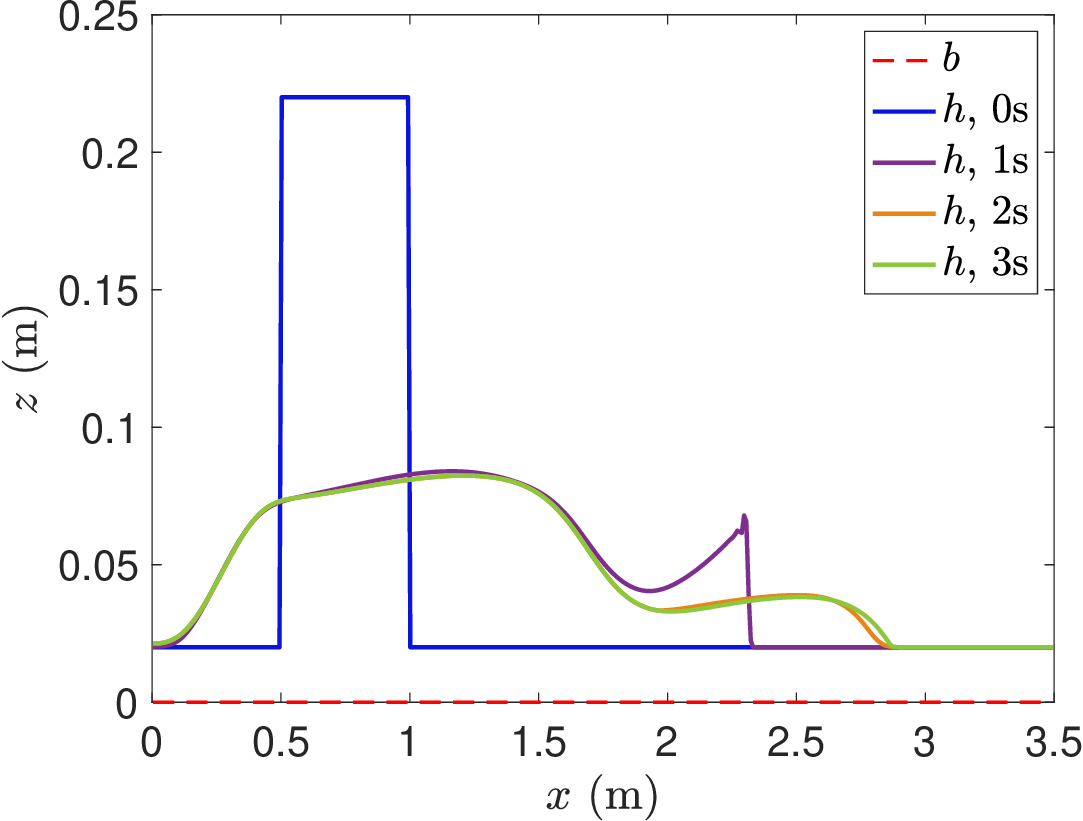}
		\label{fig_dambreak_nodry_model1_h}
	\end{subfigure}
	\begin{subfigure}[t]{0.49\textwidth}
		\centering      
	    \textbf{Model $2$}\par\medskip
		\includegraphics[width=\textwidth]{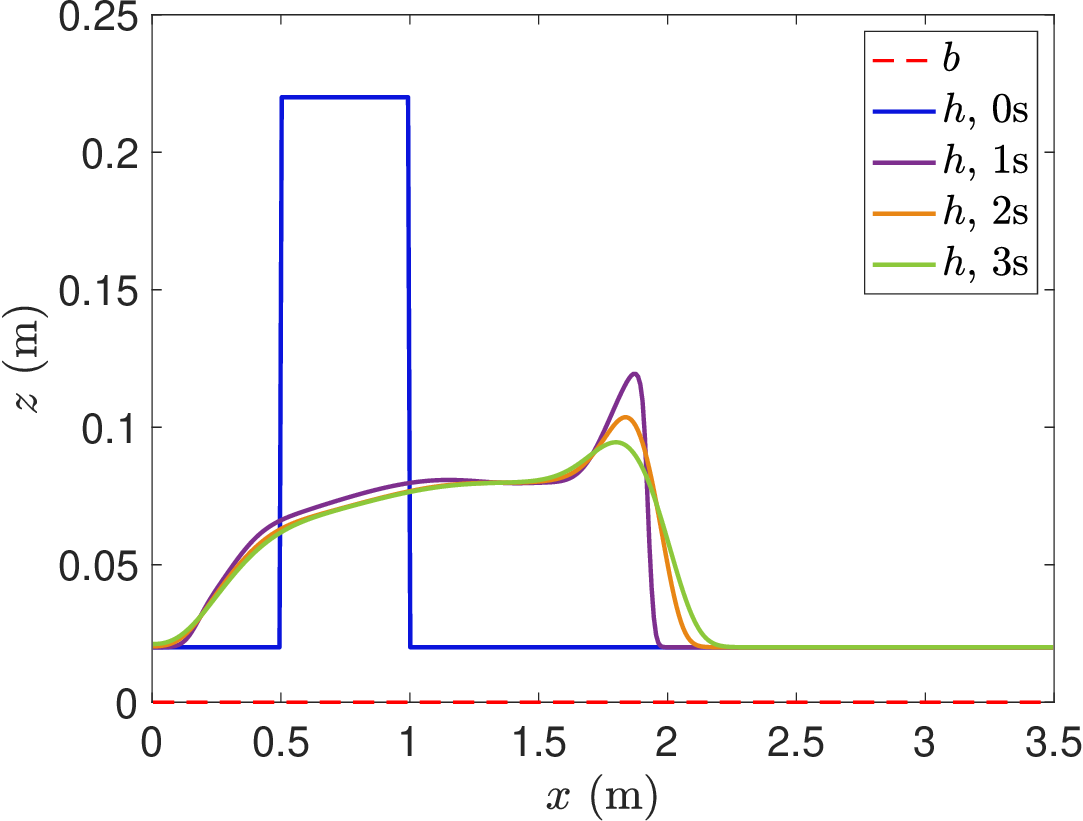}
		\label{fig_dambreak_nodry_model2_h}
	\end{subfigure}\\[1mm]
	\begin{subfigure}[t]{0.49\textwidth}
		\centering      
		\includegraphics[width=\textwidth]{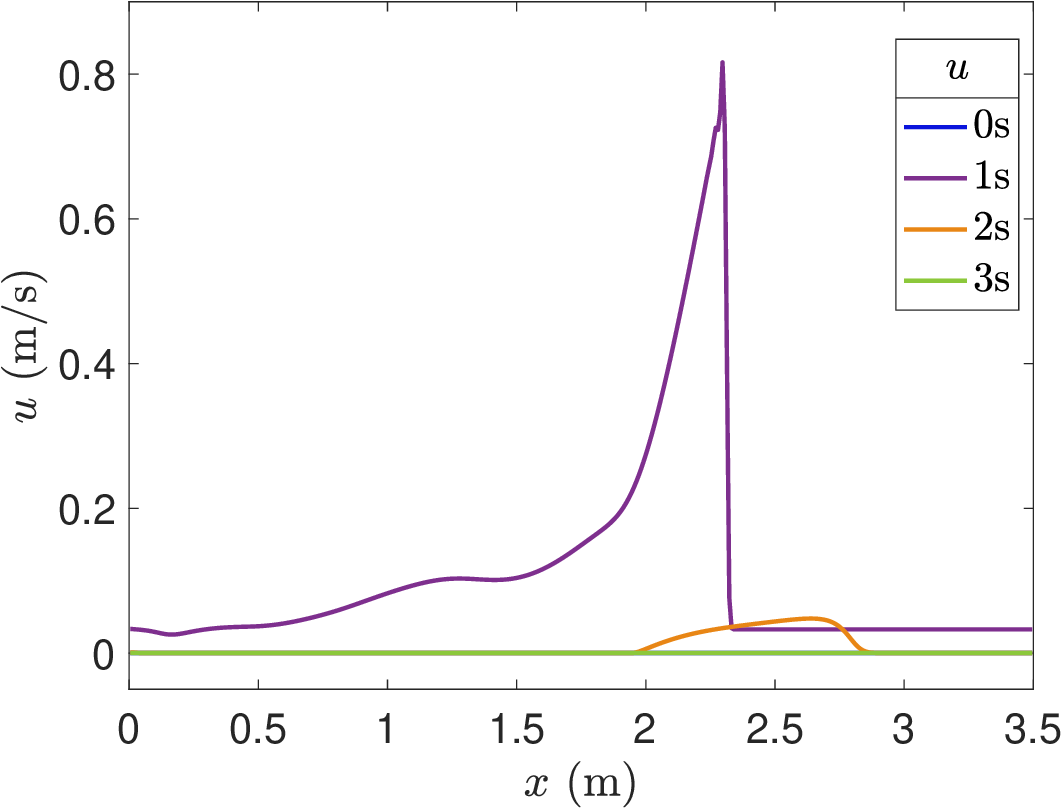}
		\label{fig_dambreak_nodry_model1_u}
	\end{subfigure}
	\begin{subfigure}[t]{0.49\textwidth}
		\centering      
		\includegraphics[width=\textwidth]{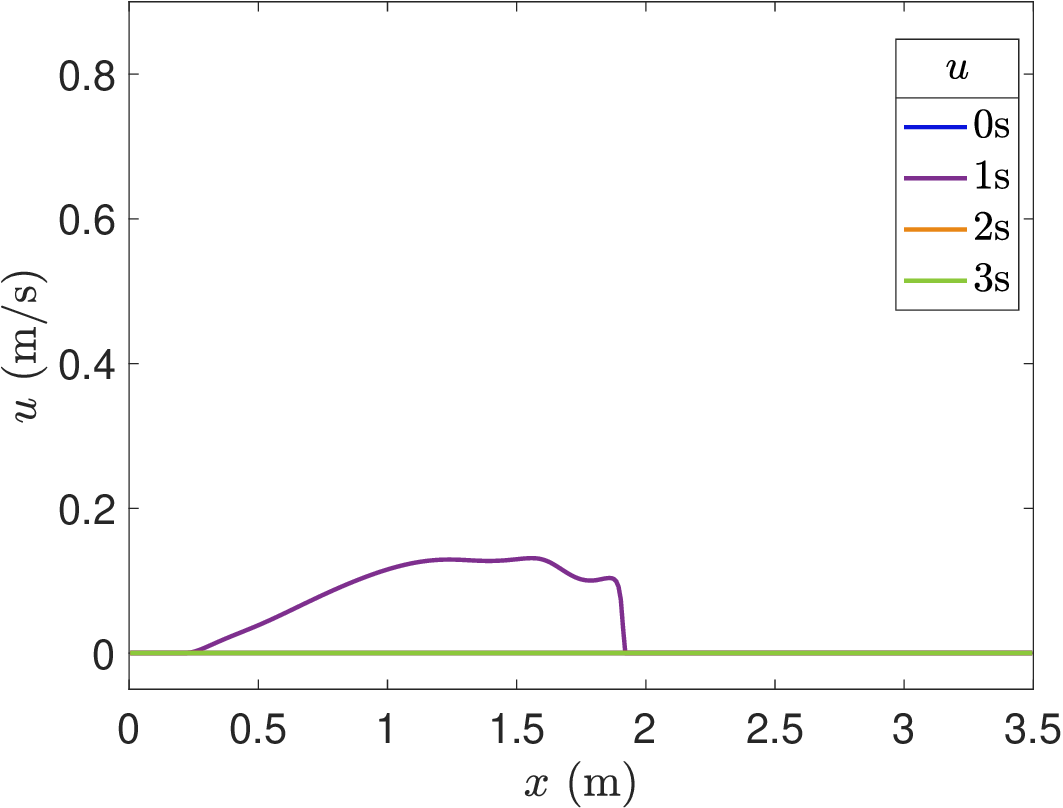}
		\label{fig_dambreak_nodry_model2_u}
	\end{subfigure}\\[1mm]
	\begin{subfigure}[t]{0.49\textwidth}
		\centering      
		\includegraphics[width=\textwidth]{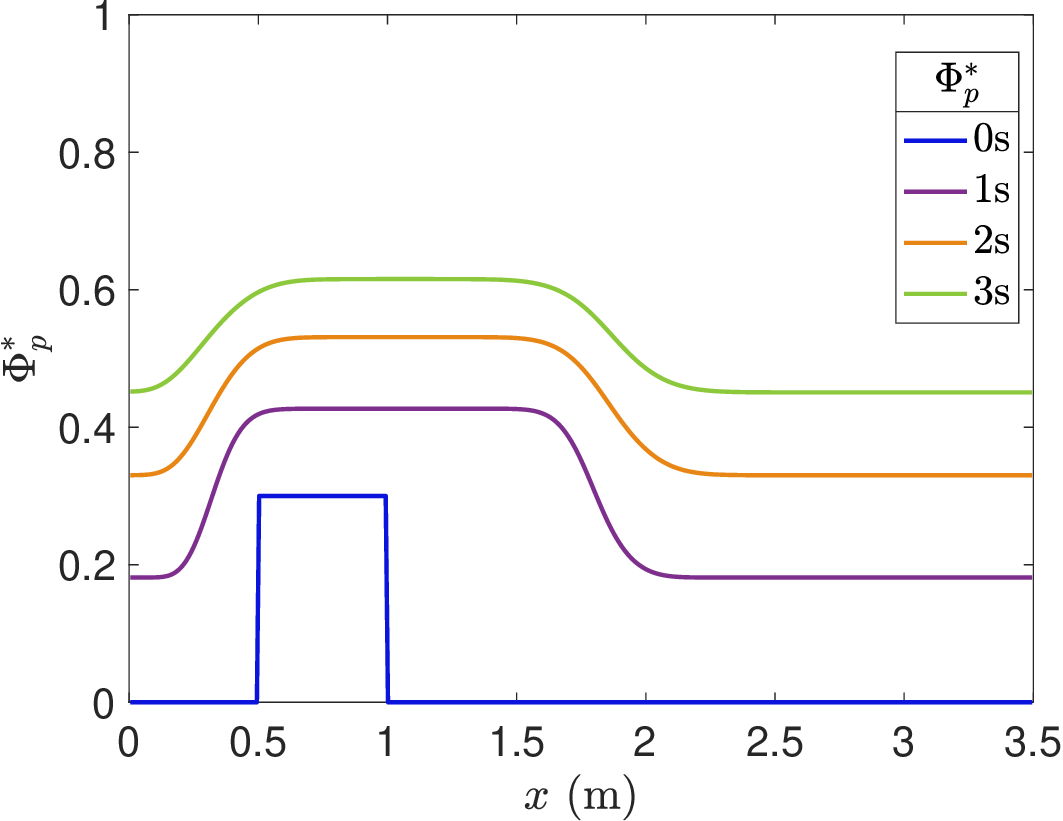}
		\label{fig_dambreak_nodry_model1_phi}
	\end{subfigure}
	\begin{subfigure}[t]{0.49\textwidth}
		\centering      
		\includegraphics[width=\textwidth]{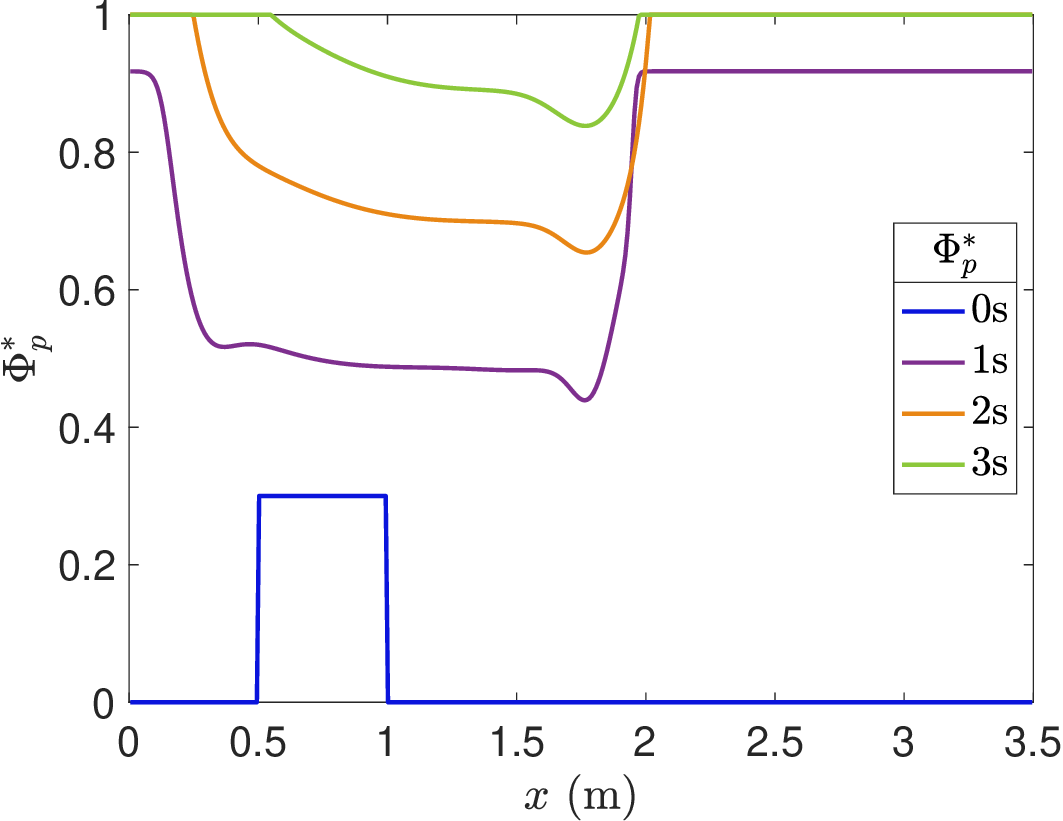}
		\label{fig_dambreak_nodry_model2_phi}
	\end{subfigure}
	\caption{Test $4$, $b=0$, dam break. Left Model $1$, right Model $2$. $\tau_{crys} = 0.5\,$s, $\Phi_p|_{t=0} = 0.3$.}
	\label{fig_dambreak_nodry}
\end{figure}

This test consists in Models $1$ and $2$ simulations of a dam break on an inclined plane without wet/dry front. The numerical parameters are given in Table \ref{tab_paramnum_simu_nonuniform}, and the model parameters in Table \ref{tab_param_test_scheme} (except for the crystallization time, which is $\tau_{crys}=0.5 \,$s). The initial amount of fluid in the domain $[0,3.5]$ is
\begin{align*}
    h(x)|_{t=0} = \left\{ \begin{array}{ll}
        \displaystyle 0.2 &\text{if } 0.5 \leq x \leq 1,\\
        \displaystyle 0.02 & \text{else}.
    \end{array} \right.
\end{align*}
The initial discharge is zero, and the initial crystal fraction is given by:
\begin{align*}
    \Phi_p(x)|_{t=0} = \left\{ \begin{array}{ll}
        \displaystyle 0.3 & \text{if } 0.5 \leq x \leq 1,\\
        \displaystyle 0 & \text{else}.
    \end{array} \right.
\end{align*}
Neumann boundary conditions are used at the top and bottom boundary of the domain.\\ \\
Results for a flat bottom $b=0$ are shown in Figure \ref{fig_dambreak_nodry}, with Model $1$ on the left panel and Model $2$ on the right panel. The crystal fraction increases more with Model $2$ where $h$ is smaller, and increases less at larger $h$ values. This is due to the friction term that is proportional to $1/h$ in Model $2$, whereas in Model $1$ the fluid height does not influence crystallization. As expected, the lava flow down the slope and its velocity decrease as crystal fraction increases. The lava velocity vanishes completely when the crystal fraction reaches a critical level.

\subsection{Test 5: Dam break with wet/dry front}

\begin{figure}[hbtp]
	\centering
	\begin{subfigure}[t]{0.49\textwidth}
		\centering      
	    \textbf{Model $1$}\par\medskip
		\includegraphics[width=\textwidth]{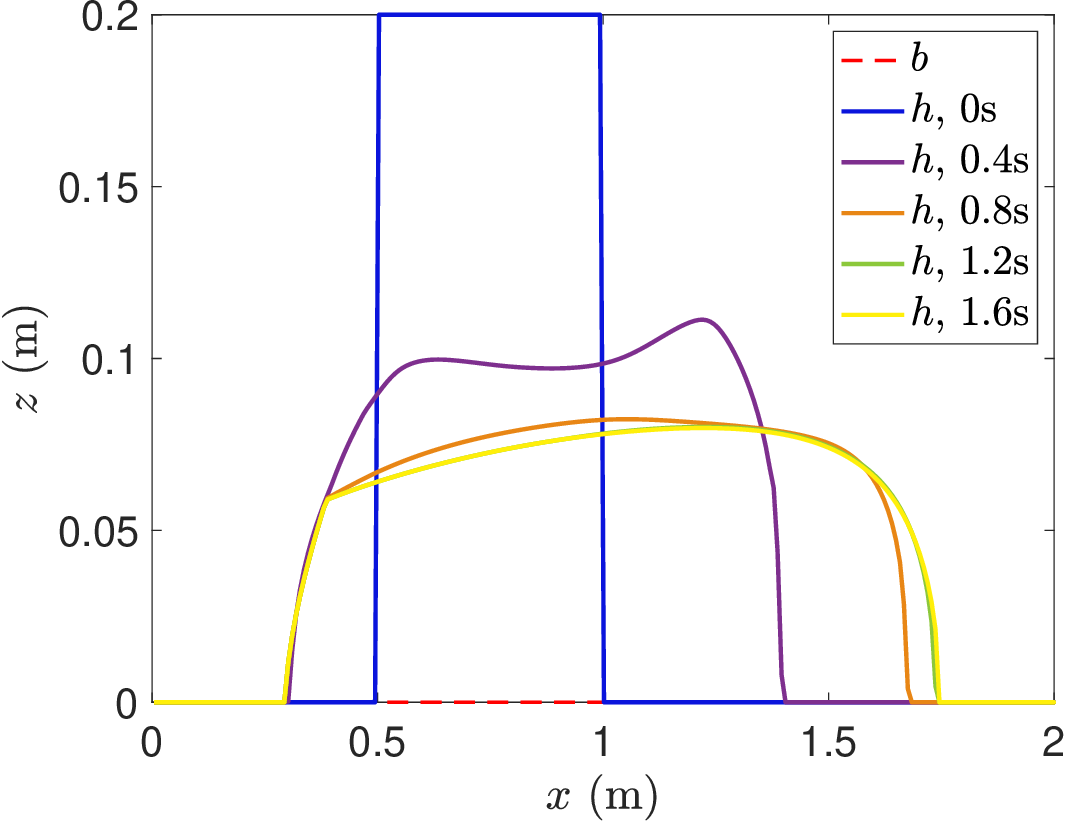}
		\label{fig_dambreak0_model1_h}
	\end{subfigure}
	\begin{subfigure}[t]{0.49\textwidth}
		\centering      
	    \textbf{Model $2$}\par\medskip
		\includegraphics[width=\textwidth]{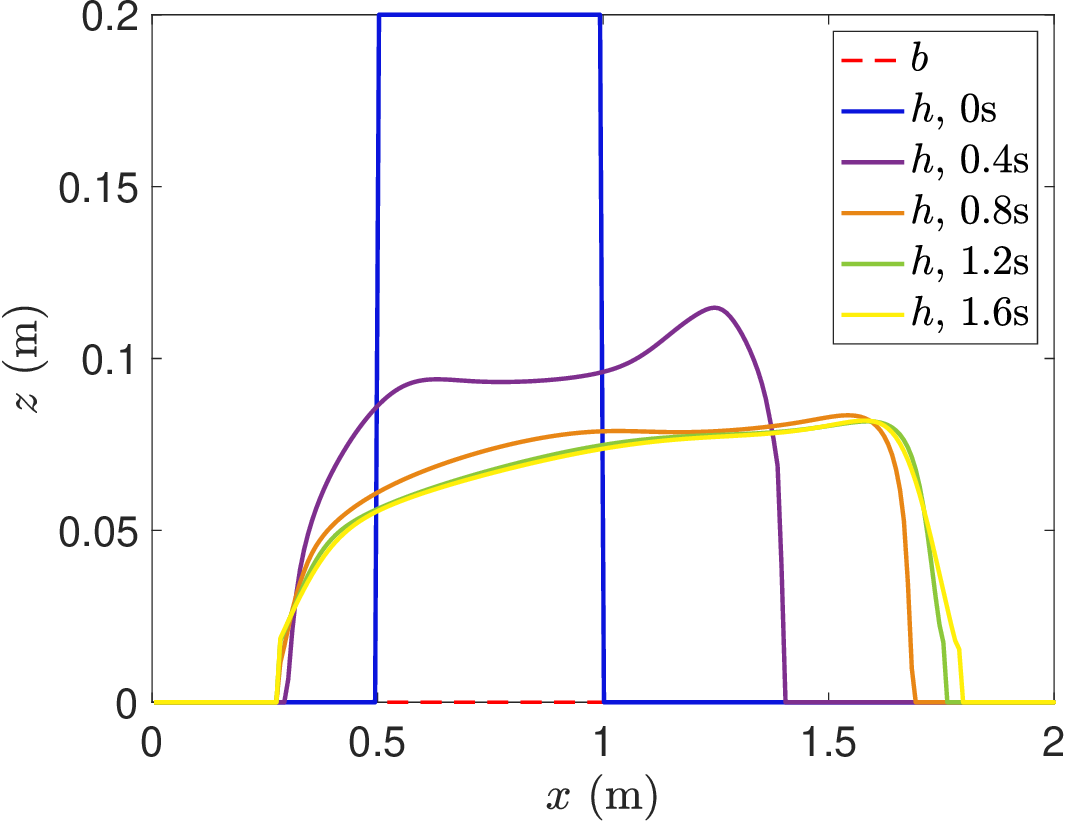}
		\label{fig_dambreak0_model2_h}
	\end{subfigure}\\[1mm]
	\begin{subfigure}[t]{0.49\textwidth}
		\centering      
		\includegraphics[width=\textwidth]{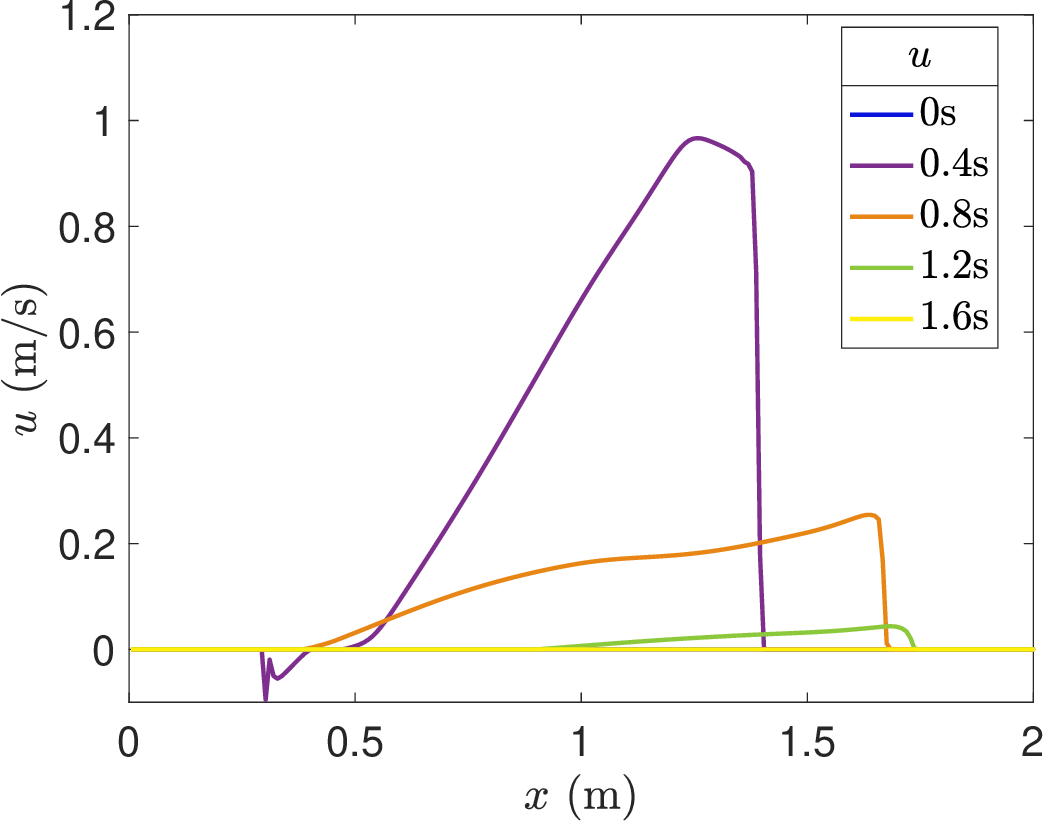}
		\label{fig_dambreak0_model1_u}
	\end{subfigure}
	\begin{subfigure}[t]{0.49\textwidth}
		\centering      
		\includegraphics[width=\textwidth]{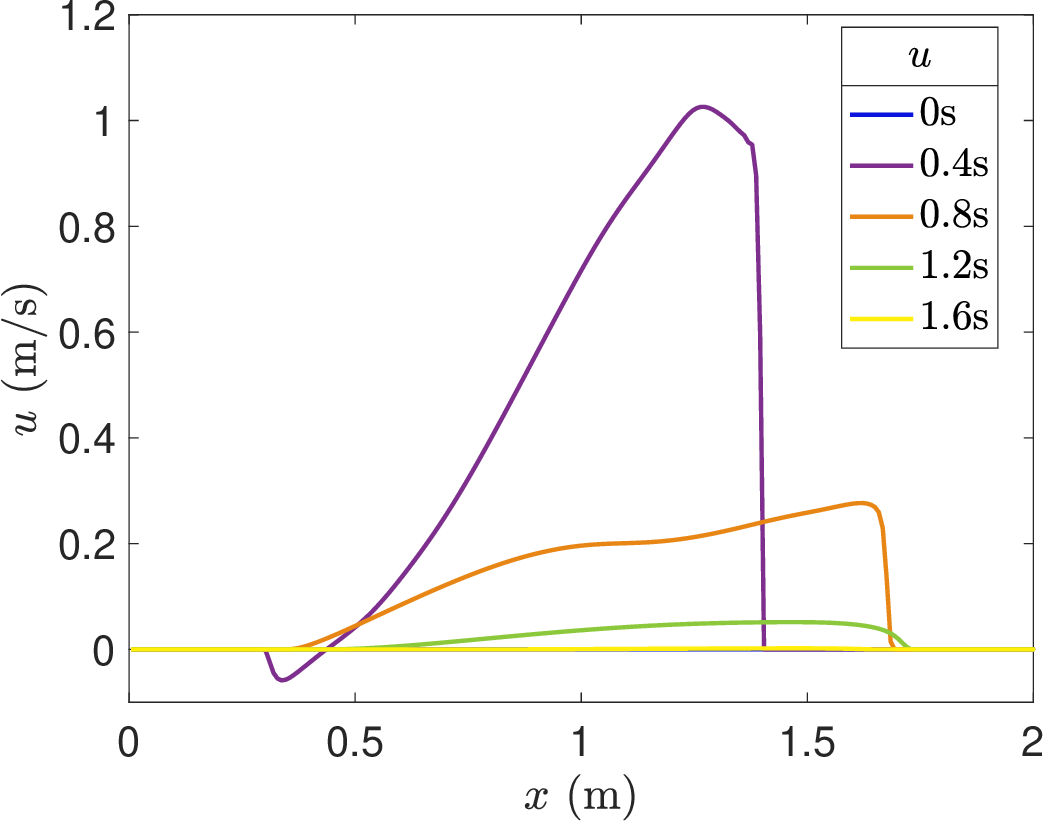}
		\label{fig_dambreak0_model2_u}
	\end{subfigure}\\[1mm]
	\begin{subfigure}[t]{0.49\textwidth}
		\centering      
		\includegraphics[width=\textwidth]{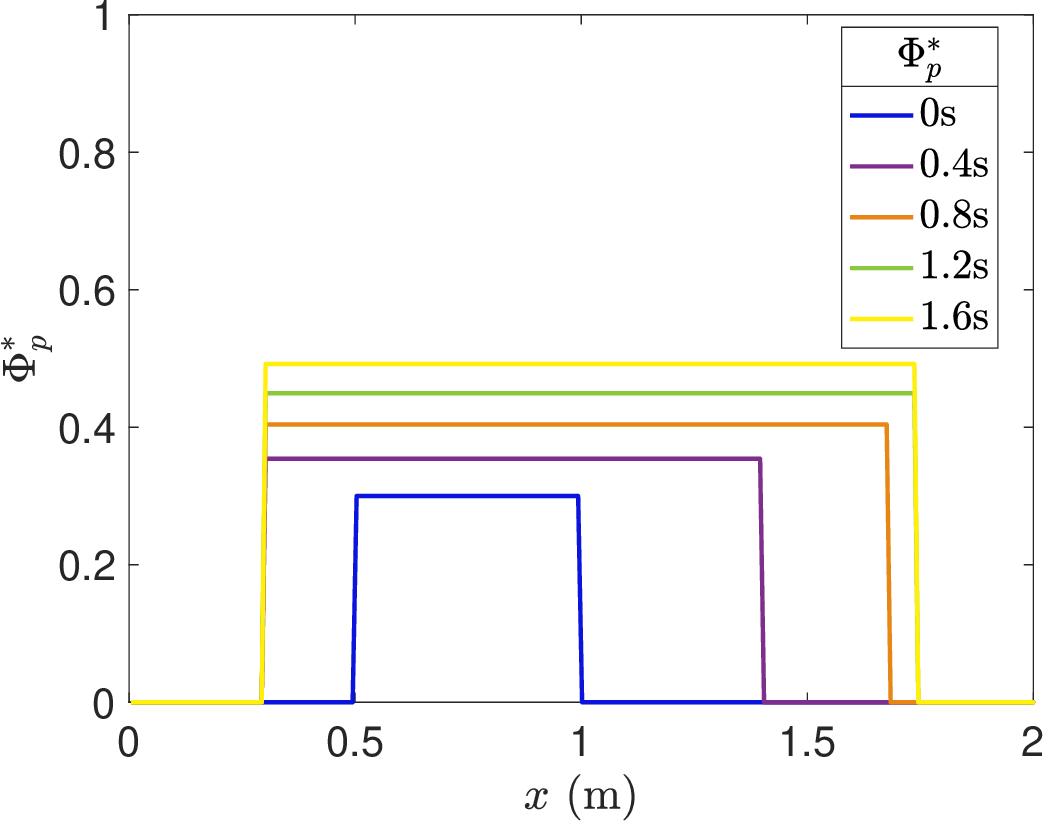}
		\label{fig_dambreak0_model1_phi}
	\end{subfigure}
	\begin{subfigure}[t]{0.49\textwidth}
		\centering      
		\includegraphics[width=\textwidth]{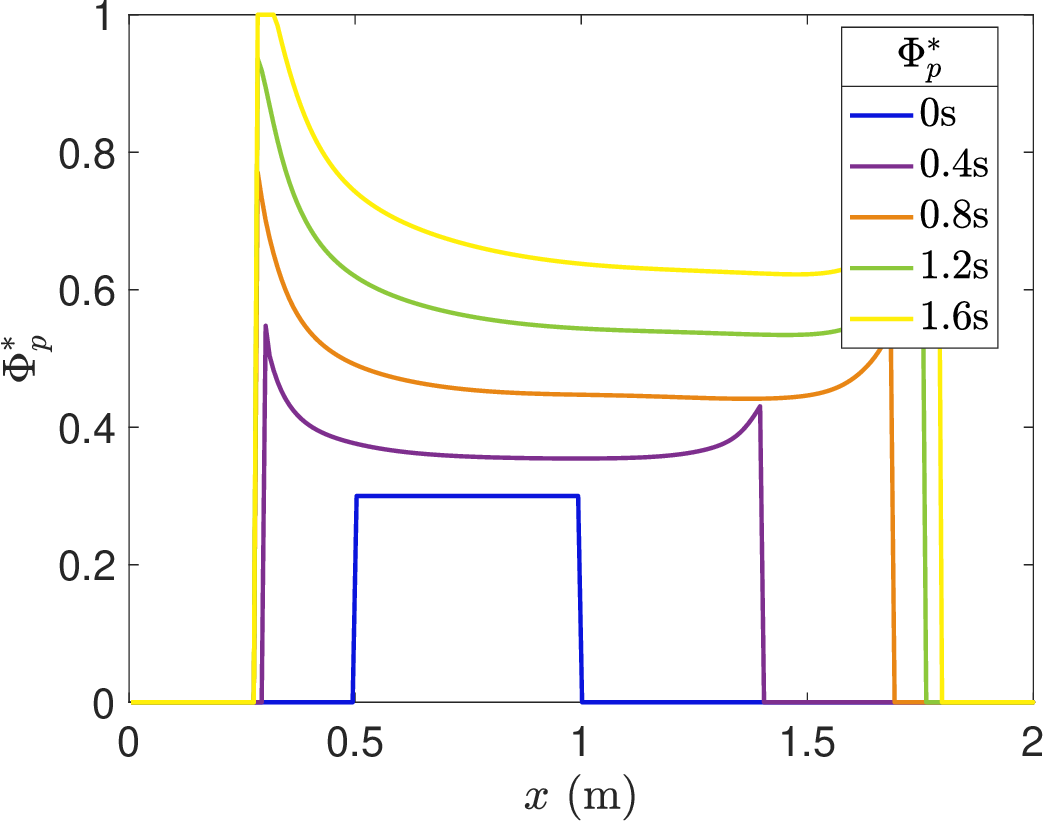}
		\label{fig_dambreak0_model2_phi}
	\end{subfigure}
	\caption{Test $5$, $b=0$, dam break. Left Model $1$, right Model $2$. $\tau_{crys} = 0.5\,$s, $\Phi_p|_{t=0} = 0.3$.}
	\label{fig_dambreak0}
\end{figure}

\begin{figure}[hbtp]
	\centering
	\begin{subfigure}[t]{0.49\textwidth}
	    \centering      
	    \textbf{Model $1$}\par\medskip
		\includegraphics[width=\textwidth]{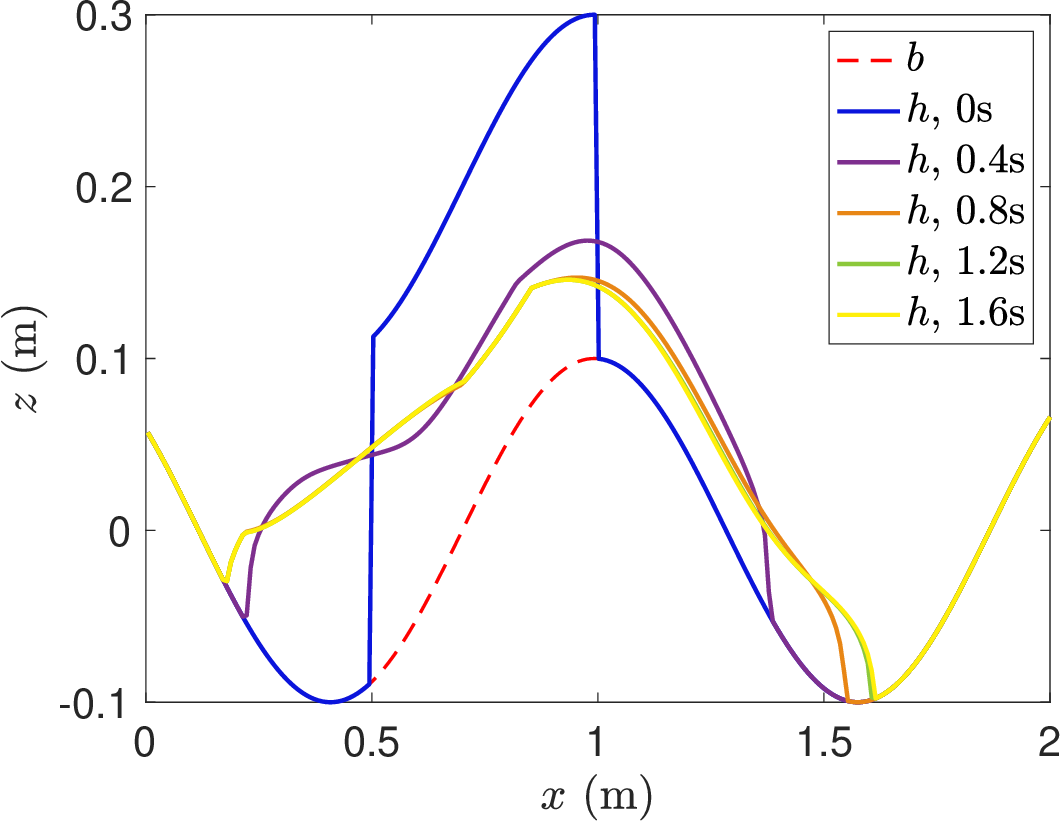}
		\label{fig_dambreak1_model1_h}
	\end{subfigure}
	\begin{subfigure}[t]{0.49\textwidth}
		\centering      
	    \textbf{Model $2$}\par\medskip
		\includegraphics[width=\textwidth]{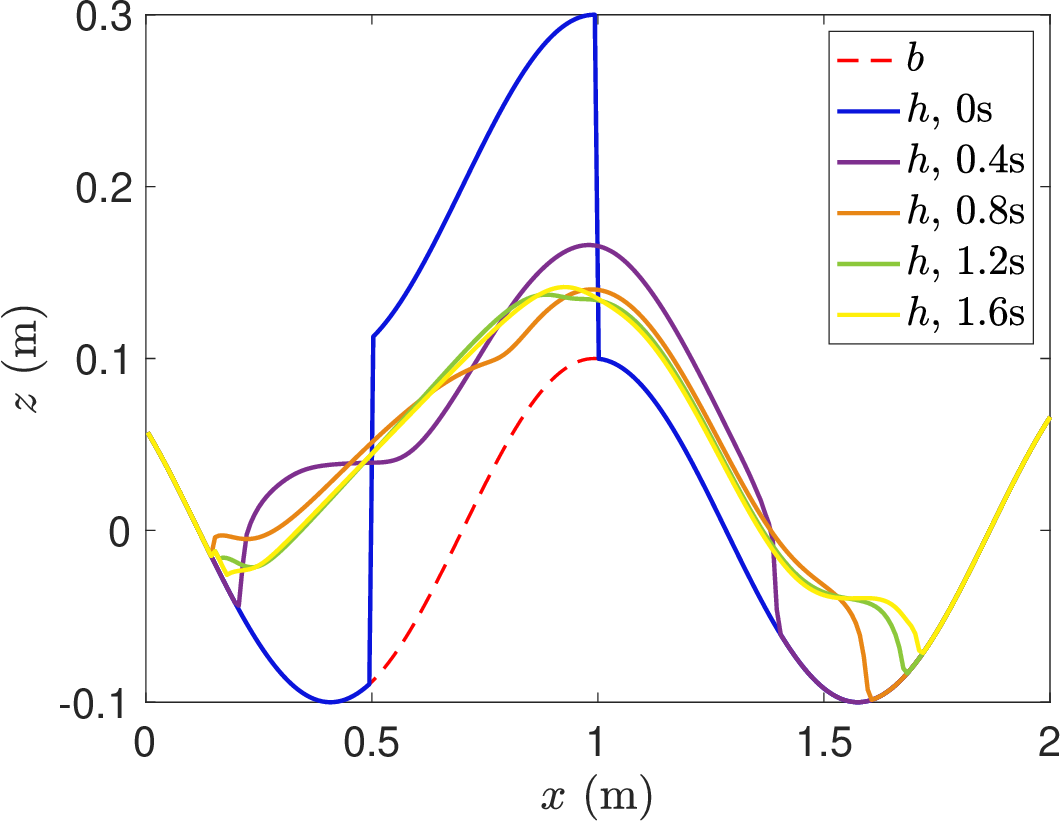}
		\label{fig_dambreak1_model2_h}
	\end{subfigure}\\[1mm]
	\begin{subfigure}[t]{0.49\textwidth}
		\centering      
		\includegraphics[width=\textwidth]{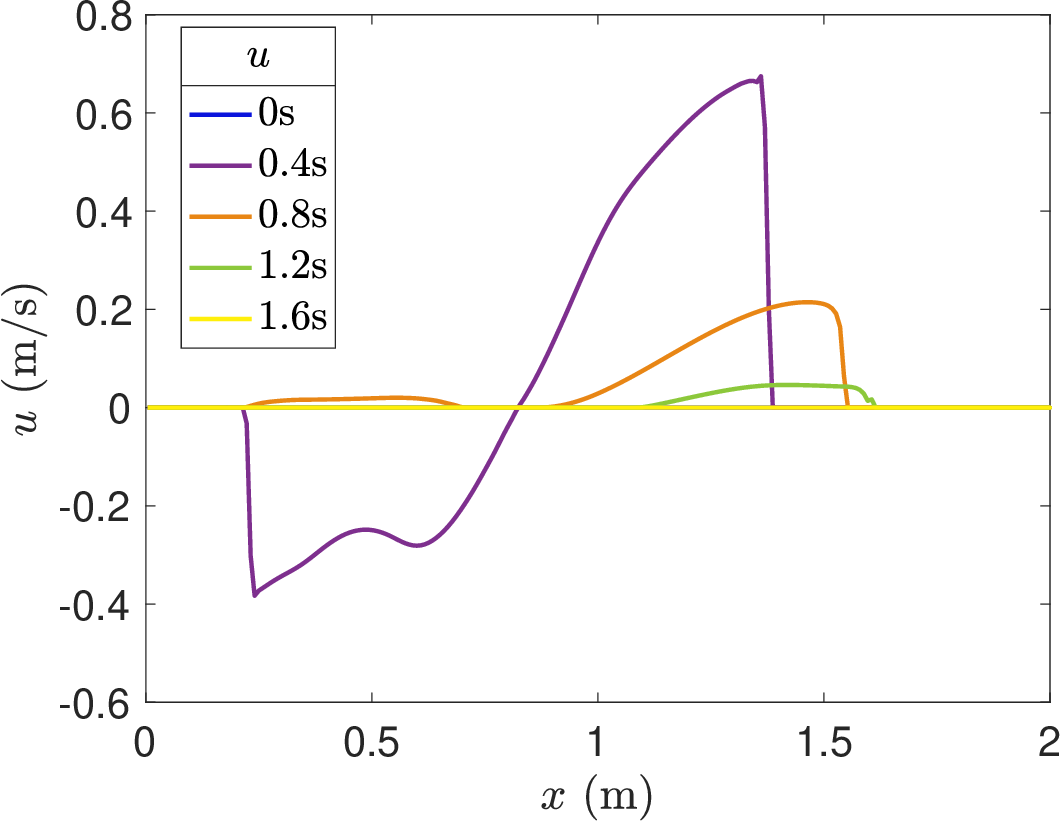}
		\label{fig_dambreak1_model1_u}
	\end{subfigure}
	\begin{subfigure}[t]{0.49\textwidth}
		\centering      
		\includegraphics[width=\textwidth]{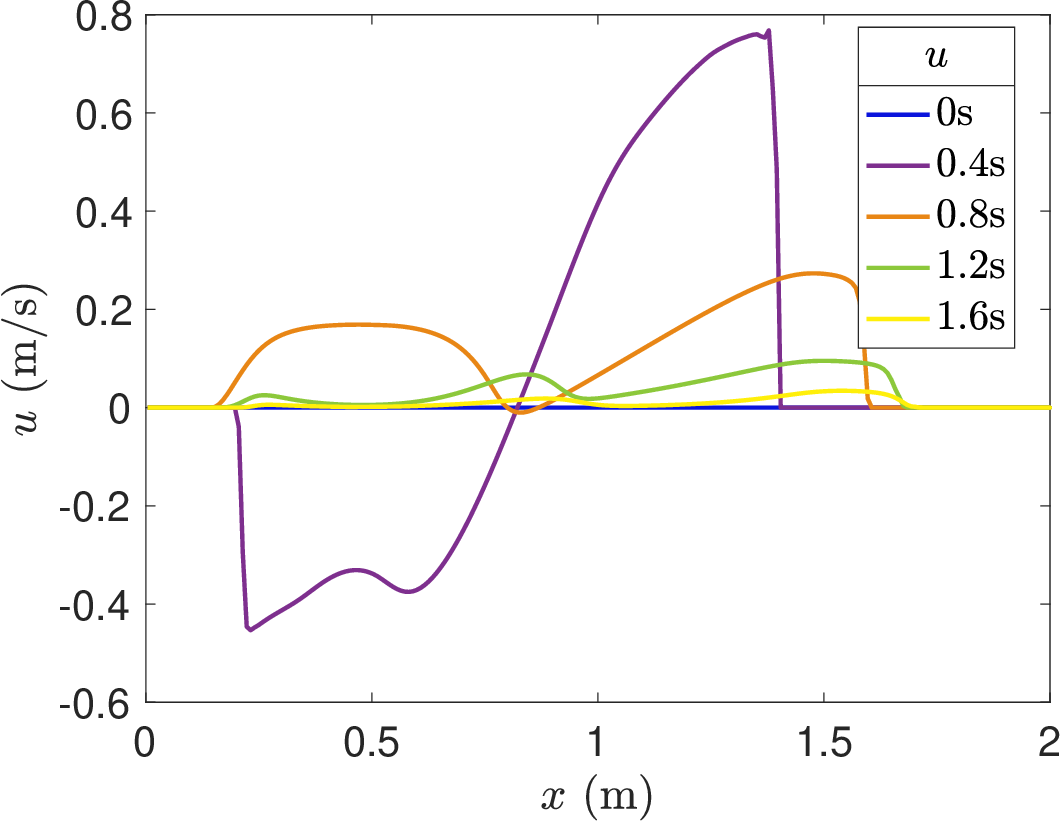}
		\label{fig_dambreak1_model2_u}
	\end{subfigure}\\[1mm]
	\begin{subfigure}[t]{0.49\textwidth}
		\centering      
		\includegraphics[width=\textwidth]{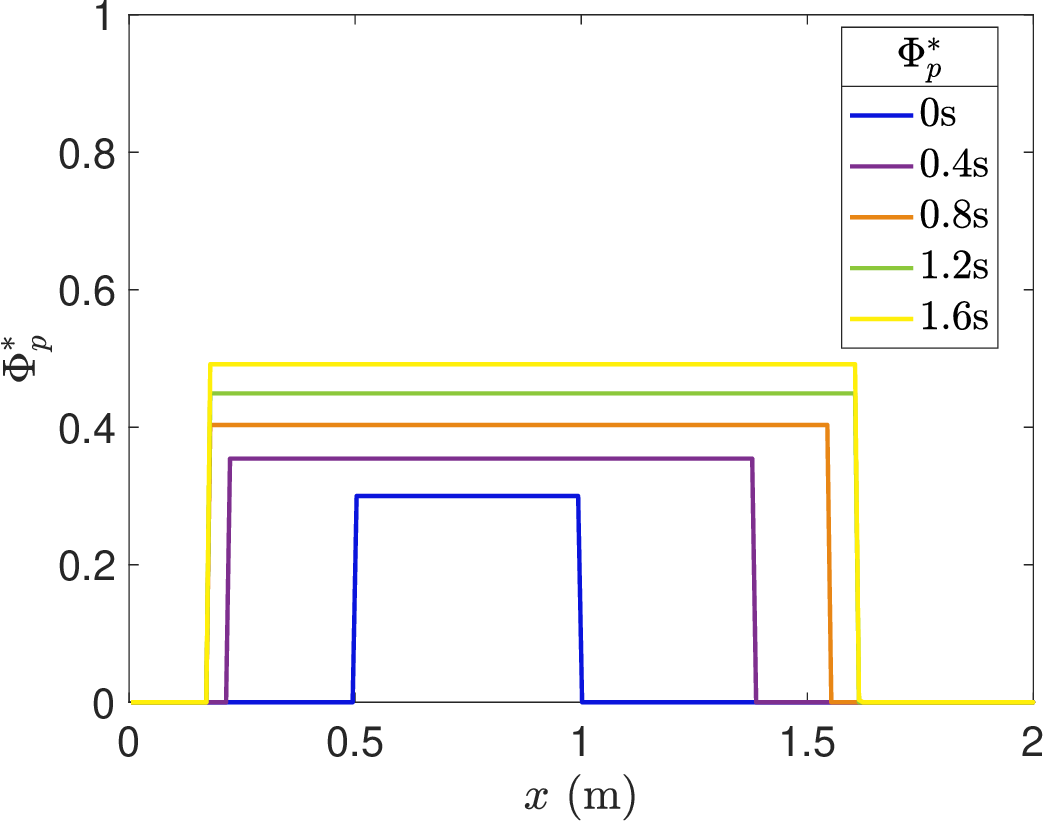}
		\label{fig_dambreak1_model1_phi}
	\end{subfigure}
	\begin{subfigure}[t]{0.49\textwidth}
		\centering      
		\includegraphics[width=\textwidth]{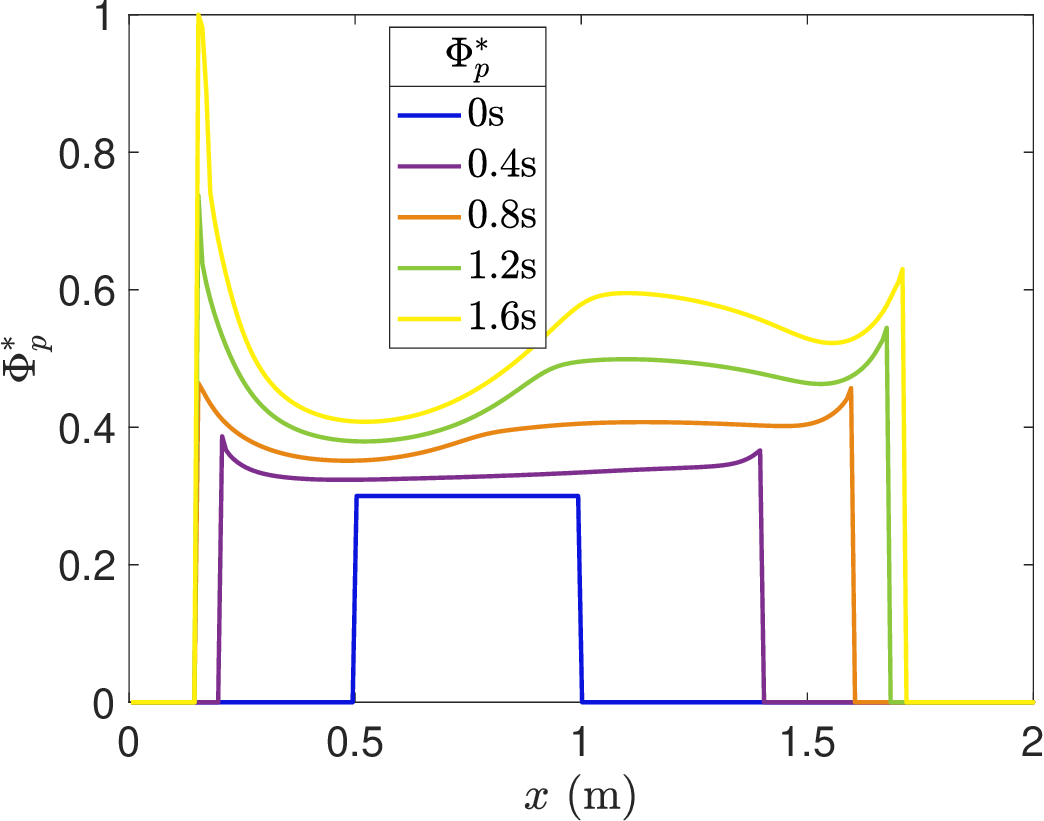}
		\label{fig_dambreak1_model2_phi}
	\end{subfigure}
	\caption{Test $5$, $b>0$, dam break. Left Model $1$, right Model $2$. $\tau_{crys} = 0.5\,$s, $\Phi_p|_{t=0} = 0.3$.}
	\label{fig_dambreak1}
\end{figure}

This numerical test is the same as the previous one, but with dry parts: $h|_{t=0} = 0.2$ for $0.5 \leq x \leq 1$ and vanishes in the rest of the domain. Results for $b=0$ of both models are shown in Figure \ref{fig_dambreak0}. One can observe that the front moves downslope, with a peak of velocity at the front. The fluid velocity decreases with time and vanishes after $1.6$ s. The crystal fraction increases constantly in Model $1$, and increases faster at the wet/dry fronts in Model $2$ where the fluid height is smaller.\\
We then use a wavy bottom given by:
\begin{align*}
    b(x) = 0.1 \sin \frac{6 \pi (x-L/5)}{L},
\end{align*}
and results are shown in Figure \ref{fig_dambreak1}. The fluid, initially located on the side of a wave, spreads out across the domain until stops when the yield stress prevents it to move because of the high crystal fraction.\\ \\
In all these simulations, the fluid velocity vanished at $t=1.6\,$s due to the increase of crystal content or the decrease of temperature, as expected from the models. The solutions behave well at the wet/dry front and there are no numerical oscillations.

\clearpage

\section{Recovering  viscous shear effects in 1D models}\label{app:friction}

The mathematical models describing the fluid flow are defined by 3D equations, usually by the Navier-Stokes equations where the stress tensor takes into account the viscous effects in the fluid. In the incompressible case for a fluid of density $\rho$, they read 
\begin{equation}
\left\{\begin{array}{l}
\nabla\cdot \mathbf{U}=0\\
\rho (\partial_t \mathbf{U}+\nabla(\mathbf{U}\otimes \mathbf{U}))+\nabla p -\nabla\cdot \sigma = \mathbf{E}
\end{array}\right.
\end{equation}
where $\mathbf{U}=(u,v,w)$ is the velocity vector, $p$ is the pressure, $\sigma$ is the viscous stress tensor and $\mathbf{E}=(E_x,E_y,E_z)$ are the external forces (our system in \eqref{eq_2D0} is the 2D version of this system).
The viscous stress tensor is defined in terms of viscosity $\mu$ and deformation tensor $D(\mathbf{U})=\nabla \mathbf{U} + \nabla \mathbf{U}^t$, as $\sigma=\mu D(\mathbf{U})$. The viscosity is constant for a Newtonian fluid and follows different laws for non-Newtonian fluids depending on the rheological properties.\\

Solving such a 3D system being notoriously difficult, a common simplification is to solve associated 1D model by considering the predominant direction of the flow and a unique velocity. The other directions are then neglected. For example, in a flow in direction $x$, we consider $u$ to be the main velocity, thus neglecting the velocities $v$ and $w$. As a consequence, the effects in directions $y$ and $z$ are also neglected. The previous system reduces to the following 1D system 
\begin{equation}\label{system1dapp}
\left\{\begin{array}{l}
\partial_x u =0,\\
\rho (\partial_t u + \partial_x ( u^2)) + \partial_x p = E_x 
\end{array}\right.
\end{equation}
This system \eqref{system1dapp} has no trace of the viscous effects because the first equation yields $\sigma=2\mu \partial_x u=0$. 
In the case, for example, for flows in a pipe, where the velocity follows the direction $x$, such neglect of the friction with the lateral walls implies that one of the largest resistance to motion is also ignored.\\ 

To recover wall friction, the common approach is to replace the viscous stress terms of the momentum equation by a boundary layer approximation, see \cite[equations 6.4-6.5]{dobran:2001}. That is,
$$
-\nabla\cdot\sigma\quad \text{is replaced with}\quad F_S=\frac{1}{R_H}\tau_w,
$$
where $R_H$ is the hydraulic radius and $\tau_w$ is the viscous shear stress,  which are defined as
$$
R_H=\frac{A}{P},\qquad \tau_w = \frac18 f_D \rho u |u|
$$
with $A$ being the flow cross-sectional area, $P$ the wetted perimeter, and $f_D$ the Darcy-Weisbach friction factor (in \cite{dobran:2001} it is written in terms of the Fanning friction factor $f$, which satisfies the relation $f_D=4f$).
The system \eqref{system1dapp} is then reformulated as
\begin{equation}
\left\{\begin{array}{l}
\partial_x u =0,\\
\rho (\partial_t u + \partial_x ( u^2)) + \partial_x p = E_x -F_S
\end{array}\right.
\end{equation}
where $F_S$ represents the added shear stress.
\\
The Darcy-Weisbach friction factor is defined in terms of the Reynolds number as
$$
f_D=\frac{a}{Re}+b.
$$
The first term is dominant in a laminar flow, whereas the second term dominates in a turbulent flow. For a laminar flow in a pipe, the coefficient $a=64$ is obtained from the Poiseuille's law. A few values for other geometries are found in the literature. In \cite{dobran:2001} the value $a=96$ ($a=24=96/4$ in \cite{dobran:2001} corresponding to the Fanning factor) is proposed for fissure or channels. The Reynolds number $Re$ is defined in terms of the hydraulic diameter, $D_H=4R_H$, as
$$
Re=\frac{\rho |u| D_H}{\mu}.
$$
Then, for a laminar flow ($b=0$) the term $F_S$ becomes
$$
F_S=\frac18\frac{1}{R_H} f_D \rho u |u| 
=\frac18\frac{1}{R_H} \frac{a}{Re} \rho u |u|
=\frac{a}{8}\frac{1}{R_H} \frac{\mu}{\rho |u| D_H} \rho u |u|
=\frac{a}{2 D_H^2} \mu u
$$
For a flow in a pipe of radius $r_c$, $a=64$ and 
$
R_H=\frac{\pi r_c^2}{2\pi r_c}=\frac{r_c}{2},
$ and so $D_H=2r_c$ and
$
F_S=\frac{8}{ r_c^2} \mu u.
$
For a flow in an open channel of thickness $h$ and width $w$, $a=96$ and
$$
R_H=\frac{hw}{w+2h}
$$
then
$$
F_S=\frac{3(w+2h)^2}{h^2 w^2} \mu u = \frac{3(1+\frac{2h}{w})^2}{h^2} \mu u.
$$
This term coincides with the friction term in \cite[page 3]{wilson:1993}, which is defined as $F_S^{Wilson}=\frac{K}{D^2}\mu u$. In this term, $K$ is a shape factor and $D$ is the flow thickness. For a semicircular channel, they are introduced as  $K=8$ and $D=r_c$, and for a rectangular channel, $K=3(1+\frac{2h}{w})^2$ and $D=h$. For flows which are very much wider than they are deep, $K=3$ (no lateral friction is considered).\\

Considering volume conservation, an averaged velocity is calculated in \cite[equation 3]{wilson:1993} as $V_m=\frac{D^2 \rho g \sin\theta}{K\mu}$. In \cite{chevrel:2018} for the application in PyFlowgo, the authors introduce also a mean velocity as $V_{mean}=\frac{D^2 \rho g \sin\theta}{n_{shape}\mu}\left(1-\frac32\frac{\tau_c}{\tau_b}+\frac12\frac{\tau_c^3}{\tau_b^3}\right)$, with the shape factor defined as $n_{shape}=3(1+\frac{h}{w})^2$. In the Newtonian case, $\tau_c=0$, this velocity coincides with the one proposed in \cite{wilson:1993}, but there is a difference in the definition of the shape factor, $K=3(1+\frac{2h}{w})^2$ while $n_{shape}=3(1+\frac{h}{w})^2$.  The shape factor $n_{shape}$ implies that the wetted perimeter of the rectangular channel is $2(w+h)$, which corresponds to a closed channel. In a lava flow, it corresponds to assuming flow in a lava tube instead of an open channel. \\

In this work, we consider a depth-integrated model, where the friction effect is kept thanks to the boundary condition at the bottom \eqref{bc_bottom}. 
Nevertheless, lateral friction is not considered. In order to add this friction to our model, we follow the studies presented above. First, we replace $u$ with $h\bar u$  in the expression of $F_S$ and we propose a definition for the viscosity $\mu$ according to the Herschel-Bulkley rheology, with the stress tensor as in \eqref{eq_Herschel-Bulkley}, which is, under deformation: 
$$\sigma = \mu D(\mathbf{U}), \quad\text{with}
		\quad \mu = \left(\eta |D(\mathbf{U})|^{n} + \tau_c\right)\frac{1}{|D(\mathbf{U})|}.
$$
Taking into account the approximation $|D(\mathbf{U})|\sim |\partial_z u|\sim  \frac{|\bar u|}{h}$, we propose to define the viscosity in the Reynolds number above as
$$
\mu\sim \left(\eta \left|\frac{\bar u}{h}\right|^n+\tau_c\right) \frac{h}{|\bar u|}.
$$
Then, the additional shear stress term becomes
$$
\bar F_S=3\left(1+\frac{2h}{w}\right)^2\left(\eta \left|\frac{\bar u}{h}\right|^n+\tau_c\right) \frac{\bar u}{|\bar u|}.
$$ 
The coefficient $s_b$ appearing in the term on the right-hand side of the velocity equation in the system \eqref{eq_syst1D_application} comes from the approximation of the shear rate at the bottom in \eqref{dzu_aprox}, $(\partial_z u)_{|b}\sim s_b \frac{\bar u}{h}$. 
For simplicity, we also apply this approximation for the shear stress against the lateral boundaries. As a result, we obtain an explicit definition of the coefficient $a_f$ that is linked to the shape of the channel: $a_f=3\left(1+\frac{2h}{w}\right)^2$.

\end{document}